\documentclass[10pt]{amsart}
\usepackage[margin=2.5cm]{geometry}
\usepackage{tikz-cd}
\usepackage{graphicx}					
\usepackage{amssymb}
\usepackage{amsmath}
\usepackage{relsize}
\usepackage{subfig, graphicx}
\usepackage{mathrsfs}
\usepackage{amssymb, amsthm, indentfirst}
\usepackage{relsize}
\usepackage{hyperref}
\usepackage{stmaryrd}
\usepackage{lineno}
\usepackage{mathabx}

\numberwithin{equation}{section}
\usepackage{color}
\theoremstyle{plain}
\newtheorem{thm}{Theorem}[section]

\newtheorem{conj}[thm]{Conjecture}
\newtheorem{lemma}[thm]{Lemma}
\newtheorem{prop}[thm]{Proposition}

\theoremstyle{definition}

\newtheorem{rmk}[thm]{Remark}
\def\Gal{\operatorname{Gal}}

\newtheorem*{hypothesis*}{Hypothesis}

\author{SHIH-YU CHEN}
\address{Institute of Mathematics~\\Academia Sinica~\\ 6F, Astronomy-Mathematics Building, No.\,1, Sec.\,4, Roosevelt Road, Taipei 10617, Taiwan, ROC}
\email{sychen0626@gate.sinica.edu.tw}

\def\GL{{\rm{GL}}}
\def\GSp{{\rm GSp}}

\def\Sp{{\rm Sp}}

\def\o{\frak{o}}
\def\c{\frak{c}}
\def\A{{\mathbb A}}
\def\C{{\mathbb C}}
\def\E{{\mathbb E}}
\def\F{{\mathbb F}}
\def\K{{\mathbb K}}
\def\L{{\mathbb L}}
\def\R{{\mathbb R}}
\def\Q{{\mathbb Q}}
\def\Z{{\mathbb Z}}

\def\<{\langle}
\def\>{\rangle}
\def\G{\mathbf{G}}

\def\bp{\begin{pmatrix}}
\def\ep{\end{pmatrix}}

\def\<{\langle}
\def\>{\rangle}

\def\GL{\operatorname{GL}}

\def\GSp{\operatorname{GSp}}

\def\Sp{\operatorname{Sp}}

\def\1{\mathbf{1}}

\def\itPi{\mathit{\Pi}}
\def\itPsi{\mathit{\Psi}}
\def\itSigma{\mathit{\Sigma}}

\makeatletter
\newcommand{\exterior}[1]{\mathop{\mathpalette\exterior@{#1}}}
\newcommand{\exterior@}[2]{%
  \raisebox{\depth}{%
  \fontsize{\sf@size}{0}%
  \m@th
  $\ifx#1\displaystyle\textstyle\else#1\fi\bigwedge$}%
  ^{\mspace{-2mu}#2}%
  \kern-\scriptspace
}
\makeatother

\title{On Deligne's conjecture for symmetric sixth $L$-functions of Hilbert modular forms}
\begin{document}

\begin{abstract}
In this paper, we prove Deligne's conjecture for symmetric sixth $L$-functions of Hilbert modular forms. We extend the result of Morimoto based on a different approach. We define automorphic periods associated to globally generic $C$-algebraic cuspidal automorphic representations of $\GSp_4$ over totally real number fields whose archimedean components are (limits of) discrete series representations. We show that the algebraicity of critical $L$-values for $\GSp_4 \times \GL_2$ can be expressed in terms of these periods. In the case of Kim--Ramakrishnan--Shahidi lifts of $\GL_2$, we establish period relations between the automorphic periods and powers of Petersson norm of Hilbert modular forms. The conjecture for symmetric sixth $L$-functions then follows from these period relations and our previous work on the algebraicity of critical values for the adjoint $L$-functions for $\GSp_4$.
\end{abstract}

\maketitle

\section{Introduction}

In \cite{Deligne1979}, Deligne proposed a conjecture on the algebraicity of values of motivic $L$-functions at critical points. 
For example, attach to a normalized elliptic newform $f$ of weight $\kappa \geq 2$ and nebentypus $\omega$, we can define the symmetric $n$-th power $L$-function $L(s,{\rm Sym}^n(f))$ for each integer $n \geq 1$. In this case, the conjectural Deligne's period $c^\pm({\rm Sym}^n(f))$ was computed in \cite[Proposition 7.7]{Deligne1979}. It is predicted that
\begin{align}\label{E:introduction}
\sigma \left( \frac{L(m,{\rm Sym}^n(f))}{(2\pi\sqrt{-1})^{d_\pm m}\cdot c^\pm({\rm Sym}^n(f))} \right) = \frac{L(m,{\rm Sym}^n({}^\sigma\!f))}{(2\pi\sqrt{-1})^{d_\pm m}\cdot c^\pm({\rm Sym}^n({}^\sigma\!f))}
\end{align}
for all $\sigma \in {\rm Aut}(\C)$ and critical points $m \in \Z$ of $L(s,{\rm Sym}^n(f))$. Here $\pm = (-1)^m$, $d_+=d_- = \tfrac{n+1}{2}$ if $n$ is odd, and $d_+=1+\tfrac{n}{2}$, $d_-=\tfrac{n}{2}$ if $n$ is even.
For $n=1$ and $n=2$, the conjecture follows from the results of Shimura \cite{Shimura1976}, \cite{Shimura1977} and Sturm \cite{Sturm1980}, \cite{Sturm1989}, respectively.
In \cite{GH1993}, Garrett and Harris studied the algebraicity of critical values of triple product $L$-functions in the balanced case. As a consequence, they proved the conjecture for $n=3$ under the assumptions that $\kappa \geq 6$ and $m$ is on the right and away from the center of the critical strip. We complement the result of Garrett and Harris in \cite{Chen2021d} and relaxed the assumption to $\kappa \geq 3$. Based on symmetric power functoriality and various algebraicity results on critical $L$-values in the literature, especially \cite{GL2021}, Morimoto proved the algebraicity of the critical $L$-values for $n=4,6$ in \cite{Morimoto2021}. More precisely, under the assumptions that $\kappa \geq 6$, $\omega=1$, and ${\rm Aut}(\C)$ replaced by ${\rm Aut}(\C/\E)$ for some bi-quadratic extension $\E/\Q$, Morimoto proves that (\ref{E:introduction}) holds up to square. In \cite{Chen2021}, we refine the result of Grobner and Lin, and prove the conjecture for $n=4$ assuming $\kappa \geq 3$.
Actually, we have analogous conjecture for symmetric power $L$-functions associated to normalized Hilbert cusp newforms over a totally real number field $\F$. The conjecture holds for $n=1,2,3,4,6$ by authors mentioned above under similar assumptions.
The goal of this paper is to prove the conjecture for $n=6$. We extend the result of Morimoto based on a different approach. In particular, specializing our main result Theorem \ref{T:main} to $\F=\Q$, we show that (\ref{E:introduction}) holds for $n=6$ assuming $\kappa \geq 6$. We also refer to \cite[\S\,7]{KT2016} for numerical examples in this case.


\subsection{Main result}
Let $\F$ be a totally real number field with $[\F:\Q]=d$.
Let $\itPi$ be a regular algebraic irreducible cuspidal automorphic representation of $\GL_2(\A_\F)$ with central character $\omega_\itPi$. We have $|\omega_\itPi| = |\mbox{ }|_{\A_\F}^{\sf w}$ for some ${\sf w} \in \Z$.
Let $f_\itPi$ be the normalized newform of $\itPi$.
The Petersson norm of $f_\itPi$ is defined by
\begin{align}\label{E:Petersson norm}
\Vert f_\itPi \Vert = \int_{\A_\F^\times\GL_2(\F)\backslash \GL_2(\A_\F)} |f_\itPi(g)|^2|\det(g)|_{\A_\F}^{-{\sf w}}\,dg^{\rm Tam}.
\end{align}
Here $dg^{\rm Tam}$ is the Tamagawa measure on $\A_\F^\times \backslash \GL_2(\A_\F)$.
For each archimedean place $v$ of $\F$, there exists $\kappa_v \geq 2$ such that $\itPi_v$ is isomorphic to the discrete series representation of $\GL_2(\R)$ with minimal weights $\pm \kappa_v$ and has central character ${\rm sgn}^{\sf w}|\mbox{ }|^{\sf w}$.
For a finite order Hecke character $\chi$ of $\A_\F^\times$, let 
\[
L(s,\itPi,{\rm Sym}^6 \otimes \chi)
\]
be the twisted symmetric sixth $L$-function of $\itPi$ by $\chi$. We denote by $L^{(\infty)}(s,\itPi,{\rm Sym}^6 \otimes \chi)$ the $L$-function obtained by excluding the archimedean $L$-factors. A critical point for the $L$-function $L(s,\itPi,{\rm Sym}^6 \otimes \chi)$ is an integer $m$ which is not a pole of the archimedean local factors $L(s,\itPi_v,{\rm Sym}^6 \otimes \chi_v)$ and $L(1-s,\itPi_v^\vee,{\rm Sym}^6 \otimes \chi_v^{-1})$ for all archimedean places $v$ of $\F$.
More precisely, $L(s,\itPi,{\rm Sym}^6 \otimes \chi)$ has no critical points if the signature of $\chi$ is not parallel. Suppose $\chi$ has parallel signature ${\rm sgn}(\chi)\in \{\pm1\}$, then the set of critical points for $L(s,\itPi,{\rm Sym}^6 \otimes \chi)$ is the union of the set of right-half critical points 
\[
\left.\left\{1-3{\sf w} \leq  m \leq \min_{v \mid \infty}\{\kappa_v\}-1-3{\sf w}\,\right\vert\, (-1)^{m+1}={\rm sgn}(\chi)\right\}
\]
and the set of left-half critical points 
\[
\left.\left\{-\min_{v \mid \infty}\{\kappa_v\}+2-3{\sf w} \leq  m \leq -3{\sf w}\,\right\vert\, (-1)^m={\rm sgn}(\chi)\right\}.
\]
The Deligne's periods of $L(s,\itPi,{\rm Sym}^6)$ are given by
\[
c^\pm(\itPi,{\rm Sym}^6) = \begin{cases}
(2\pi\sqrt{-1})^{\sum_{v \mid \infty}6\kappa_v}\cdot G(\omega_\itPi)^{12}\cdot \Vert f_\itPi \Vert^6 & \mbox{ if $\pm = +$},\\
|D_\F|^{1/2}\cdot (2\pi\sqrt{-1})^{\sum_{v \mid \infty}6\kappa_v}\cdot G(\omega_\itPi)^9 \cdot \Vert f_\itPi \Vert^6 & \mbox{ if $\pm = -$}.
\end{cases}
\]
Here $D_\F$ is the discriminant of $\F$ and $G(\omega_\itPi)$ is the Gauss sum of $\omega_\itPi$.
We refer to \cite[Proposition 2.2]{Yoshida1994} for the appearance of $|D_\F|^{1/2}$.
For $\sigma \in {\rm Aut}(\C)$, let ${}^\sigma\!\itPi$ be the unique regular algebraic irreducible cuspidal automorphic representation of $\GL_2(\A_\F)$ such that ${}^\sigma\!\itPi_f$ is the $\sigma$-conjugate of $\itPi_f = \bigotimes_{v \nmid \infty}\itPi_v$.
We have the following conjecture proposed by Deligne \cite[\S\,7]{Deligne1979} on the algebraicity of the critical values of $L(s,\itPi,{\rm Sym}^6 \otimes \chi)$.

\begin{conj}[Deligne]\label{C:DC 2}
Let $\chi$ be a finite order Hecke character of $\A_\F^\times$ and $m-3{\sf w} \in {\rm Crit}(\itPi,{\rm Sym}^6\otimes\chi)$ a critical point. For all $\sigma \in {\rm Aut}(\C)$, we have
\[
\sigma \left( \frac{L^{(\infty)}(m-3{\sf w},\itPi,{\rm Sym}^6\otimes\chi)}{(2\pi\sqrt{-1})^{d_\pm dm}\cdot G(\chi)^{d_\pm}\cdot c^\pm(\itPi,{\rm Sym}^6)}\right) = \frac{L^{(\infty)}(m-3{\sf w},{}^\sigma\!\itPi,{\rm Sym}^6\otimes{}^\sigma\!\chi)}{(2\pi\sqrt{-1})^{d_\pm dm}\cdot G({}^\sigma\!\chi)^{d_\pm}\cdot c^\pm({}^\sigma\!\itPi,{\rm Sym}^6)}.
\]
Here $\pm = (-1)^{m+{\sf w}+1}{\rm sgn}(\chi)$, $d_+=4$, $d_-=3$.
\end{conj}
\begin{rmk}
The conjecture is true when $\itPi$ is of CM-type.
In \cite{Morimoto2021}, Conjecture \ref{C:DC 2} was proved (up to square when $\F=\Q$) by Morimoto under the following assumptions:
\begin{itemize}
\item[(1)] $\F\cap\Q(\zeta_5)=\Q$.
\item[(2)] $\kappa_v \geq 6$ for all archimedean places $v$ of $\F$.
\item[(3)] $\omega_\itPi = |\mbox{ }|_{\A_\F}^{{\sf w}}$.
\item[(4)] ${\rm Aut}(\C)$ replaced by ${\rm Aut}(\C/\E^{\Gal})$ for some CM-extension $\E/\F$, where $\E^{\Gal}$ is the Galois closure of $\E$ in $\C$.
\end{itemize}
The first assumption is imposed in order to apply the result \cite{CT2017} on the functoriality for ${\rm Sym}^6$. For the second assumption, it is required so that we can apply the result \cite{Morimoto2018}. Assumptions (3) and (4) can be lifted if one can refine and generalize the result \cite{GL2021} in the case $\GL_4 \times \GL_3$ over CM-fields.  
\end{rmk}

Following is the main result of this paper.
We extend the result of Morimoto based on a different approach.
Instead of using the result \cite{GL2021}, we define automorphic periods for $\GSp_4$ and establish period relations for Kim--Ramakrishnan--Shahidi lifts for $\GL_2$.
\begin{thm}\label{T:main}
Assume $\F\cap\Q(\zeta_5)=\Q$ and $\kappa_v \geq 6$ for all archimedean places $v$ of $\F$.
\begin{itemize}
\item[(1)] Conjecture \ref{C:DC 2} holds with ${\rm Aut}(\C)$ replaced by ${\rm Aut}(\C/\F^{\Gal})$.
\item[(2)] If we assume further that $\kappa_v = \kappa_w$ for all archimedean places $v,w$ of $\F$, then Conjecture \ref{C:DC 2} holds.
\end{itemize}
\end{thm}

\begin{rmk}
For assertion (1), to descend the equivariance from ${\rm Aut}(\C/\F^{\Gal})$ to ${\rm Aut}(\C)$, it suffices to improve the equivariance property in \cite[Lemma 3.3.1]{GL2016} (cf.\,Remark \ref{R:5.2}).
\end{rmk}

\subsection{An outline of the proof}
There are two key ingredients in the proof:
\begin{itemize}
\item[(1)] The algebraicity of critical values of twisted standard $L$-functions for $\GSp_6(\A_\F)$ (Theorem \ref{T:Liu}).
\item[(2)] The period relations for automorphic periods of the Kim--Ramakrishnan--Shahidi lifts and Petersson norms of Hilbert cusp newforms (Theorem \ref{P:period relation main}).
\end{itemize}
We may assume $\itPi$ is non-CM, that is, $\itPi$ is not an automorphic induction of a Hecke character over some CM-extension of $\F$.
Based on the results of Arthur \cite{Arthur2013}, Clozel--Throne \cite{CT2017}, and Patrikis \cite{Patrikis2019}, in Proposition \ref{P:descent} we show that $\itPi$ transfer weakly to an $C$-algebraic irreducible cuspidal automorphic representation of $\GSp_6(\A_\F)$ with respect to the symmetric sixth power representation of $\GL_2(\C)$. More precisely, we consider the following diagram on the Galois side:
\[
\begin{tikzcd}[row sep=1.8em, column sep=6em]
  L_\F \arrow[r,"\phi_\itPi"] \arrow[ddrr, dashrightarrow]
    & \GL_2(\C) \arrow[r,"{\rm Sym}^6\otimes\,\omega_\itPi^{-3}"] & \GL_7(\C) \\
& & {\rm SO}_7(\C)\arrow[u, hook] \\
& & {\rm GSpin}_7(\C)\arrow[u, two heads]
\end{tikzcd}
\]
Here $L_\F$ and $\phi_\itPi$ are the conjectural Langlands group of $\F$ and the conjectural Langlands parameter associated to $\itPi$, respectively. We twist $\omega_\itPi^{-3}$ so that the image of $({\rm Sym}^6\circ\phi_\itPi)\otimes\omega_\itPi^{-3}$ lies in ${\rm SO}_7(\C)$.
Moreover, the transfer is strong at archimedean places. 
By ingredient (1), when $\F \neq \Q$, we then conclude that Conjecture \ref{C:DC 2} holds for all $\chi$ and $m$ if and only if it holds for some $\chi$ and $m$.
Theorem \ref{T:main}-(2) for $\F=\Q$ follows from the cases for $\F\neq \Q$ and a base change trick (cf.\,\S\,\ref{SS:3.2.2}). We show that Conjecture \ref{C:DC 2} holds for 
\[
\chi = |\mbox{ }|_{\A_\F}^{3{\sf w}}\omega_\itPi^{-3},\quad m=1-3{\sf w}.
\]
For this specific case, we prove the conjecture by ingredient (2).
In \S\,\ref{SS:periods for GSp_4}, for a globally generic $C$-algebraic irreducible cuspidal automorphic representation $\itSigma$ of $\GSp_4(\A_\F)$ whose archimedean components are (limits of) discrete series representations, we attach to an automorphic period $p^I(\itSigma) \in \C^\times$ for each (admissible) subset $I$ of $S_\infty$. Here $S_\infty$ denote the set of archimedean places of $\F$. The periods are obtained by comparing rational structures on $\itSigma_f$ given by the Whittaker model and by the coherent cohomology of certain automorphic vector bundles on the Hilbert--Siegel modular variety associated to $\GSp_{4,\F}$.
By the Serre duality for coherent cohomology and our previous result \cite{Chen2021c}, we prove in Theorem \ref{T:Ichino-Chen} the period relation that
\[
L(1,\itSigma,{\rm Ad})\sim p^I(\itSigma)\cdot p^{S_\infty\smallsetminus I}(\itSigma^\vee).
\]
Here $L(s,\itSigma,{\rm Ad})$ is the adjoint $L$-function of $\itSigma$ and $\alpha \sim \beta$ if $\alpha=\gamma\cdot\beta$ for some $\gamma \in \pi^\Z\cdot\overline{\Q}^\times$.
In particular, we take $\itSigma$ be the Kim--Ramakrishnan--Shahidi lift of $\itPi$ (cf.\,\S\,\ref{SS:KRS lift}). 
In this case, ingredient (2) roughly says that
\begin{align}\label{E:1.3}
p^\varnothing(\itSigma)\sim \Vert f_\itPi \Vert^3,\quad p^{S_\infty}(\itSigma)\sim \Vert f_\itPi \Vert^4.
\end{align}
We have the factorization of $L$-functions:
\[
L(s,\itSigma,{\rm Ad}) = L(s,\itPi,{\rm Sym}^6\otimes\omega_\itPi^{-3})\cdot L(s,\itPi,{\rm Sym}^2\otimes\omega_\itPi^{-1}).
\]
Since $L(1,\itPi,{\rm Sym}^2\otimes\omega_\itPi^{-1})\sim \Vert f_\itPi \Vert$, we thus conclude that
\[
L(1,\itPi,{\rm Sym}^6\otimes\omega_\itPi^{-3})\sim \Vert f_\itPi\Vert^6.
\]
Theorem \ref{T:main} for the critical value $L(1,\itPi,{\rm Sym}^6\otimes\omega_\itPi^{-3})$ then follows by gathering the fudge factors appear in the above relations $\sim$.
We sketch the stretergy to establish the period relation (\ref{E:1.3}). We consider auxiliary algebraic irreducible cuspidal automorphic representations $\itPi_1'$ and $\itPi_2'$ of $\GL_2(\A_\F)$ satisfying certain conditions.\
Based on the cohomological interpretation of Novodvorsky's global zeta integral for Rankin--Selberg $L$-functions for $\GSp_4 \times \GL_2$ (cf.\,Proposition \ref{P:Galois equiv. global}), we prove the following result on algebraicity of the rightmost critical $L$-values:
\begin{align*}
L(m_1,\itSigma \times \itPi_1') \sim p^\varnothing(\itSigma)\cdot \Vert f_{\itPi_1'} \Vert,\quad L(m_2,\itSigma \times \itPi_2') \sim p^{S_\infty}(\itSigma).
\end{align*}
Here $m_1,m_2$ are the rightmost critical points.
Note that the results are special cases of Theorem \ref{T:algebraicity GSp_4 x GL_2}.
On the other hand, when we take $\itSigma$ be the Kim--Ramakrishnan--Shahidi lift of $\itPi$, the algebraicity of these critical values can be expressed in terms of powers of $\Vert f_\itPi \Vert$ using the results \cite{GH1993}, \cite{GL2016}, \cite{Morimoto2018}, \cite{Chen2021d}, \cite{Chen2021}, \cite{Harris2021}, \cite{JST2021}. We then obtain the period relation (\ref{E:1.3}).
The details for the proof of Theorem \ref{P:period relation main} are given in \S\,\ref{S:period relations}.

\subsection{Notation}\label{SS:notation}
Fix a totally real number field $\F$ with $[\F:\Q]=d$.
Let $D_\F$ be the discriminant of $\F$.
Let $\A_\F$ be the ring of adeles of $\F$ and $\A_{\F,f}$ its finite part. 
Let $\hat{\o}_\F$ be the maximal compact subring of $\A_{\F,f}$.
Let $\zeta_\F(s)$ be the completed Dedekind zeta function of $\F$.
Let $\psi_\Q=\bigotimes_v\psi_{v}$ be the standard additive character of $\Q\backslash \A_\Q$ defined so that
\begin{align*}
\psi_{p}(x) & = e^{-2\pi \sqrt{-1}\,x} \mbox{ for }x \in \Z[p^{-1}],\\
\psi_{\infty}(x) & = e^{2\pi \sqrt{-1}\,x} \mbox{ for }x \in \R.
\end{align*}
Let $\psi_\F = \psi_\Q\circ{\rm tr}_{\F/\Q}$ and call it the standard additive character of $\F \backslash \A_\F$.
Let $S_\infty$ be the set of archimedean places of $\F$. 
For $v \in S_\infty$, let $\iota_v$ be the real embedding of $\F$ associated to $v$ and identify $\F_v$ with $\R$ via $\iota_v$.
Let $v$ be a place of $\F$. If $v$ is a finite place, let $\frak{o}_{v}$, $\varpi_v$, and $q_v$ be the maximal compact subring of $\F_v$, a generator of the maximal ideal of $\frak{o}_{v}$, and the cardinality of $\frak{o}_{v} / \varpi_v\frak{o}_{v}$. Let $|\mbox{ }|_v = |\mbox{ }|$ be the absolute value on $\F_v$ normalized so that $|\varpi_v|_v = q_v^{-1}$. If $v \in S_\infty$, let $|\mbox{ }|_v$ be the be the ordinary absolute value on $\R$.
Let $|\mbox{ }|_{\A_\F} = \prod_v|\mbox{ }|_v$ be the adelic norm on $\A_\F$.
Let $\Gamma_\R(s)$ and $\Gamma_\C(s)$ be the archimedean gamma functions defined by
\[
\Gamma_\R(s) = \pi^{-s/2}\Gamma(\tfrac{s}{2}),\quad \Gamma_\C(s) = 2(2\pi)^{-s}\Gamma(s).
\]

Let $\chi$ be an algebraic Hecke character of $\A_\F^\times$.
The signature of $\chi$ at $v \in S_\infty$ is the value $\chi_v(-1) \in \{\pm1\}$. The signature ${\rm sgn}(\chi)$ of $\chi$ is the sequence of signs $(\chi_v(-1))_{v \in S_\infty}$.
We say $\chi$ has parallel signature if it has the same signature at all real places.
The Gauss sum $G(\chi)$ of $\chi$ is defined by
\[
G(\chi) = |D_\F|^{-1/2}\prod_{v \nmid \infty}\varepsilon(0,\chi_v,\psi_v),
\]
where $\psi_\F = \bigotimes_v\psi_v$ and $\varepsilon(s,\chi_v,\psi_v)$ is the $\varepsilon$-factor of $\chi_v$ with respect to $\psi_v$ defined in \cite{Tate1979}.
For $\sigma \in {\rm Aut}(\C)$, let ${}^\sigma\!\chi$ be the unique algebraic Hecke character  of $\A_\F^\times$ such that ${}^\sigma\!\chi(x) = \sigma(\chi(x))$ for $x \in \A_{\F,f}^\times$. Note that ${\rm sgn}(\chi) = {\rm sgn}({}^\sigma\!\chi)$.
It is easy to verify that 
\begin{align}\label{E:Galois Gauss sum}
\begin{split}
\sigma(G(\chi)) = {}^\sigma\!\chi(u_\sigma)G({}^\sigma\!\chi),\quad
\sigma\left(\frac{G(\chi\chi')}{G(\chi)G(\chi')}\right) = \frac{G({}^\sigma\!\chi{}^\sigma\!\chi')}{G({}^\sigma\!\chi)G({}^\sigma\!\chi')}
\end{split}
\end{align}
for algebraic Hecke characters $\chi,\chi'$ of $\A_\F^\times$,
where $u_\sigma \in \prod_p \Z_p^\times$ is the unique element such that $\sigma(\psi_{\Q}(x)) = \psi_{\Q}(u_\sigma x)$ for $x \in \A_{\Q,f}$.

Let $\sigma \in {\rm Aut}(\C)$. Define the $
\sigma$-linear action on $\C(X)$, which is the field of formal Laurent series in variable $X$ over $\C$, as follows:
\[
{}^\sigma\!P(X) = \sum_{n \gg -\infty}^\infty\sigma(a_n) X^n
\]
for $P(X) = \sum_{n \gg -\infty}^\infty a_n X^n \in \C(X)$.
For a complex representation $\itPi$ of a group $G$ on the space ${V}_\itPi$ of $\itPi$, let ${}^\sigma\!\itPi$ be the representation of $G$ defined
\begin{align}\label{E:sigma action}
{}^\sigma\!\itPi(g) = t \circ \itPi(g) \circ t^{-1},
\end{align}
where $t:{V}_\itPi \rightarrow {V}_\itPi$ is a $\sigma$-linear isomorphism. Note that the isomorphism class of ${}^\sigma\!\itPi$ is independent of the choice of $t$. We call ${}^\sigma\!\itPi$ the $\sigma$-conjugate of $\itPi$. 

\section{Automorphic periods for $\GL_2$ and $\GSp_4$}\label{S:automorphic period}

\subsection{Automorphic vector bundles on Shimura varieties}

Let $(G,X)$ be a Shimura datum, that is, $G$ is a connected reductive algebraic group over $\Q$ and $X$ is a $G(\R)$-conjugacy class of homomorphisms $h: \mathbb{S} \rightarrow G_\R$ satisfying conditions (1.1.1)-(1.1.3) in \cite{Harris1985}. Here $\mathbb{S} = {\rm Res}_{\C / \R}\mathbb{G}_m$ is the Deligne torus. We have the associated Shimura variety 
\[
{\rm Sh}(G,X) = \varprojlim_{\mathcal{K}} {\rm Sh}_{\mathcal K}(G,X) = \varprojlim_{\mathcal{K}}G(\Q)\backslash X\times G(\A_{\Q,f}) / \mathcal{K},
\]
where $\mathcal{K}$ runs through neat open compact subgroups of $G(\A_{\Q,f})$.
It is a pro-algebraic variety over $\C$ with continuous $G(\A_{\Q,f})$-action and admits canonical model over the reflex field $E(G,X)$ of $(G,X)$.

Fix $h \in X$. Let $K_h$ be the stabilizer of $h$ in $G(\R)$. The Hodge decomposition induced by ${\rm Ad}\circ h$ on the complexified  Lie algebra $\frak{g}_\C$ of $G(\R)$ is given by
\begin{align*}\label{E:Hodge}
\frak{g}_\C = \frak{g}_\C^{(-1,1)} \oplus \frak{g}_\C^{(0,0)} \oplus \frak{g}_\C^{(1,-1)}.
\end{align*}
Here
\[
\frak{g}_\C^{(p,q)} = \left.\{X \in \frak{g}_\C \, \right\vert \, h(z)^{-1}Xh(z) = z^{-p}\overline{z}^{-q}X \mbox{ for $z \in \C$}\}.
\]
We write $\frak{p}_h^{\pm} = \frak{g}_\C^{\pm(-1,1)}$. Note that $\frak{p}_h^+$ and $\frak{p}_h^-$ are the holomorphic and anti-holomorphic tangent spaces of $X$ at $h$, respectively, and $\frak{g}_\C^{(0,0)}$ is the complexified Lie algebra $\frak{k}_{h,\C}$ of $K_h$.
Let $\frak{P}_h$ be the subalgebra of $\frak{g}_\C$ defined by
\begin{align*}\label{E:subalgebra}
\frak{P}_h = \frak{k}_{h,\C} \oplus \frak{p}_h^-,
\end{align*}
and $P_h$ the parabolic subgroup of $G(\C)$ with Lie algebras $\frak{P}_h$. We write
\[
\check{X}_h = G(\C) / P_h.
\]
Note that the flag variety $\check{X}_h$ has a natural structure over $E(G,X)$. The inclusion $G(\R)\subset G(\C)$ induces an embedding
\[
\beta_h: X \simeq G(\R) / K_h \hookrightarrow \check{X}_h.
\]
For a $G_\C$-vector bundle $\mathcal{V}$ on $\check{X}_h$, we denote by $\beta_h^*(\mathcal{V})$ the pullback $G(\R)$-bundle on $X$. 
We say $\mathcal{V}$ is motivic if the following condition is satisfied:
\begin{align*}\label{E:motivic}
\mbox{The action of $G_\C$ on $\mathcal{V}$ factors through $(G/Z_s)_\C$}.
\end{align*}
Here $Z_s$ is the largest subtorus of the center $Z_G$ of $G$ that is split over $\R$ but that has no subtorus split over $\Q$.
For motivic $\mathcal{V}$, the automorphic vector bundle $[\mathcal{V}]$ on ${\rm Sh}(G,X)$ is defined by
\[
[\mathcal{V}] = \varprojlim_{\mathcal{K}} G(\Q)\backslash \,\beta_h^*(\mathcal{V}) \times G(\A_{\Q,f}) / \mathcal{K},
\]
where $\mathcal{K}$ runs through neat open compact subgroups of $G(\A_{\Q,f})$. We have the following result due to Harris \cite[Theorem 3.3]{Harris1985} and Milne \cite[Theorem 5.1]{Milne1990}.

\begin{thm}
The functor $\mathcal{V} \mapsto [\mathcal{V}]$, from motivic $G_\C$-vector bundles on $\check{X}_h$ to $G(\A_{\Q,f})$-vector bundles on ${\rm Sh}(G,X)$, is rational over $E(G,X)$.
\end{thm}
\noindent
Note that we have an equivalence of categories (cf.\,\cite[Remark 1.2]{IP2018})
\[
\mbox{$G_\C$-vector bundles on $\check{X}_h$} \longleftrightarrow \mbox{finite-dimensional algebraic representations of $P_h$ over $\C$}.
\]

\subsection{Coherent cohomology groups on Shimura varieties}

We keep the notation of the previous section. 
Let $\mathcal{A}(G(\A_\Q))$, $\mathcal{A}_{(2)}(G(\A_\Q))$, and $\mathcal{A}_0(G(\A_\Q))$ be the spaces of $K_h$-finite automorphic forms, essentially square-integrable automorphic forms, and cusp forms, respectively, on $G(\A_\Q)$.
Let $(\rho,V)$ be an irreducible algebraic representation of $K_h$. Then $\rho$ extends to an algebraic representation of $P_h$ so that the action factors through the reductive quotient $K_{h,\C}$. We denote by $\mathcal{V}_\rho$ the corresponding $G_\C$-vector bundle on $\check{X}_h$.
We say $(\rho,V)$ is motivic if $\mathcal{V}_\rho$ is motivic.
We have the $(\frak{P}_h,K_h)$-modules 
\[
\mathcal{A}(G(\A_\Q))\otimes V,\quad \mathcal{A}_{(2)}(G(\A_\Q))\otimes V,\quad\mathcal{A}_0(G(\A_\Q))\otimes V,
\]
where the action of $\frak{P}_h$ on $V$ factors through $\frak{k}_\C$.
Consider the complexes with respect to the Lie algebra differential operator (cf.\,\cite[Chapter I]{BW2000}):
\begin{align}\label{E:complexes}
\begin{split}
C_{\rho}^\bullet &= \left( \mathcal{A}(G(\A_\Q)) \otimes \exterior{\bullet} \frak{p}_h^+ \otimes V\right)^{K_h},\\
C_{(2),\rho}^\bullet &= \left( \mathcal{A}_{(2)}(G(\A_\Q)) \otimes \exterior{\bullet} \frak{p}_h^+ \otimes V\right)^{K_h},\\
C_{{\rm cusp},\,\rho}^\bullet &= \left( \mathcal{A}_{0}(G(\A_\Q)) \otimes \exterior{\bullet} \,\frak{p}_h^+ \otimes  V\right)^{K_h}.
\end{split}
\end{align}
The corresponding $q$-th $(\frak{P}_h,K_h)$-cohomology groups of the above complexes are denoted respectively by 
\[
H^q(\frak{P}_h,K_h; \mathcal{A}(G(\A_\Q))\otimes V),\quad H^q(\frak{P}_h,K_h;\mathcal{A}_{(2)}(G(\A_\Q)) \otimes V),\quad H^q(\frak{P}_h,K_h;\mathcal{A}_0(G(\A_\Q)) \otimes V).
\]
Note that $G(\A_{\Q,f})$ acts on the above complexes by right translation. This in turn defines $G(\A_{\Q,f})$-module structures on the cohomology groups.
Suppose $V$ is motivic, by the result of Harris \cite[Corollary 3.4]{Harris1990} and Jun \cite[Theorem 6.7]{Su2018}, the cohomology group $H^q(\frak{P}_h,K_h;\mathcal{A}(G(\A_\Q)) \otimes V)$ is canonically isomorphic as $G(\A_{\Q,f})$-module to the $q$-th coherent cohomology group of the canonical extension $[\mathcal{V}_\sigma]^{\rm can}$ of the automorphic vector bundle $[\mathcal{V}_\rho]$ to a good toroidal compactification of ${\rm Sh}(G,X)$. 
Hereafter we identify the coherent cohomology groups with the relative Lie algebra $(\frak{P}_h,K_h)$-cohomology groups and write
\[
H^q([\mathcal{V}_\rho]^{\rm can}) = H^q(\frak{P}_h,K_h;\mathcal{A}(G(\A_\Q)) \otimes V).
\]
We also write 
\[
H^q_{(2)}([\mathcal{V}_\rho]) = H^q(\frak{P}_h,K_h; \mathcal{A}_{(2)}(G(\A_\Q))\otimes V),\quad H^q_{\rm cusp}([\mathcal{V}_\rho]) = H^q(\frak{P}_h,K_h; \mathcal{A}_0(G(\A_\Q))\otimes V).
\]
The natural inclusions $\mathcal{A}_0(G(\A_\Q)) \subset \mathcal{A}_{(2)}(G(\A_\Q)) \subset \mathcal{A}(G(\A_\Q))$ induce the $G(\A_{\Q,f})$-module homomorphisms
\begin{align}\label{E:diag 1}
H^q_{\rm cusp}([\mathcal{V}_\rho]) \longrightarrow H^q_{(2)}([\mathcal{V}_\rho]) \longrightarrow H^q([\mathcal{V}_\rho]^{\rm can}).
\end{align}
Let $[\mathcal{V}_\rho]^{\rm sub}$ be the subcanonical extension of $[\mathcal{V}_\rho]$ on the good toroidal compactification over which $[\mathcal{V}_\rho]^{\rm can}$ is defined.
We denote by $H_!^q([\mathcal{V}_\rho])$ the image of the homomorphism $H^q([\mathcal{V}_\rho]^{\rm sub}) \rightarrow H^q([\mathcal{V}_\rho]^{\rm can})$ induced by the exact sequence \cite[(2.2.4)]{Harris1990}.
We call $H_!^q([\mathcal{V}_\rho])$ the $q$-th interior cohomology group of $[\mathcal{V}_\rho]$.
In the following theorem, we recall some results of Harris \cite{Harris1985}, \cite{Harris1990} and Milne \cite{Milne1983} on the coherent cohomology groups.
Fix a Cartan subalgebra $\frak{h}$ of $\frak{k}_h$ and a positive system of $(\frak{g}_\C,\frak{h}_\C)$ such that the set of non-compact positive roots gives the root spaces in $\frak{p}_h^+$.
Denote by $\Lambda_\rho \in \frak{h}_\C^*$ the corresponding highest weight of $\rho$. 
Let $\mathcal{A}_{(2)}(G(\A_\Q),\rho)$ (resp.\,$\mathcal{A}_0(G(\A_\Q),\rho)$) be the space consisting of essentially square-integrable automorphic forms (resp.\,cusp forms) on $G(\A_\Q)$ which are eigenfunctions of the Casimir operator of $\frak{g}_\C$ with eigenvalue
\[
\<\Lambda_\rho+ \delta, \,\Lambda_\rho+\delta\>_{\frak{h}^*} - \<\delta,\,\delta\>_{\frak{h}^*}.
\]
Here $\<\,,\,\>_{\frak{h}^*}$ is the Killing form on $\frak{h}^* = {\rm Hom}_\R(\frak{h},\R)$ and $\delta$ is the half-sum of positive roots.
 
\begin{thm}\label{T:Harris}
Let $(\rho,V)$ be an irreducible motivic algebraic representation of $K_h$ defined over some finite extension $\Q(\rho)$ of $E(G,X)$.
\begin{itemize}
\item[(1)] For $\sigma \in {\rm Aut}(\C / E(G,X))$, conjugation by $\sigma$ induces natural $\sigma$-linear $G(\A_{\Q,f})$-module isomorphisms
\[
T_\sigma: H^q([\mathcal{V}_\rho]^{\rm sub}) \longrightarrow H^q([\mathcal{V}_{{}^\sigma\!\rho}]^{\rm sub}),\quad T_\sigma: H^q([\mathcal{V}_\rho]^{\rm can}) \longrightarrow H^q([\mathcal{V}_{{}^\sigma\!\rho}]^{\rm can}),
\]
and such that the diagram 
\[
\begin{tikzcd}
H^q([\mathcal{V}_\rho]^{\rm sub}) \arrow[r, "T_\sigma"] \arrow[d] & H^q([\mathcal{V}_{{}^\sigma\!\rho}]^{\rm sub}) \arrow[d]\\
H^q([\mathcal{V}_\rho]^{\rm can}) \arrow[r, "T_\sigma"]  & H^q([\mathcal{V}_{{}^\sigma\!\rho}]^{\rm can})
\end{tikzcd}
\]
is commutative. 
Moreover, $H^q([\mathcal{V}_\rho]^{\rm sub})$ and $H^q([\mathcal{V}_\rho]^{\rm can})$ are admissible $G(\A_{\Q,f})$-modules and have canonical rational structures over $\Q(\rho)$ given by taking the Galois invariants with respect to $T_\sigma$ for $\sigma \in {\rm Aut}(\C / \Q(\rho))$.
\item[(2)] We have
\begin{align*}
H_{(2)}^q([\mathcal{V}_\rho]) &= \left(\mathcal{A}_{(2)}(G(\A_\Q),\rho)\otimes \exterior{q} \,\frak{p}_h^+\otimes V \right)^{K_h},\\
H_{\rm cusp}^q([\mathcal{V}_\rho]) &= \left(\mathcal{A}_{0}(G(\A_\Q),\rho)\otimes \exterior{q} \,\frak{p}_h^+\otimes V \right)^{K_h}.
\end{align*}
\item[(3)] The homomorphism $H_{\rm cusp}^q([\mathcal{V}_\rho])\rightarrow H^q([\mathcal{V}_\rho]^{\rm can})$ in (\ref{E:diag 1}) is injective and its image is contained in $H_!^q([\mathcal{V}_\rho])$.
\item[(4)] The interior cohomology group $H_!^q([\mathcal{V}_\rho])$ is contained in the image of the homomorphism $H_{{(2)}}^q([\mathcal{V}_\rho])\rightarrow H^q([\mathcal{V}_\rho]^{\rm can})$ in (\ref{E:diag 1}). In particular, $H_!^q([\mathcal{V}_\rho])$ is a semisimple $G(\A_{\Q,f})$-module.
\end{itemize}
\end{thm}

We identify $H_{\rm cusp}^q([\mathcal{V}_\rho])$ with a $G(\A_{\Q,f})$-submodule of $H_!^q([\mathcal{V}_\rho])$ via the injection in Theorem \ref{T:Harris}-(3).
By Theorem \ref{T:Harris}-(1), we have
\begin{align}\label{E:Galois equiv. class}
T_\sigma(H_!^q([\mathcal{V}_\rho])) = H_!^q([\mathcal{V}_{{}^\sigma\!\rho}])
\end{align}
for all $\sigma \in {\rm Aut}(\C/E(G,X))$.

\subsection{Hilbert--Siegel modular varieties}\label{SS:Hilbert--Siegel}
For $n\geq 1$, Let $\GSp_{2n}$ be the symplectic similitude group defined by
\[
 \GSp_{2n} =
 \left\{ g \in \GL_{2n} \, \left| \, g
 \begin{pmatrix}
  0      & \1_n \\
  - \1_n & 0
 \end{pmatrix}
 {}^t \! g = 
  \nu(g) \begin{pmatrix}
  0      & \1_n \\
  - \1_n & 0
 \end{pmatrix}, \,\nu(g) \in \GL_1
 \right. \right\}.
\]
Let $\GSp_{2n}^+(\R)$ be the closed subgroup of $\GSp_{2n}(\R)$ consisting of elements with positive similitude.
Let $K_n = \R_+\times{\rm U}(n)$ and regard it as a closed subgroup of $\GSp_{2n}(\R)$ by the homomorphism
\[
(a, A+\sqrt{-1}\,B)\longmapsto a{\bf 1}_{2n}\cdot\bp A&B\\-B&A\ep.
\]
We write $\frak{g}_n\subset \frak{gl}_{2n}$ for the Lie algebra of $\GSp_{2n}(\R) \subset \GL_{2n}(\R)$.  Define $\frak{k}_n \subset \frak{g}_n$ and $\frak{p}_n^\pm,\frak{P}_n \subset \frak{g}_{n,\C}$ by
\begin{align*}
\frak{k}_n = {\rm Lie}(K_n), \quad
\frak{p}_{n}^\pm = \left\{  \left.\bp -\sqrt{-1}\,A & \pm A \\ \pm A & \sqrt{-1}\,A\ep \, \right\vert \, A \in {\rm Sym}_n(\C)\right\},\quad \frak{P}_n = \frak{k}_{n,\C}\oplus\frak{p}_n^-.
\end{align*}
We will identify $\frak{p}_n^\pm$ with ${\rm Sym}_n(\C)$ by the map
\[
\bp -\sqrt{-1}\,A & \pm A \\ \pm A & \sqrt{-1}\,A\ep \longmapsto A.
\]

Let $(G_n,X_n)$ be the Shimura datum defined by 
\[
G_n = {\rm Res}_{\F/\Q}\GSp_{2n,\F},
\]
and $X_n$ is the $G_n(\R)$-conjugacy class containing the morphism $h_n : \mathbb{S} \rightarrow G_{n,\R}$ with
\[
h_n(x+\sqrt{-1}\,y) = \left(\bp x{\bf 1}_n & y{\bf 1}_n \\ -y{\bf 1}_n & x{\bf 1}_n\ep ,\cdots,\bp x{\bf 1}_n & y{\bf 1}_n \\ -y{\bf 1}_n & x{\bf 1}_n\ep\right)
\]
on $\R$-points.
The associated Shimura variety is called the Hilbert--Siegel modular variety. Note that the reflex field $E(G_n,X_n)$ is equal to $\Q$. 
Under the identification of $\F_v$ with $\R$ for each $v \in S_\infty$, we have
\[
K_{h_n} = K_n^{S_\infty},\quad \frak{k}_{h_n} = \frak{k}_n^{S_\infty},\quad \frak{p}_{h_n}^\pm = (\frak{p}_{n}^\pm)^{S_\infty},\quad \frak{P}_{h_n} = \frak{P}_{n}^{S_\infty}.
\]
Any motivic irreducible algebraic representation of $K_n^{S_\infty}$ is of the form:
\[
(\rho,V_\rho) = \left( \bigotimes_{v \in S_\infty}\rho_v,\bigotimes_{v \in S_\infty}V_{\rho_v}\right)
\]
for some irreducible algebraic representation $(\rho_v,V_{\rho_v})$ of $K_n$ such that $\rho_v \vert_{\R_+} \simeq \rho_w\vert_{\R_+}$ for all $v, w \in S_\infty$.
For $\sigma \in {\rm Aut}(\C)$, let $({}^\sigma\!\rho,V_{{}^\sigma\!\rho})$ be the motivic irreducible algebraic representation defined by 
\[
({}^\sigma\!\rho,V_{{}^\sigma\!\rho}) = \left( \bigotimes_{v \in S_\infty}\rho_{\sigma^{-1}\circ v},\bigotimes_{v \in S_\infty}V_{\rho_{\sigma^{-1}\circ v}}\right).
\]

\subsection{Automorphic periods for $\GL_2$}

In this section, we recall the automorphic periods of $C$-algebraic irreducible cuspidal automorphic representations of $\GL_2(\A_\F)$ defined by Harris \cite{Harris1989I} through the coherent cohomology for automorphic line bundles on ${\rm Sh}(G_1,X_1)$.

\subsubsection{(Limits of) discrete series representations of $\GL_2(\R)$}

Let $(\kappa;\,{\sf w}) \in \Z\times\Z$ such that $\kappa \equiv {\sf w}\,({\rm mod}\,2)$.
Let $(\rho_{(\kappa;\,{\sf w})}, V_{(\kappa;\,{\sf w})})$ be the algebraic character of $K_1$ defined by $V_{(\kappa;\,{\sf w})} = \C$ and 
\[
\rho_{(\kappa;\,{\sf w})}(au)\cdot z = a^{\sf w}u^\kappa\cdot z 
\]
for $a \in \R^\times$ and $u \in {\rm U}(1)$. Note that $\rho_{(\kappa;\,{\sf w})}^\vee = \rho_{(-\kappa;\,-{\sf w})}$.
Assume further that
\[
\kappa \geq 1.
\]
Let $D_{(\pm\kappa;\,{\sf w})}$ be the (limit of) discrete series representation of $\GL_2^+(\R)$ with minimal $K_1$-type $\rho_{(\pm\kappa;\,{\sf w})}$.
Let $\itPi_{(\kappa;{\sf w})}$ be the (limit of) discrete series representation of $\GL_2(\R)$ defined by
\[
\itPi_{(\kappa;\,{\sf w})} = D_{(\kappa;\,{\sf w})} \oplus D_{(-\kappa;\,{\sf w})},
\]
with 
\[
{\rm diag}(-1,1)\cdot({\bf v}_1,{\bf v}_2) =  ({\bf v}_2,{\bf v}_1).
\]
Conversely, up to central twists, any (limit of) discrete series representation of $\GL_2(\R)$ is obtained in this way. 
It is well-known that $\itPi_{(\kappa;\,{\sf w})}$ is generic. Let $\psi_\R$ be the standard additive character of $\R$, that is, $\psi_\R(x)=e^{2\pi\sqrt{-1}\,x}$.
We denote by $\mathcal{W}(\itPi_{(\kappa;\,{\sf w})},\psi_\R)$ the space of Whittaker functions of $\itPi_{(\kappa;\,{\sf w})}$ with respect to $\psi_\R$.
Let $W_{(\kappa;\,{\sf w})}^+ \in \mathcal{W}(\itPi_{(\kappa;\,{\sf w})},\psi_\R)$ be the Whittaker function of weight $\kappa$ normalized so that
\begin{align}\label{E:GL_2 Whittaker}
W_{(\kappa;\,{\sf w})}^+ ({\rm diag}(a,1)) = |a|^{(\kappa+{\sf w})/2}e^{-2\pi|a|}\cdot\mathbb{I}_{\R_+}(a).
\end{align}

\subsubsection{Rational structure via the Whittaker model}

Let $(\itPi,V_\itPi)$ be an $C$-algebraic irreducible cuspidal automorphic representation of $\GL_2(\A_\F)$ with central character $\omega_\itPi$. Here we follow \cite[Definition 5.11]{BG2014} for the notion of $C$-algebraicity, it is refer to algebraic automorphic representation in \cite{Clozel1990}. 
We assume further that the following condition is satisfied:
\begin{align}\label{E:discrete condition 1}
\mbox{$\itPi_v$ is a (limit of) discrete series representation for all $v\in S_\infty$.}
\end{align}
Then there exists $(\underline{\kappa};\,{\sf w}) \in \Z^{S_\infty}\times\Z$ with $\underline{\kappa} = (\kappa_v)_{v\in S_\infty}$ satisfying the following conditions:
\begin{itemize}
\item $|\omega_\itPi| = |\mbox{ }|_{\A_\F}^{\sf w}$.
\item $\kappa_v \geq 1$ and $\kappa_v \equiv {\sf w}\,({\rm mod}\,2)$ for all $v \in S_\infty$.
\item $\itPi_v = \itPi_{(\kappa_v;\,{\sf w})}$ for all $v \in S_\infty$.
\end{itemize}
We say $\itPi$ is regular if $\kappa_v \geq 2$ for all $v \in S_\infty$.
We call $(\underline{\kappa};\,{\sf w})$ the weight of $\itPi$. 
Let $\psi$ be a non-trivial additive character of $\F\backslash\A_\F$. 
Let $\itPi_f = \bigotimes_{v \nmid \infty}\itPi_v$ be the finite part of $\itPi$ and $\mathcal{W}(\itPi_f,\psi_f)$ the space of Whittaker functions of $\itPi_f$ with respect to $\psi_f$.
For $\varphi \in V_\itPi$, let $W_{\varphi,\psi}$ be the Whittaker function of $\varphi$ with respect to $\psi$ defined by
\[
W_{\varphi,\psi}(g) = \int_{\F \backslash \A_\F}\varphi\left(\bp 1 & x \\ 0 & 1\ep g \right)\overline{\psi(x)}\,dx^{\rm Tam}.
\]
Here $dx^{\rm Tam}$ is the Tamagawa measure on $\A_\F$.
Let $V_\itPi^+$ be the space of cusp forms in $V_\itPi$ of weight $\underline{\kappa}$, that is, $\varphi \in V_\itPi^+$ if and only if 
\[
\varphi(gu) = \prod_{v \in S_\infty}u_v^{\kappa_v}\cdot \varphi(g)
\]
for all $ u = (u_v)_{v \in S_\infty} \in {\rm U}(1)^{S_\infty}$ and $g \in \GL_2(\A_\F)$. 
Take $\psi=\psi_\F$ to be the standard additive character. For $\varphi \in V_\itPi^+$, let $W_\varphi^{(\infty)} \in \mathcal{W}(\itPi_f,\psi_f)$ be the unique Whittaker function so that
\[
W_{\varphi,\psi} = \prod_{v \in S_\infty}W_{(\kappa_v;\,{\sf w})}^+\cdot W_\varphi^{(\infty)}.
\]
Then the map $\varphi \mapsto W_\varphi^{(\infty)}$ defines a $\GL_2(\A_{\F,f})$-equivariant isomorphism from $V_\itPi^+$ to $\mathcal{W}(\itPi_f,\psi_f)$.
For $\sigma \in {\rm Aut}(\C)$, let 
\[
{}^\sigma\!\itPi = {}^\sigma\!\itPi_\infty \otimes {}^\sigma\!\itPi_f
\]
be the irreducible admissible representation of $\GL_2(\A_\F)$ defined so that the $v$-component of ${}^\sigma\!\itPi_\infty$ is isomorphic to $\itPi_{\sigma^{-1}\circ v}$ for each $v \in S_\infty$.
It is known that ${}^\sigma\!\itPi$ is cuspidal automorphic. Moreover, it is clear that ${}^\sigma\!\itPi$ is $C$-algebraic of weight $({}^\sigma\!\underline{\kappa};\,{\sf w})$, where ${}^\sigma\!\underline{\kappa} = (\kappa_{\sigma^{-1}\circ v})_{v \in S_\infty}$.
Let $\Q(\underline{\kappa})$ and $\Q(\itPi)$ be the fixed fields of $\left\{\sigma \in {\rm Aut}(\C)\,\vert\,{}^\sigma\!\underline{\kappa} = \underline{\kappa}\right\}$ and $\left\{\sigma \in {\rm Aut}(\C)\,\vert\,{}^\sigma\!\itPi_f = \itPi_f\right\}$, respectively.
Note that $\Q(\underline{\kappa})\subset \Q(\itPi)$ by the strong multiplicity one theorem for $\GL_2$.
Let 
\[
t_\sigma : \mathcal{W}(\itPi_f,\psi_f) \longrightarrow \mathcal{W}({}^\sigma\!\itPi_f,\psi_f)
\]
be the $\sigma$-linear $\GL_2(\A_{\F,f})$-equivariant isomorphism  defined by
\[
t_\sigma W(g) = \sigma \left( W({\rm diga}(u_\sigma^{-1},1)g)\right).
\]
Here $u_\sigma \in \widehat{\Z}^\times \subset \hat{\o}_\F^\times$ is the unique element such that $\sigma(\psi(x)) = \psi(u_\sigma x)$ for all $x \in \A_{\F,f}$.
Let 
\begin{align*}
V_\itPi^+ \longrightarrow V_{{}^\sigma\!\itPi}^+,\quad \varphi\longmapsto {}^\sigma\!\varphi
\end{align*}
be the $\sigma$-linear $\GL_2(\A_{\F,f})$-equivariant isomorphism defined such that 
\[
W_{{}^\sigma\!\varphi}^{(\infty)} = t_{\sigma}W_{\varphi}^{(\infty)}.
\]
We thus obtain a $\Q(\itPi)$-rational structure $(V_\itPi^+)^{{\rm Aut}(\C/\Q(\itPi))}$ on $V_\itPi^+$ given by taking the Galois invariants:
\begin{align}\label{E:rational structure 1}
(V_\itPi^+)^{{\rm Aut}(\C/\Q(\itPi))} = \left.\left\{\varphi\in V_\itPi^+\,\right\vert\,{}^\sigma\!\varphi=\varphi\mbox{ for $\sigma \in {\rm Aut}(\C/\Q(\itPi))$}\right\}.
\end{align}

\subsubsection{Rational structure via the coherent cohomology}
Now we recall the rational structures on $V_\itPi^+$ given by the coherent cohomology of automorphic lines bundles on ${\rm Sh}(G_1,X_1)$.
For $I\subset S_\infty$, define $\underline{\kappa}(I) = (\kappa_v(I))_{v \in S_\infty} \in \Z^{S_\infty}$ by
\[
\kappa_v(I) = \begin{cases}
\kappa_v-2 & \mbox{ if $v \in I$},\\
-\kappa_v & \mbox{ if $v \notin I$}.
\end{cases}
\]
By definition, we have ${}^\sigma\!(\underline{\kappa}(I)) = {}^\sigma\!\underline{\kappa}({}^\sigma\!I)$ for $\sigma\in{\rm Aut}(\C)$.
Let $\Q(I)$ be the fixed field of $\left\{\sigma \in {\rm Aut}(\C)\,\vert\,{}^\sigma\!I=I\right\}$.
It is clear that $\underline{\kappa}(I)$ is invariant by ${\rm Aut}(\C/\Q(\underline{\kappa})\Q(I))$.
We say $I$ is $admissible$ with respect to $\underline{\kappa}$ if $I$ is uniquely determined by ${}^\sharp I$ and $\underline{\kappa}(I)$, that is, if $J\subset S_\infty$ with ${}^\sharp J = {}^\sharp I$ and $\underline{\kappa}(J) = \underline{\kappa}(I)$, then $J=I$.
For instance, $I$ is admissible when $I \in \{\varnothing, S_\infty\}$ or $\kappa_v \geq 2$ for all $v \in I$.
Let 
\[
[\mathcal{V}_{(\underline{\kappa}(I);\,-{\sf w})}]
\]
denote the automorphic line bundle on ${\rm Sh}(G_1,X_1)$ defined by the motivic algebraic representation 
\[
(\rho_{(\underline{\kappa}(I);\,-{\sf w})},V_{(\underline{\kappa}(I);\,-{\sf w})}) = \left( \bigotimes_{v \in S_\infty}\rho_{({\kappa}_v(I);\,-{\sf w})}, \bigotimes_{v \in S_\infty}V_{({\kappa}_v(I);\,-{\sf w})}\right)
\]
of $K_1^{S_\infty}$. Let $H_{\rm cusp}^q([\mathcal{V}_{(\underline{\kappa}(I);\,-{\sf w})}])[\itPi_f]$ be the $\itPi_f$-isotypic component of $H_{\rm cusp}^q([\mathcal{V}_{(\underline{\kappa}(I);\,-{\sf w})}])$.
We have the following result proved in \cite[Lemmas 1.4.3 and 1.4.5]{Harris1989I} (see also \cite[Lemma 2.5 and Remark 2.7]{Chen2021e}).

\begin{lemma}\label{L:rational structure 1}
Let $I \subset S_\infty$.
\begin{itemize}
\item[(1)] For all $\sigma \in {\rm Aut}(\C)$ and $q\geq 0$, we have
\[
T_\sigma(H_{\rm cusp}^q([\mathcal{V}_{(\underline{\kappa}(I);\,-{\sf w})}])[\itPi_f]) = H_{{\rm cusp}}^q([\mathcal{V}_{({}^\sigma\!\underline{\kappa}({}^\sigma\!I);\,-{\sf w})}])[{}^\sigma\!\itPi_f].
\]
In particular, we have a $\Q(\itPi)\Q(I)$-rational structure on $H_{\rm cusp}^q([\mathcal{V}_{(\underline{\kappa}(I);\,-{\sf w})}])[\itPi_f]$ given by taking the Galois invariants:
\begin{align*}
&H_{\rm cusp}^q([\mathcal{V}_{(\underline{\kappa}(I);\,-{\sf w})}])[\itPi_f]^{{\rm Aut}(\C/\Q(\itPi)\Q(I))}\\
& = \left.\left\{c\in H_{\rm cusp}^q([\mathcal{V}_{(\underline{\kappa}(I);\,-{\sf w})}])[\itPi_f]\,\right\vert\,T_\sigma c=c\mbox{ for $\sigma \in {\rm Aut}(\C/\Q(\itPi)\Q(I))$}\right\}.
\end{align*}
\item[(2)]
Assume $I$ is admissible, then $H_{\rm cusp}^{{}^\sharp I}([\mathcal{V}_{(\underline{\kappa}(I);\,-{\sf w})}])[\itPi_f] \simeq \itPi_f$.
In this case, we have a $\GL_2(\A_{\F,f})$-equivariant isomorphism
\[
V_\itPi^+ \longrightarrow H_{\rm cusp}^{{}^\sharp I}([\mathcal{V}_{(\underline{\kappa}(I);\,-{\sf w})}])[\itPi_f],\quad \varphi \longmapsto [\varphi]_I
\]
defined by
\begin{align*}
[\varphi]_I = \varphi^I \otimes \text{\LARGE$\wedge$}_{v \in I}X_{+,v} \otimes \bigotimes_{v \in S_\infty} 1 \in \mathcal{A}_0(\GL_2(\A_\F)) \otimes \exterior{{}^\sharp I}(\frak{p}_1^+)^{S_\infty} \otimes V_{(\underline{\kappa}(I);\,-{\sf w})}.
\end{align*}
Here $\varphi^I$ is defined by
\[
\varphi^I(g) = \varphi\left(g\cdot\prod_{v \in I}{\rm diag}(-1,1)\right),
\]
and $X_+ = \bp\sqrt{-1} &  -1 \\ -1 & -\sqrt{-1}\ep \in \frak{p}_1^+$ and $X_{+,v} \in \frak{p}_1^{S_\infty}$ is defined so that its $v$-component is $X_+$ and zero otherwise. We also fix an ordering of the wedge $\text{\LARGE$\wedge$}_{v \in I}X_{+,v} \in \exterior{{}^\sharp I}(\frak{p}_1
^+)^{S_\infty}$ once and for all.
\end{itemize}
\end{lemma}

\subsubsection{Automorphic periods and period relations}

By comparing the rational structures in (\ref{E:rational structure 1}) and Lemma \ref{L:rational structure 1}-(1), we 
have the following lemma/definition for the automorphic periods of $\itPi$.

\begin{lemma}
Let $I\subset S_\infty$ be admissible with respect to $\underline{\kappa}$. There exists a sequence of non-zero complex numbers $\left(p^{{}^\sigma\!I}({}^\sigma\!\itPi)\right)_{\sigma \in {\rm Aut}(\C)}$ such that
\[
T_\sigma\left(\frac{[\varphi]_I}{p^I(\itPi)}\right) = \frac{[{}^\sigma\!\varphi]_{{}^\sigma\!I}}{p^{{}^\sigma\!I}({}^\sigma\!\itPi)}
\]
for all $\sigma \in {\rm Aut}(\C)$ and $\varphi \in V_\itPi^+$. Here $T_\sigma : H_!^{{}^\sharp I}([\mathcal{V}_{(\underline{\kappa}(I);\,-{\sf w})}]) \rightarrow H_!^{{}^\sharp I}([\mathcal{V}_{({}^\sigma\!\underline{\kappa}({}^\sigma\!I);\,-{\sf w})}])$ is the $\sigma$-linear isomorphism in (\ref{E:Galois equiv. class}).
\end{lemma}

For the automorphic periods, we have the following period relations proved in \cite[Propositions 1.3.3 and 1.5.6]{Harris1989I}

\begin{lemma}\label{L:period relation 1}
Let $f_\itPi \in V_\itPi^+$ be the normalized newform of $\itPi$ and $\Vert f_\itPi \Vert$ the Petersson norm of $f_\itPi$ defined as in (\ref{E:Petersson norm}). 
\begin{itemize}
\item[(1)] For $\sigma \in {\rm Aut}(\C)$, we have
\[
\sigma \left( \frac{p^\varnothing(\itPi)}{(2\pi\sqrt{-1})^{-\sum_{v \in S_\infty}(\kappa_v+{\sf w})/2}} \right) = \frac{p^\varnothing({}^\sigma\!\itPi)}{(2\pi\sqrt{-1})^{-\sum_{v \in S_\infty}(\kappa_v+{\sf w})/2}}.
\]
\item[(2)] Let $I\subset S_\infty$ be admissible with respect to $\underline{\kappa}$. For $\sigma \in {\rm Aut}(\C)$, we have
\[
\sigma \left( \frac{\Vert f_\itPi \Vert}{p^I(\itPi)\cdot p^{S_\infty \smallsetminus I}(\itPi^\vee)} \right) = \frac{\Vert f_{{}^\sigma\!\itPi} \Vert}{p^{{}^\sigma\!I}({}^\sigma\!\itPi)\cdot p^{S_\infty \smallsetminus {{}^\sigma\!I}}({}^\sigma\!\itPi^\vee)}.
\]
\end{itemize}
\end{lemma}

\subsection{Holomorphic Eisenstein series on $\GL_2$}

In this section, we recall the algebraicity of holomorphic Eisenstein series on $\GL_2(\A_\F)$.
Let $\chi$ be an algebraic Hecke character of $\A_\F^\times$. 
Let ${\sf w} \in \Z$ be the integer such that $|\chi| = |\mbox{ }|^{{\sf w}}_{\A_\F}$.
We assume $\chi$ has parallel signature with ${\rm sgn}(\chi) = (-1)^{\sf w}$.
For $s \in \C$, let
\[
I(\chi,s) = {\rm Ind}_{B_2(\A_\F)}^{\GL_2(\A_\F)}(|\mbox{ }|_{\A_\F}^{s-1/2}\boxtimes \chi |\mbox{ }|_{\A_\F}^{-s+1/2})
\]
be the induced representation consisting of smooth right $(K_1^{S_\infty}\times\GL_2(\A_{\F,f}))$-finite functions $f : \GL_2(\A_\F)\rightarrow \C$ such that
\[
f\left(\bp a &*\\0&d \ep g\right) = \chi(d)\left|\frac{a}{d}\right|_{\A_\F}^s\cdot f(g).
\]
Here $B_2$ is the standard Borel subgroup of $\GL_2$ consisting of upper triangular matrices.
We denote by $\rho$ the right translation action of $K_1^{S_\infty}\times\GL_2(\A_{\F,f})$ on $I(\chi,s)$.
A function
\[
\C \times \GL_2(\A_\F) \longrightarrow \C,\quad (s,g)\longmapsto f^{(s)}(g)
\]
is called a holomorphic section of $I(\chi,s)$ if it satisfies the following conditions:
\begin{itemize}
\item For each $s \in \C$, the function $g \mapsto f^{(s)}(g)$ belongs to $I(\chi ,s)$.
\item For each $g \in \GL_2(\A_\F)$, the function $s \mapsto f^{(s)}(g)$ is holomorphic.
\item $f^{(s)}$ is right $(K_1^{S_\infty}\times\GL_2(\A_{\F,f}))$-finite.
\end{itemize}
A function $f^{(s)}$ on $\C \times \GL_2(\A_\F)$ is called a meromorphic section of $I(\chi,s)$ if there exists a non-zero entire function $\beta$ such that $\beta(s)f^{(s)}$ is a holomorphic section.
For a place $v$ of $\F$, we define $I(\chi_v,s)$ and the notion of holomorphic and meromorphic sections in a similar way.
\begin{itemize}
\item When $v$ is finite, a meromorphic section $f_v^{(s)}$ of $I(\chi_v,s)$ is called a rational section if the map $s\mapsto f_v^{(s)}(g)$ is a rational function in $q_v^{-s}$ for all $g \in \GL_2(\F_v)$. For a rational section $f_v^{(s)}$ and $\sigma \in {\rm Aut}(\C)$, let ${}^\sigma\!f_v^{(s)}$ be the rational section of $I({}^\sigma\!\chi_v,s)$ defined by
\[
{}^\sigma\!f_v^{(s)}\left(\bp a &*\\0&d \ep k\right) = {}^\sigma\!\chi_v(d)\left|\frac{a}{d}\right|_{\A_\F}^s\cdot {}^\sigma\!(f_v^{(s)}(k))
\]
for $a,d \in \F_v^\times$ and $k \in \GL_2(\o_v)$.
\item When $v$ is finite and $\chi_v$ is unramified, let $f_{v,\circ}^{(s)}$ be the $\GL_2(\o_v)$-invariant rational section of $I(\chi_v,s)$ normalized so that 
\[
f_{v,\circ}^{(s)}(1) = L(2s,\chi_v^{-1}).
\]
\item When $v \in S_\infty$ and $\kappa \geq 1$ with $(-1)^\kappa = {\rm sgn}(\chi)$, let $f_{v,\kappa}^{(s)}$ be the holomorphic section of $I(\chi_v,s)$ of weight $\kappa$ normalized so that $f_{v,\kappa}^{(s)}(1)=1$.
\end{itemize}
Let $f^{(s)}$ be a holomorphic section of $I(\chi,s)$. We define the associated Eisenstein series $E(f^{(s)})$ on $\GL_2(\A_\F)$ by the absolutely convergent series 
\[
E(g,f^{(s)}) = \sum_{\gamma \in B_2(\F)\backslash\GL_2(\A_\F)}f^{(s)}(\gamma g)
\]
for ${\rm Re}(s)>1+\tfrac{\sf w}{2}$, and by meromorphic continuation otherwise.
We have the following result on the algebraicity of holomorphic Eisenstein series. When $\kappa>2$, the Eisenstein series converges absolutely at $s=\tfrac{\kappa+{\sf w}}{2}$. In this case, the result is a special case of the result of Harris \cite{Harris1984}.
Based on explicit computation of Fourier coefficients, Shimura \cite{Shimura1978} proved the algebraicity for any $\kappa \geq 1$, except for $\F=\Q$, $\kappa=2$, and $\chi=|\mbox{ }|_{\A_\Q}^{\sf w}$.

\begin{prop}[Harris, Shimura]\label{P:Eisenstein series}
Let $f^{(s)} = \bigotimes_{v \nmid \infty}f_v^{(s)}$ be a meromorphic section of $I(\chi_f,s)$ and $\kappa \geq 1$ with $(-1)^\kappa = {\rm sgn}(\chi)$.
Assume the following conditions are satisfied:
\begin{itemize}
\item $f^{(s)}$ is holomorphic for ${\rm Re}(s)>\tfrac{\sf w}{2}$.
\item $f_v^{(s)}$ is a rational section for all $v \nmid \infty$.
\item $f_v^{(s)} = f_{v,\circ}^{(s)}$ for almost all $v$.
\item If $\F=\Q$, then $\kappa \neq 2$ or $\chi \neq |\mbox{ }|_{\A_\Q}^{\sf w}$.
\end{itemize}
Then the following assertions hold:
\begin{itemize}
\item[(1)] The Eisenstein series $E\left(\bigotimes_{v \in S_\infty} f_{v,\kappa}^{(s)}\otimes f^{(s)}\right)$ is holomorphic at $s=\tfrac{\kappa+{\sf w}}{2}$.
\item[(2)] The automorphic form $E^{[\kappa]}(f^{(s)}) = E\left(\bigotimes_{v \in S_\infty} f_{v,\kappa}^{(s)}\otimes f^{(s)}\right)\vert_{s=\tfrac{\kappa+{\sf w}}{2}}$ defines a global section in the coherent cohomology $H^0([\mathcal{V}_{(-\underline{\kappa};\,-{\sf w})}]^{\rm can})$ given by
\[
[E^{[\kappa]}(f^{(s)})] = E^{[\kappa]}(f^{(s)}) \otimes \bigotimes_{v \in S_\infty}1 \in \mathcal{A}(\GL_2(\A_\F))\otimes V_{(-\underline{\kappa};\,-{\sf w})}.
\]
\item[(3)] For $\sigma \in {\rm Aut}(\C)$, we have
\[
T_\sigma \left( \frac{[E^{[\kappa]}(f^{(s)})]}{|D_\F|^{1/2}\cdot (2\pi\sqrt{-1})^{d(\kappa-{\sf w})/2}\cdot G(\chi)^{-1}}\right) = \frac{[E^{[\kappa]}({}^\sigma\!f^{(s)})]}{|D_\F|^{1/2}\cdot (2\pi\sqrt{-1})^{d(\kappa-{\sf w})/2}\cdot G({}^\sigma\!\chi)^{-1}}.
\]
\end{itemize}
\end{prop}

\subsection{Automorphic periods for $\GSp_4$}\label{SS:periods for GSp_4}

In this section, we define automorphic periods of globally generic $C$-algebraic irreducible cuspidal automorphic representations of $\GSp_4(\A_\F)$ through the coherent cohomology for automorphic vector bundles on ${\rm Sh}(G_2,X_2)$. 
Let $U$ be the maximal unipotent subgroup of $\GSp_4$ defined by
\[
U = \left\{\bp 1 & * & *&* \\0&1&*&*\\0&0&1&0\\0&0&*&1 \ep\in\GSp_4\right\}.
\] 
For a non-trivial additive character $\psi$ of $\F \backslash \A$, let $\psi_U$ be the associated additive character of $U(\F) \backslash U(\A)$ defined by
\[
\psi_U\left(\bp 1 & x & *&* \\0&1&*&y\\0&0&1&0\\0&0&-x&1 \ep\right)= \psi(-x-y).
\]
Similar notation applies to additive character of $\F_v$.

\subsubsection{Algebraic representations of $K_2$}

Let $(\underline{\lambda};\,{\sf u}) \in \Z^2 \times \Z$ with $\underline{\lambda} = (\lambda_1,\lambda_2)$ such that 
\[
\lambda_1 \geq \lambda_2,\quad \lambda_1+\lambda_2 \equiv{\sf u}\,({\rm mod}\,2).
\]
Let $(\rho_{(\underline{\lambda};\,{\sf u})},V_{(\underline{\lambda};\,{\sf u})})$ be the irreducible algebraic representation of $K_2$ defined as follows: $V_{(\underline{\lambda};\,{\sf u})}$ is the space of homogeneous polynomials over $\C$ of degree $\lambda_1-\lambda_2$ in variables $x$ and $y$. The action is given by
\[
(\rho_{(\underline{\lambda};\,{\sf u})}(au)\cdot P)(x,y) = a^{\sf u}(\det u)^{\lambda_2}\cdot P((x,y)u)
\]
for $a \in \R^\times$, $u \in {\rm U}(2)$, and $P \in V_{(\underline{\lambda};\,{\sf u})}$. Put $\underline{\lambda}^\vee = (-\lambda_2,-\lambda_1)$. Note that $\rho_{(\underline{\lambda};\,{\sf u})}^\vee = \rho_{(\underline{\lambda}^\vee;\,-{\sf u})}$. Let $(V_{(\underline{\lambda};\,{\sf u})})_\Q$ be the $\Q$-structure on $V_{(\underline{\lambda};\,{\sf u})}$ consisting of polynomials over $\Q$.
For example, we have an isomorphism 
\[
\frak{p}_2^+ \rightarrow V_{(2,0;\,0)}
\]
given by
\begin{align}\label{E:p_2}
\bp 1 & 0 \\ 0 & 0\ep \longmapsto x^2,\quad \bp 0&1 \\ 1 & 0\ep\longmapsto 2xy,\quad \bp 0&0\\0&1\ep\longmapsto y^2.
\end{align}
Here we identify $\frak{p}_2^+$ with ${\rm Sym}_2(\C)$ as in \S\,\ref{SS:Hilbert--Siegel}.
In this way we also fix a $\Q$-rational structure $(\frak{p}_2^+)_\Q$ on $\frak{p}_2^+$.

We write $\rho_{\underline{\lambda}} = \rho_{(\underline{\lambda};\,{\sf u})} \vert_{{\rm U}(2)}$ and $V_{\underline{\lambda}} = V_{(\underline{\lambda};\,{\sf u})}$ when we consider only the action of ${\rm U}(2)$.
Let $c_{\underline{\lambda}} : V_{\underline{\lambda}} \rightarrow V_{\underline{\lambda}^\vee}$ be the ${\rm U}(2)$-conjugate-equivariant $\C$-linear isomorphism normalized so that
\begin{align*}
c_{\underline{\lambda}}(x^{\lambda_1-\lambda_2}) = y^{\lambda_1-\lambda_2}.
\end{align*}
Let $\<\,,\,\>_{\underline{\lambda}} : V_{\underline{\lambda}} \times V_{\underline{\lambda}^\vee} \rightarrow \C$ be the ${\rm U}(2)$-equivariant bilinear pairing normalized so that
\[
\<x^{\lambda_1-\lambda_2},\,y^{\lambda_1-\lambda_2}\>_{\underline{\lambda}}=1.
\]
Then it is easy to see that
\begin{align}\label{E:U(2) relation}
c_{\underline{\lambda}}(x^{\lambda_1-\lambda_2-i}y^i) = (-1)^i\cdot x^iy^{\lambda_1-\lambda_2-i},\quad
\<x^{\lambda_1-\lambda_2-i}y^i,\,x^iy^{\lambda_1-\lambda_2-i}\>_{\underline{\lambda}} = (-1)^i{\lambda_1-\lambda_2 \choose i}^{-1}
\end{align}
for $0 \leq i \leq \lambda_1-\lambda_2$.
Assume further that $\lambda_1-\lambda_2 \geq 2$.
We fix ${\rm U}(2)$-equivariant embeddings 
\begin{align}\label{E:U(2) embeddings}
\xi_{\underline{\lambda}}^+ : V_{\underline{\lambda}}\longrightarrow\exterior{2}\frak{p}_2^+ \otimes V_{\underline{\lambda}-(3,1)},\quad\xi_{\underline{\lambda}}^- : V_{\underline{\lambda}^\vee}\longrightarrow\frak{p}_2^+ \otimes V_{\underline{\lambda}^\vee-(2,0)}.
\end{align}
The existence of $\xi_{\underline{\lambda}}^\pm$ is a direct consequence of the Clebsch--Gordan formula for ${\rm U}(2)$ and they are unique up to scalar multiples.
We normalize them so that
\[
\xi_{\underline{\lambda}}^+(x^{\lambda_1-\lambda_2}) = \bp 1 & 0 \\ 0 & 0\ep\wedge\bp 0 & 1 \\ 1 & 0\ep \otimes x^{\lambda_1-\lambda_2-2},\quad \xi_{\underline{\lambda}}^-(x^{\lambda_1-\lambda_2}) = \bp 1 & 0 \\ 0 & 0\ep\otimes x^{\lambda_1-\lambda_2-2}.
\]
In particular, $\xi_{\underline{\lambda}}^\pm$ preserve the $\Q$-rational structures.
The following lemma will be used in the proof of Proposition \ref{P:Galois equiv. global}.
\begin{lemma}\label{L:auxiliary}
Assume $\lambda_1-\lambda_2 \geq 2$. Let $0 \leq i \leq \lambda_1-\lambda_2-2$. Then 
\[
\bp 1 & 0 \\ 0 & 0\ep\wedge\bp 0 & 0 \\ 0 & 1\ep \otimes x^{\lambda_1-\lambda_2-2-i}y^i,\quad \bp 1 & 0 \\ 0 & 0\ep \otimes x^{\lambda_1-\lambda_2-2-i}y^i
\]
appear in the decompositions of $\xi_{\underline{\lambda}}^+(x^{\lambda_1-\lambda_2-i-1}y^{i+1})$ and $\xi_{\underline{\lambda}}^-(x^{\lambda_1-\lambda_2-i}y^i)$ into $({\rm U}(1) \times {\rm U}(1))$-eigenvectors.
\end{lemma}

\begin{proof}
Let $N_- \in \frak{k}_{2,\C}$ defined by
\[
N_- = \frac{1}{2}\bp 0 & 1 & 0 &\sqrt{-1} \\ -1&0&\sqrt{-1}&0 \\ 0&-\sqrt{-1}&0&1 \\ -\sqrt{-1}&0&-1&0\ep.
\]
For $0 \leq i \leq \lambda_1-\lambda_2$, we have
\[
N_-\cdot x^{\lambda_1-\lambda_2-i}y^i = -(\lambda_1-\lambda_2-i)\cdot x^{\lambda_1-\lambda_2-i-1}y^{i+1}.
\]
Note that
\[
X\cdot({\bf v}\otimes{\bf w}) = X\cdot {\bf v}\otimes{\bf w}+{\bf v}\otimes X\cdot{\bf w}
\]
for all $X \in \frak{k}_{2,\C}$.
The assertion follows from induction on $0 \leq i \leq \lambda_1-\lambda_2-2$ based on these relations.
We leave the details to the readers.
\end{proof}

\subsubsection{(Limits of) discrete series representations of $\GSp_4(\R)$}

We regard ${\rm U}(1) \times {\rm U}(1)$ as the maximal torus of ${\rm U}(2)$ consisting of diagonal matrices.
Let $\frak{h}$ be the Lie algebra of $\R_+ \times ({\rm U}(1) \times {\rm U}(1))$. Then $\frak{h}$ is a Cartan subalgebra of both $\frak{k}_2$ and $\frak{g}_2$.
Note that 
\[
\frak{h} = \R\cdot\bp 0&0&1&0\\0&0&0&0\\-1&0&0&0\\0&0&0&0 \ep \oplus \R\cdot\bp 0&0&0&0\\0&0&0&1\\0&0&0&0\\0&-1&0&0 \ep \oplus \R\cdot{\bf 1}_4
\]
We identify $\frak{h}_\C^*$ with $\C^3$ by the map
\[
a\cdot\bp 0&0&1&0\\0&0&0&0\\-1&0&0&0\\0&0&0&0 \ep^* \oplus b\cdot\bp 0&0&0&0\\0&0&0&1\\0&0&0&0\\0&-1&0&0 \ep^* \oplus c\cdot{\bf 1}_4^* \longmapsto (\sqrt{-1}\,a,\sqrt{-1}\,b;c).
\]
Here the script $*$ refers to dual basis.
Let $\Delta^+$ and $\Delta_c^+$ be the positive systems of $(\frak{g}_{2,\C},\frak{h}_\C)$ and $(\frak{k}_{2,\C},\frak{h}_\C)$, respectively, given by
\begin{align}\label{E:positive system}
\Delta^+ = \left\{ (1,-1;0),(2,0;0),(0,2;0),(1,1;0)\right\},\quad
\Delta_c^+ = \left\{ (1,-1;0)\right\}.
\end{align}
Note that $\Delta^+ \smallsetminus \Delta_c^+$ corresponds to the root spaces in $\frak{p}_2^+$.
Let $(\underline{\lambda};\,{\sf u}) \in \Z^2 \times \Z$ with $\underline{\lambda} = (\lambda_1,\lambda_2)$ such that 
\[
1-\lambda_1 \leq \lambda_2 \leq 0,\quad \underline{\lambda} \neq (1,0),\quad \lambda_1+\lambda_2 \equiv{\sf u}\,({\rm mod}\,2).
\]
With respect to the choice of positive systems in (\ref{E:positive system}), let $D_{(\underline{\lambda};\,{\sf u})}$ and $D_{(\underline{\lambda}^\vee;\,{\sf u})}$ be the generic (limit of) discrete series representations of $\GSp_4^+(\R)$ with minimal $K_2$-type $\rho_{(\underline{\lambda};\,{\sf u})}$ and $\rho_{(\underline{\lambda}^\vee;\,{\sf u})}$, respectively, defined as in \cite[XII, \S\,7]{Knapp1986}. When $\lambda_1+\lambda_2 \geq 2$, let $D_{(\lambda_1,-\lambda_2+2;\,{\sf u})}$ and $D_{(\lambda_2-2,-\lambda_1;\,{\sf u})}$ be the holomorphic and anti-holomorphic (limit of) discrete series representations of $\GSp_4^+(\R)$ with minimal $K_2$-type $\rho_{(\lambda_1,-\lambda_2+2;\,{\sf u})}$ and $\rho_{(\lambda_2-2,-\lambda_1;\,{\sf u})}$, respectively.
We also write
\begin{align}\label{E:discrete series}
D_{(\underline{\lambda};\,{\sf u})}^{(0)} = D_{(\lambda_1,-\lambda_2+2;\,{\sf u})},\quad D_{(\underline{\lambda};\,{\sf u})}^{(1)} = D_{(\underline{\lambda};\,{\sf u})},\quad 
D_{(\underline{\lambda};\,{\sf u})}^{(2)} =
D_{(\underline{\lambda}^\vee;\,{\sf u})},\quad
D_{(\underline{\lambda};\,{\sf u})}^{(3)} =
D_{(\lambda_2-2,-\lambda_1;\,{\sf u})}.
\end{align}
Let $\itSigma_{(\underline{\lambda};\,{\sf u})}^{\rm gen}$ and $\itSigma_{(\underline{\lambda};\,{\sf u})}^{\rm hol}$ be the generic and holomorphic (limit of) discrete series representations of $\GSp_4(\R)$, respectively, defined by
\[
\itSigma_{(\underline{\lambda};\,{\sf u})}^{\rm gen} = D_{(\underline{\lambda};\,{\sf u})}^{(1)}\oplus D_{(\underline{\lambda};\,{\sf u})}^{(2)},\quad \itSigma_{(\underline{\lambda};\,{\sf u})}^{\rm hol} = D_{(\underline{\lambda};\,{\sf u})}^{(0)}\oplus D_{(\underline{\lambda};\,{\sf u})}^{(3)},
\]
with
\[
{\rm diag}(-1,-1,1,1)\cdot({\bf v}_1,{\bf v}_2) = ({\bf v}_2,{\bf v}_1). 
\]
Then the following set is an $L$-packet of $\GSp_4(\R)$:
\begin{align}\label{E:L-packet}
L_{(\underline{\lambda};\,{\sf u})} = \begin{cases}
\left\{ \itSigma_{(\underline{\lambda};\,{\sf u})}^{\rm gen},\itSigma_{(\underline{\lambda};\,{\sf u})}^{\rm hol}\right\} & \mbox{ if $\lambda_1+\lambda_2 \geq 2$},\\
\left\{\itSigma_{(\underline{\lambda};\,{\sf u})}^{\rm gen}\right\} & \mbox{ if $\lambda_1+\lambda_2=1$}.
\end{cases}
\end{align}
Conversely, up to similitude twists, any $L$-packet of $\GSp_4(\R)$ that contains a (limit of) discrete series representation is obtained in this way.
Note that $\itSigma_{(\underline{\lambda};\,{\sf u})}^{\rm gen}$ is a limit of discrete series representation if and only if $\lambda_1+\lambda_2=1$ or $\lambda_2=0$.
We have the following explicit formula for Whittaker functions on $\GSp_4(\R)$ due to Moriyama \cite[Proposition 7]{Moriyama2004}.
Let $\psi_\R$ be the standard additive character of $\R$ and denote by $\mathcal{W}(\itSigma_{(\underline{\lambda};\,{\sf u})}^{\rm gen},\psi_{\R,U})$ the space of Whittaker functions of $\itSigma_{(\underline{\lambda};\,{\sf u})}^{\rm gen}$ with respect to $\psi_{\R,U}$.

\begin{thm}[Moriyama]\label{T:Moriyama}
There exists a unique $K_2$-equivariant homomorphism
\[
V_{(\underline{\lambda}^\vee;\,{\sf u})} \longrightarrow \mathcal{W}(\itSigma_{(\underline{\lambda};\,{\sf u})}^{\rm gen},\psi_{\R,U}),\quad x^iy^{\lambda_{1}-\lambda_{2}-i} \longmapsto W_{((\underline{\lambda};\,{\sf u}),\,i)}^+
\] 
such that
\begin{align*}
&W_{((\underline{\lambda};\,{\sf u}),\,i)}^+({\rm diag}(a_1a_2,a_1,a_2^{-1},1))\\
 &=(2\pi\sqrt{-1})^{-i}e^{-2\pi a_1}\int_{c_1-\sqrt{-1}\infty}^{c_1+\sqrt{-1}\infty}\frac{ds_1}{2\pi \sqrt{-1}}\,\int_{c_2-\sqrt{-1}\infty}^{c_2+\sqrt{-1}\infty}\frac{ds_2}{2\pi \sqrt{-1}}\,2^{-s_1-s_2}\Gamma_\R(s_1+\lambda_{1}+1)\Gamma_\R(-s_2-\lambda_{2})\\
&\quad\quad\quad\quad\quad\quad\quad\quad\quad\quad\times\Gamma_\R(s_1+s_2+\lambda_{1}-\lambda_{2}+2)\Gamma_\R(s_1+s_2+\lambda_{1}+\lambda_{2}+2)\frac{\Gamma(s_1+\lambda_{1}+1+i)}{\Gamma(s_1+\lambda_{1}+1)}\\
&\quad\quad\quad\quad\quad\quad\quad\quad\quad\quad\quad\quad\quad\quad\quad\quad\quad\quad\quad\quad\quad\quad\quad\quad\quad\quad\quad\quad\quad\quad\times a_1^{(-s_1-s_2+{\sf u})/{2}}|a_2|^{-s_1-i}
\end{align*}
for $a_1,a_2 \in \R^\times$, $a_1>0$.
Here $c_1,c_2 \in \R$ satisfy
\[
c_1+c_2+\lambda_{1}+\lambda_{2}+2>0,\quad c_1+\lambda_{1}+1>0>c_2+\lambda_{2}. 
\]
\end{thm}

In the following lemma, we specialize the result 
\cite[Theorem 3.2.1]{BHR1994} to $\GSp_4(\R)$ (see also \cite[Theorem 4.6.2]{Harris1990}).
\begin{lemma}\label{T:BHR}
For the (limit of) discrete series representation $D_{(\underline{\lambda};\,{\sf u})}^{(i)}$ of $\GSp_4^+(\R)$, there exist a unique $q\geq0$ and irreducible algebraic representation $(\rho,V_\rho)$ of $K_2$ such that
\[
H^q(\frak{P}_2,K_2;D_{(\underline{\lambda};\,{\sf u})}^{(i)}\otimes V_\rho) \neq0.
\]
More precisely, $q=i$ and 
\[
\rho = \begin{cases}
\rho_{(\lambda_2-2,-\lambda_1;\,-{\sf u})} & \mbox{ if $i=0$},\\
\rho_{(\underline{\lambda}^\vee-(2,0);\,-{\sf u})} & \mbox{ if $i=1$},\\
\rho_{(\underline{\lambda}-(3,1);\,-{\sf u})} & \mbox{ if $i=2$},\\
\rho_{(\lambda_1-3,-\lambda_2-1;\,-{\sf u})} & \mbox{ if $i=3$}.\\
\end{cases}
\]
In this case, $H^q(\frak{P}_2,K_2;D_{(\underline{\lambda};\,{\sf u})}^{(i)}\otimes V_\rho)$ is one-dimensional and
\[
H^q(\frak{P}_2,K_2;D_{(\underline{\lambda};\,{\sf u})}^{(i)}\otimes V_\rho) = \left(D_{(\underline{\lambda};\,{\sf u})}^{(i)} \otimes \exterior{q}\frak{p}_2^+ \otimes V_\rho\right)^{K_2}.
\]
\end{lemma}

\begin{proof}
We have four positive systems of $(\frak{g}_{2,\C},\frak{h}_\C)$ that containing $\Delta_c^+$, which are given as follows:
\begin{align*}
\Delta_0^+ &= \Delta^+,\\
\Delta_1^+ &= \left\{ (1,-1;0),(2,0;0),(0,-2;0),(1,1;0)\right\},\\
\Delta_2^+ &= \left\{ (1,-1;0),(2,0;0),(0,-2;0),(-1,-1;0)\right\},\\
\Delta_3^+ &= \left\{ (1,-1;0),(-2,0;0),(0,-2;0),(-1,-1;0)\right\}.
\end{align*}
The set $\mathcal{F}\subset \frak{h}_\C^*$ of algebraic differentials is given by
\[
\mathcal{F} = \left\{ (a,b;c) \in \Z^3\,\vert\,a+b\equiv c \,({\rm mod}\,2)\right\}.
\]
Let $\delta = (2,1;0)$ be the half-sum of positive roots in $\Delta^+$.
For $\Lambda \in \mathcal{F}+\delta$ and $\Delta_i^+$ such that $\Lambda$ is dominant with respect to $\Delta_i^+$ and non-singular with respect to $\Delta_c^+$, we have the (limit of) discrete series representation $D(\Lambda,\Delta_i^+)$ of $\GSp_4^+(\R)$ defined as in \cite[XII, \S\,7]{Knapp1986}.
By \cite[Theorem 3.2.1]{BHR1994}, there exist a unique $q\geq 0$ and $(\rho,V_\rho)$ such that the assertions hold for $D(\Lambda,\Delta_i^+)$. Moreover, we have 
\[
q=i,\quad\rho = \rho_{\Lambda^\vee-\delta}.
\]
Finally, note that $D_{(\underline{\lambda};\,{\sf u})}^{(i)} = D(\Lambda,\Delta_i^+)$ with
\[
\Lambda = \begin{cases}
(\lambda_1-1,-\lambda_2;\,{\sf u})& \mbox{ if $i=0$},\\
(\lambda_1-1,\lambda_2;\,{\sf u}) & \mbox{ if $i=1$},\\
(-\lambda_2,-\lambda_1+1;\,{\sf u}) & \mbox{ if $i=2$},\\
(\lambda_2,-\lambda_1+1;\,{\sf u}) & \mbox{ if $i=3$}.
\end{cases}
\]
\end{proof}

\subsubsection{Rational structure via the Whittaker model}

Let $(\itSigma,V_\itSigma)$ be a globally generic $C$-algebraic irreducible cuspidal automorphic representation of $\GSp_4(\A_\F)$ with central character $\omega_\itSigma$. Here we follow \cite[Definition 5.11]{BG2014} for the notion of $C$-algebraicity. 
We assume further that the following condition is satisfied:
\begin{align}\label{E:discrete condition 2}
\mbox{$\itSigma_v$ is a (limit of) discrete series representation for all $v\in S_\infty$.}
\end{align}
Then there exists $(\underline{\lambda};\,{\sf u}) \in (\Z^2)^{S_\infty} \times \Z$ with $\underline{\lambda} = (\underline{\lambda}_v)_{v \in S_\infty}$ and $\underline{\lambda}_v = (\lambda_{1,v},\lambda_{2,v})$ 
satisfying the following conditions:
\begin{itemize}
\item $|\omega_\itSigma| = |\mbox{ }|_{\A_\F}^{\sf u}$.
\item $1-\lambda_{1,v} \leq \lambda_{2,v} \leq 0$,  $\underline{\lambda}_v \neq (1,0)$, and $\lambda_{1,v}+\lambda_{2,v} \equiv{\sf u}\,({\rm mod}\,2)$
for all $v \in S_\infty$.
\item $\itSigma_v = \itSigma_{(\underline{\lambda}_v;\,{\sf u})}^{\rm gen}$ for all $v \in S_\infty$.
\end{itemize}
We say $\itSigma$ is regular if $2-\lambda_{1,v} \leq \lambda_{2,v} \leq -1$ for all $v \in S_\infty$.
We call $(\underline{\lambda};\,{\sf u})$ the weight of $\itSigma$. 
Let $(\rho_{(\underline{\lambda};\,-{\sf u})},V_{(\underline{\lambda};\,-{\sf u})})$ be the irreducible motivic algebraic representation of $K_2^{S_\infty}$ defined by
\[
(\rho_{(\underline{\lambda};\,-{\sf u})},V_{(\underline{\lambda};\,-{\sf u})}) = \left(\bigotimes_{v \in S_\infty}\rho_{(\underline{\lambda}_v;\,-{\sf u})} , \bigotimes_{v \in S_\infty}V_{(\underline{\lambda}_v;\,-{\sf u})}\right).
\]
Let $\<\,,\,\>_{\underline{\lambda}} = \prod_{v \in S_\infty}\<\,,\,\>_{\underline{\lambda}_v}$ be the ${\rm U}(2)^{S_\infty}$-equivariant pairing on $V_{\underline{\lambda}} \times V_{\underline{\lambda}^\vee}$.
Let $\psi$ be a non-trivial additive character of $\F\backslash\A_\F$. Let $\itSigma_f = \bigotimes_{v \nmid \infty}\itSigma_v$ be the finite part of $\itSigma$ and $\mathcal{W}(\itSigma_f,\psi_{U,f})$ the space of Whittaker functions of $\itSigma_f$ with respect to $\psi_{U,f}$. 
For $\varphi \in V_\itSigma$, let $W_{\varphi,\psi_U}$ be the Whittaker function of $\varphi$ with respect to $\psi_U$ defined by
\[
W_{\varphi,\psi_U}(g) = \int_{U(\F) \backslash U(\A_\F)}\varphi\left(u g \right)\overline{\psi_U(u)}\,du^{\rm Tam}.
\]
Here $du^{\rm Tam}$ is the Tamagawa measure on $U(\A_\F)$.
Let $V_\itSigma^+$ be the space of vector-valued cusp forms in $V_\itSigma$ of weight $(\underline{\lambda}^\vee;\,{\sf u})$, that is, $V_\itSigma^+$ consisting of functions $\varphi : \GSp_4(\A_\F) \rightarrow V_{(\underline{\lambda};\,-{\sf u})}$ satisfying the following conditions:
\begin{itemize}
\item $\varphi(gk) = \rho_{(\underline{\lambda};\,-{\sf u})}(k)^{-1}\varphi(g)$
for all $ k  \in K_2^{S_\infty}$ and $g \in \GSp_4(\A_\F)$.
\item The function $g \mapsto \<\varphi(g),{\bf v}\>_{\underline{\lambda}}$ is a cusp form in $V_\itSigma$ for all ${\bf v} \in V_{(\underline{\lambda}^\vee;\,{\sf u})}$.
\end{itemize}
For $\varphi \in V_\itSigma^+$ and $\underline{i} = (i_v)_{v \in S_\infty}$ with $0 \leq i_v \leq \lambda_{1,v}-\lambda_{2,v}$, let ${\rm pr}_{\underline{i}}(\varphi) \in V_\itSigma$ be the $\underline{i}$-th component of $\varphi$ defined by
\begin{align}\label{E:i-th component}
{\rm pr}_{\underline{i}}(\varphi)(g) = \left\<\varphi(g),\, \bigotimes_{v \in S_\infty}x^{i_v}y^{\lambda_{1,v}-\lambda_{2,v}-i_v}\right\>_{\underline{\lambda}}.
\end{align}
Under the diagonal embedding ${\rm U}(1) \times {\rm U}(1) \rightarrow {\rm U}(2)$, the weight of ${\rm pr}_{\underline{i}}(\varphi)$ at $v \in S_\infty$ is equal to
\[
(-\lambda_{1,v}+i_v,-\lambda_{2,v}-i_v).
\]
It is clear that we have a $\GSp_4(\A_{\F,f})$-equivariant isomorphism $V_\itSigma^+ \rightarrow (V_\itSigma \otimes V_{(\underline{\lambda};\,-{\sf u})})^{K_2^{S_\infty}}$ given by
\begin{align}\label{E:equivariant isom.}
\varphi \longmapsto \sum_{\underline{i}}\prod_{v \in S_\infty}(-1)^{i_v}{\lambda_{1,v}-\lambda_{2,v} \choose i_v}\cdot{\rm pr}_{\underline{i}}(\varphi)\otimes \bigotimes_{v \in S_\infty}x^{\lambda_{1,v}-\lambda_{2,v}-i_v}y^{i_v}.
\end{align}
Take $\psi=\psi_\F$ to be the standard additive character. For $\varphi \in V_\itSigma^+$, let $W_\varphi^{(\infty)} \in \mathcal{W}(\itSigma_f,\psi_{U,f})$ be the unique Whittaker function so that
\[
W_{{\rm pr}_{\underline{i}}(\varphi),\psi_U} = \prod_{v \in S_\infty}W_{((\underline{\lambda}_v;\,{\sf u}),\,i_v)}^+\cdot W_\varphi^{(\infty)}
\]
for all $\underline{i} = (i_v)_{v \in S_\infty}$ with $0 \leq i_v \leq \lambda_{1,v}-\lambda_{2,v}$.
Then the map $\varphi \mapsto W_\varphi^{(\infty)}$ defines a $\GSp_4(\A_{\F,f})$-equivariant isomorphism from $V_\itSigma^+$ to $\mathcal{W}(\itSigma_f,\psi_{U,f})$.
For $\sigma \in {\rm Aut}(\C)$, let 
\[
{}^\sigma\!\itSigma = {}^\sigma\!\itSigma_\infty \otimes {}^\sigma\!\itSigma_f
\]
be the irreducible admissible representation of $\GSp_4(\A_\F)$ defined so that the $v$-component of ${}^\sigma\!\itSigma_\infty$ is isomorphic to $\itSigma_{\sigma^{-1}\circ v}$ for each $v \in S_\infty$.
In Proposition \ref{P:rational structure 2} below, we will show that ${}^\sigma\!\itSigma$ is cuspidal automorphic and globally generic. It is clear that ${}^\sigma\!\itSigma$ is $C$-algebraic of weight $({}^\sigma\!\underline{\lambda};\,{\sf u})$, where ${}^\sigma\!\underline{\lambda} = (\underline{\lambda}_{\sigma^{-1}\circ v})_{v \in S_\infty}$.
Let $\Q(\underline{\lambda})$ and $\Q(\itSigma)$ be the fixed fields of $\left\{\sigma \in {\rm Aut}(\C)\,\vert\,{}^\sigma\!\underline{\lambda} = \underline{\lambda}\right\}$ and $\left\{\sigma \in {\rm Aut}(\C)\,\vert\,{}^\sigma\!\itSigma_f = \itSigma_f\right\}$, respectively.
Note that $\Q(\underline{\lambda})\subset \Q(\itSigma)$ by the strong multiplicity one theorem for $\GL_4$.
Let 
\[
t_\sigma : \mathcal{W}(\itSigma_f,\psi_{U,f}) \longrightarrow \mathcal{W}({}^\sigma\!\itSigma_f,\psi_{U,f})
\]
be the $\sigma$-linear $\GSp_4(\A_{\F,f})$-equivariant isomorphism  defined by
\[
t_\sigma W(g) = \sigma \left( W({\rm diag}(u_\sigma^{-3},u_\sigma^{-2},1,u_\sigma^{-1})g)\right).
\]
Here $u_\sigma \in \widehat{\Z}^\times \subset \hat{\o}_\F^\times$ is the unique element such that $\sigma(\psi(x)) = \psi(u_\sigma x)$ for all $x \in \A_{\F,f}$.
Let 
\begin{align*}
V_\itSigma^+ \longrightarrow V_{{}^\sigma\!\itSigma}^+,\quad \varphi\longmapsto {}^\sigma\!\varphi
\end{align*}
be the $\sigma$-linear $\GSp_4(\A_{\F,f})$-equivariant isomorphism defined so that 
\[
W_{{}^\sigma\!\varphi}^{(\infty)} = t_{\sigma}W_{\varphi}^{(\infty)}.
\]
We thus obtain a $\Q(\itSigma)$-rational structure $(V_\itSigma^+)^{{\rm Aut}(\C/\Q(\itSigma))}$ on $V_\itSigma^+$ given by taking the Galois invariants:
\begin{align}\label{E:rational structure 2}
(V_\itSigma^+)^{{\rm Aut}(\C/\Q(\itSigma))} = \left.\left\{\varphi\in V_\itSigma^+\,\right\vert\,{}^\sigma\!\varphi=\varphi\mbox{ for $\sigma \in {\rm Aut}(\C/\Q(\itSigma))$}\right\}.
\end{align}

\subsubsection{Rational structure via the coherent cohomology}
In this section, we define the rational structures on $V_\itSigma^+$ given by the coherent cohomology of automorphic vector bundles on ${\rm Sh}(G_2,X_2)$.
For $I\subset S_\infty$, define $\underline{\lambda}(I) = (\underline{\lambda}_v(I))_{v \in S_\infty}$ and $\underline{\lambda}^I = (\underline{\lambda}_v^I)_{v \in S_\infty}$ in $(\Z^2)^{S_\infty}$ by
\[
\underline{\lambda}_v(I) = \begin{cases}
\underline{\lambda}_v-(3,1) & \mbox{ if $v \in I$}\\
\underline{\lambda}_v^\vee-(2,0) & \mbox{ if $v \notin I$}
\end{cases},\quad \underline{\lambda}_v^I = \begin{cases}
\underline{\lambda}_v & \mbox{ if $v \in I$}\\
\underline{\lambda}_v^\vee & \mbox{ if $v \notin I$}
\end{cases}.
\]
By definition, we have ${}^\sigma\!(\underline{\lambda}(I)) = {}^\sigma\!\underline{\lambda}({}^\sigma\!I)$ for $\sigma\in{\rm Aut}(\C)$.
It is clear that $\underline{\lambda}(I)$ is invariant by ${\rm Aut}(\C/\Q(\underline{\lambda})\Q(I))$.
We say $I$ is $admissible$ with respect to $\underline{\lambda}$ if the following condition is satisfied:
Let $0 \leq j_v \leq 3$ for $v \in S_\infty$. We have
\begin{align}\label{E:admissibility}
H^{d+{}^\sharp I}\left(\frak{P}_{2}^{S_\infty},K_{2}^{S_\infty};\bigotimes_{v \in S_\infty}D_{(\underline{\lambda}_v;\,{\sf u})}^{(j_v)}\otimes V_{(\underline{\lambda}(I);\,-{\sf u})}\right) \neq 0
\end{align}
if and only if 
\[
j_v = \begin{cases}
2 & \mbox{ if $v \in I$},\\
1 & \mbox{ if $v \notin I$}.
\end{cases}
\]
By Lemma \ref{T:BHR}, the admissibility is a combinatorial condition on $I$. For instance, $I$ is admissible when the following condition is satisfied:
\begin{itemize}
\item For $v \in I$ (resp.\,$v \notin I$), we have $\lambda_{1,v}+\lambda_{2,v}\geq 2$ (resp.\,$\lambda_{2,v} \leq -1$).
\end{itemize} 
In particular, any subset of $S_\infty$ is admissible if $\itSigma_v$ is a discrete series representation for all $v \in S_\infty$. Also note that $I$ is admissible if and only if $S_\infty \smallsetminus I$ is admissible.
Let 
\[
[\mathcal{V}_{(\underline{\lambda}(I);\,-{\sf u})}]
\]
denote the automorphic vector bundle on ${\rm Sh}(G_2,X_2)$ defined by the motivic algebraic representation 
\[
(\rho_{(\underline{\lambda}(I);\,-{\sf u})},V_{(\underline{\lambda}(I);\,-{\sf u})}) = \left( \bigotimes_{v \in S_\infty}\rho_{({\lambda}_v(I);\,-{\sf u})}, \bigotimes_{v \in S_\infty}V_{({\lambda}_v(I);\,-{\sf u})}\right)
\]
of $K_2^{S_\infty}$. Let $H_?^q([\mathcal{V}_{(\underline{\lambda}(I);\,-{\sf u})}])[\itSigma_f]$ be the $\itSigma_f$-isotypic component of $H_?^q([\mathcal{V}_{(\underline{\lambda}(I);\,-{\sf u})}])$ for $? \in \{{\rm cusp, !, (2)}\}$.
In the following proposition, we show that being globally generic is an arithmetic property of an $C$-algebraic irreducible cuspidal automorphic representation of $\GSp_4(\A_\F)$ satisfying condition (\ref{E:discrete condition 2}). The result was proved in \cite[Lemma 3.1]{Chen2021c} when $\itSigma_v$ is a discrete series representation for all $v \in S_\infty$.

\begin{prop}\label{P:rational structure 2}
For $\sigma \in {\rm Aut}(\C)$, the representation ${}^\sigma\!\itSigma$ is cuspidal automorphic and globally generic.
\end{prop}

\begin{proof}
By Theorem \ref{T:Harris}-(2) and Lemma \ref{T:BHR}, we have 
$H_{\rm cusp}^{2d}([\mathcal{V}_{(\underline{\lambda}(S_\infty);\,-{\sf u})}])[\itSigma_f]\neq0$.
By Theorem \ref{T:Harris}-(3), this implies that $H_!^{2d}([\mathcal{V}_{(\underline{\lambda}(S_\infty);\,-{\sf u})}])[\itSigma_f]\neq0$.
Let $\sigma \in {\rm Aut}(\C)$. It follows from (\ref{E:Galois equiv. class}) that 
\[
H_!^{2d}([\mathcal{V}_{({}^\sigma\!\underline{\lambda}(S_\infty);\,-{\sf u})}])[{}^\sigma\!\itSigma_f] = T_\sigma \left( H_!^{2d}([\mathcal{V}_{(\underline{\lambda}(S_\infty);\,-{\sf u})}])[\itSigma_f]\right)\neq0.
\]
By Theorem \ref{T:Harris}-(2) and (4), there exists an irreducible discrete automorphic representation $\itSigma'$ of $\GSp_4(\A_\F)$ such that $\itSigma_f' = {}^\sigma\!\itSigma_f$ and 
\begin{align}\label{E:rational structure proof 1}
H^{2d}\left(\frak{P}_2^{S_\infty},K_2^{S_\infty};\itSigma_\infty'\otimes V_{({}^\sigma\!\underline{\lambda}(S_\infty);\,-{\sf u}})\right) \neq0.
\end{align}
In particular, by \cite[Proposition 4.3.2]{Harris1990}, the infinitesimal characters of $\itSigma_\infty'$ and ${}^\sigma\!\itSigma_\infty$ are equal.
On the other hand, consider the transfer $\itPsi'$ of $\itSigma'$ to $\GL_4(\A_\F)$ with respect to the spin representation of $\GSp_4(\C)$. 
Since $\itSigma'$ is almost locally generic, it follows that $\itPsi'$ is an isobaric automorphic representation of $\GL_4(\A_\F)$. More precisely, the global Arthur parameter associated to $\itSigma'$ is of types (a) or (b) in the notation of \cite[Remark 6.1.8]{GT2019}.
The condition on infinitesimal character of $\itSigma_\infty'$ then implies that $\itPsi'$ is $C$-algebraic.
By the purity lemma \cite[Lemme 4.9]{Clozel1990}, we deduce that $\itPsi_\infty'$ is essentially tempered. 
It then follows that $\itSigma_\infty'$ is also essentially tempered.
Combine with (\ref{E:rational structure proof 1}), by \cite[Theorem 3.5]{Harris1990b}, $\itSigma_\infty'$ is a (limit of) discrete series representation of $\GSp_4(\F_\infty)$.
By Lemma \ref{T:BHR} and (\ref{E:rational structure proof 1}) again, we conclude that
\[
\itSigma_v' \in L_{({}^\sigma\!\underline{\lambda}_v;\,{\sf u})}
\]
for all $v \in S_\infty$. Here $L_{({}^\sigma\!\underline{\lambda}_v;\,{\sf u})}$ is the $L$-packet of $\GSp_4(\R)$ defined in (\ref{E:L-packet}).
Therefore, $\itSigma'$ and ${}^\sigma\!\itSigma$ belong to the same global Arthur packet. We then deduce that ${}^\sigma\!\itSigma$ is discrete automorphic by Arthur's multiplicity formula for $\GSp_4(\A_\F)$ established by Gee and Ta\"ibi \cite[Theorem 7.4.1]{GT2019}.
Since ${}^\sigma\!\itSigma_\infty$ is essentially tempered, this implies that ${}^\sigma\!\itSigma$ is cuspidal automorphic by \cite{Wallach1984}.
Finally, as explained in \cite[Remark 7.4.7]{GT2019}, among the global Arthur packet associated to $\itSigma'$, there is a unique globally generic discrete automorphic representation, which is characterized by the condition that it is locally generic at all places. Hence ${}^\sigma\!\itSigma$ must be globally generic.
This completes the proof.
\end{proof}

Let $I\subset S_\infty$. In the following lemma, we define a $\GSp_4(\A_{\F,f})$-equivariant embedding from $V_\itSigma^+$ to the isotypic component $H_{\rm cusp}^{d+{}^\sharp I}(\mathcal{V}_{(\underline{\lambda}(I);\,-{\sf u})})[\itSigma_f]$.
Let 
\[
\xi_{\underline{\lambda}}^I : V_{(\underline{\lambda}^I;\,-{\sf u})}\longrightarrow\exterior{d+{}^\sharp I}(\frak{p}_2^+)^{S_\infty} \otimes V_{(\underline{\lambda}(I);\,-{\sf u})}
\]
be the $K_2^{S_\infty}$-equivariant homomorphism defined by
\[
\xi_{\underline{\lambda}}^I = \left(\bigotimes_{v \in I}\xi_{\underline{\lambda}_v}^+\right) \otimes \left(\bigotimes_{v \notin I}\xi_{\underline{\lambda}_v}^-\right).
\]
Here $\xi_{\underline{\lambda}_v}^\pm$ are defined in (\ref{E:U(2) embeddings}).

\begin{lemma}\label{L:Eichler-Shimura 2}
Let $I \subset S_\infty$. We have a $\GSp_4(\A_{\F,f})$-equivariant embedding 
\[
V_\itSigma^+ \longrightarrow H_{\rm cusp}^{d+{}^\sharp I}(\mathcal{V}_{(\underline{\lambda}(I);\,-{\sf u})})[\itSigma_f],\quad \varphi\longmapsto [\varphi]_I
\]
defined by
\[
[\varphi]_I = \sum_{\underline{i}} \prod_{v \in S_\infty}{\lambda_{1,v}-\lambda_{2,v} \choose i_v}\cdot {\rm pr}_{\underline{i}}(\varphi)^I\otimes \xi_{\underline{\lambda}}^I\left(\left( \bigotimes_{v \in I}(-1)^{i_v}x^{\lambda_{1,v}-\lambda_{2,v}-i_v}y^{i_v}\right)\otimes \left(\bigotimes_{v \notin I}x^{i_v}y^{\lambda_{1,v}-\lambda_{2,v}-i_v} \right)\right).
\]
Here $\underline{i} = (i_v)_{v \in S_\infty} \in \Z^{S_\infty}$ with $0 \leq i_v \leq \lambda_{1,v}-\lambda_{2,v}$ and
\[
{\rm pr}_{\underline{i}}(\varphi)^I(g) = {\rm pr}_{\underline{i}}(\varphi)\left(g \cdot \prod_{v \in S_\infty\smallsetminus I}{\rm diag}(-1,-1,1,1)\right).
\]
\end{lemma}

\begin{proof}
By Theorem \ref{T:Harris}-(2), we have
\[
\left(V_\itSigma \otimes \exterior{d+{}^\sharp I}(\frak{p}_2^+)^{S_\infty} \otimes V_{(\underline{\lambda}(I);\,-{\sf u})} \right)^{K_2^{S_\infty}} \subset H_{\rm cusp}^{d+{}^\sharp I}(\mathcal{V}_{(\underline{\lambda}(I);\,-{\sf u})})[\itSigma_f].
\]
We also have a $\GSp_4(\A_{\F,f})$-equivariant embedding 
\[
{\rm id} \otimes \xi_{\underline{\lambda}}^I : (V_\itSigma \otimes V_{(\underline{\lambda}^I;\,-{\sf u})})^{K_2^{S_\infty}} \longrightarrow \left(V_\itSigma \otimes \exterior{d+{}^\sharp I}(\frak{p}_2^+)^{S_\infty} \otimes V_{(\underline{\lambda}(I);\,-{\sf u})} \right)^{K_2^{S_\infty}} .
\]
By considering the minimal $K_2^{S_\infty}$-types of $\itSigma_\infty \vert_{\GSp_4^+(\F_\infty)}$, we easily see that $(V_\itSigma \otimes V_{(\underline{\lambda}^I;\,-{\sf u})})^{K_2^{S_\infty}}$ is one-dimensional. 
By (\ref{E:equivariant isom.}), we have a $\GSp_4(\A_{\F,f})$-equivariant isomorphism $V_\itSigma^+ \rightarrow (V_\itSigma \otimes V_{(\underline{\lambda}^I;\,-{\sf u})})^{K_2^{S_\infty}}$ given by
\[
\varphi \longmapsto \sum_{\underline{i}} \prod_{v \in S_\infty}{\lambda_{1,v}-\lambda_{2,v} \choose i_v}\cdot {\rm pr}_{\underline{i}}(\varphi)^I\otimes \left( \bigotimes_{v \in I}(-1)^{i_v}x^{\lambda_{1,v}-\lambda_{2,v}-i_v}y^{i_v}\right)\otimes \left(\bigotimes_{v \notin I}x^{i_v}y^{\lambda_{1,v}-\lambda_{2,v}-i_v} \right).
\]
Composing with ${\rm id}\otimes\xi_{\underline{\lambda}}^I$, we then obtain the homomorphism $\varphi\mapsto[\varphi]_I$. This completes the proof.
\end{proof}

\begin{lemma}\label{L:rational structure 2}
Let $I \subset S_\infty$.
\begin{itemize}
\item[(1)] For all $\sigma \in {\rm Aut}(\C)$ and $q\geq 0$, we have
\[
T_\sigma(H_{\rm cusp}^q([\mathcal{V}_{(\underline{\lambda}(I);\,-{\sf u})}])[\itSigma_f]) = H_{{\rm cusp}}^q([\mathcal{V}_{({}^\sigma\!\underline{\lambda}({}^\sigma\!I);\,-{\sf u})}])[{}^\sigma\!\itSigma_f].
\]
In particular, we have a $\Q(\itSigma)\Q(I)$-rational structure on $H_{\rm cusp}^q([\mathcal{V}_{(\underline{\lambda}(I);\,-{\sf u})}])[\itSigma_f]$ given by taking the Galois invariants:
\begin{align*}
&H_{\rm cusp}^q([\mathcal{V}_{(\underline{\lambda}(I);\,-{\sf u})}])[\itSigma_f]^{{\rm Aut}(\C/\Q(\itSigma)\Q(I))}\\
& = \left.\left\{c\in H_{\rm cusp}^q([\mathcal{V}_{(\underline{\lambda}(I);\,-{\sf u})}])[\itSigma_f]\,\right\vert\,T_\sigma c=c\mbox{ for $\sigma \in {\rm Aut}(\C/\Q(\itSigma)\Q(I))$}\right\}.
\end{align*}
\item[(2)]
Assume $I$ is admissible with respect to $\underline{\lambda}$, then $H_{\rm cusp}^{d+{}^\sharp I}([\mathcal{V}_{(\underline{\lambda}(I);\,-{\sf u})}])[\itSigma_f] \simeq \itSigma_f$.
\end{itemize}
\end{lemma}

\begin{proof}
By (\ref{E:Galois equiv. class}) and Proposition \ref{P:rational structure 2}, to prove the first assertion,
it suffices to show that
\begin{align}\label{E:rational structure 2 proof 1}
H_{\rm cusp}^q([\mathcal{V}_{(\underline{\lambda}(I);\,-{\sf u})}])[\itSigma_f] = H_!^q([\mathcal{V}_{(\underline{\lambda}(I);\,-{\sf u})}])[\itSigma_f].
\end{align}
A discrete irreducible cuspidal automorphic representation $\itSigma'$ appears in $H_{(2)}^q([\mathcal{V}_{(\underline{\lambda}(I);\,-{\sf u})}])[\itSigma_f]$ only when $\itSigma'_f = \itSigma_f$. In particular, $\itSigma'$ and $\itSigma$ belongs to the same global Arthur packet. Since $\itSigma_v$ is essentially tempered and generic for all $v\in S_\infty$, the local Arthur packets at archimedean places are actually the local $L$-packets $L_{(\underline{\lambda}_v;\,{\sf u})}$ in (\ref{E:L-packet}). Hence $\itSigma_v'$ is also essentially tempered for all $v \in S_\infty$. This implies that $\itSigma'$ is cuspidal by \cite{Wallach1984}. We conclude that the multiplicities of $\itSigma_f$ in $H_{(2)}^q([\mathcal{V}_{(\underline{\lambda}(I);\,-{\sf u})}])$ and $H_{\rm cusp}^q([\mathcal{V}_{(\underline{\lambda}(I);\,-{\sf u})}])$ are equal.
Equality (\ref{E:rational structure 2 proof 1}) then follows from Theorem \ref{T:Harris}-(3) and (4).
Also we deduce that the multiplicity of $\itSigma_f$ in $H_{\rm cusp}^{d+{}^\sharp I}([\mathcal{V}_{(\underline{\lambda}(I);\,-{\sf u})}])$ is less than or equal to the number of tuples $\underline{j} = (j_v)_{v \in S_\infty}$ with $0\leq j_v \leq 3$ and such that (\ref{E:admissibility}) hold. 
By the definition of admissibility, this number is equal to one. Therefore, $\itSigma_f$ appears in $H_{\rm cusp}^{d+{}^\sharp I}([\mathcal{V}_{(\underline{\lambda}(I);\,-{\sf u})}])$ with multiplicity one by Lemma \ref{L:Eichler-Shimura 2}. In fact, the archimedean component of ${\rm pr}_{\underline{i}}(\varphi)^I$ belongs to 
\[
\left(\bigotimes_{v \in I}D^{(2)}_{(\underline{\lambda}_v;\,{\sf u})} \right)\otimes \left(\bigotimes_{v \notin I}D^{(1)}_{(\underline{\lambda}_v;\,{\sf u})}\right).
\]
for all $\varphi \in V_\itSigma^+$ and $\underline{i} = (i_v)_{ \in S_\infty} \in \Z^{S_\infty}$ with $0 \leq i_v \leq \lambda_{1,v}-\lambda_{2,v}$.
This completes the proof.
\end{proof}





\subsubsection{Automorphic periods and period relations}

By comparing the rational structures in (\ref{E:rational structure 2}) and Lemma \ref{L:rational structure 2}-(1), we 
have the following lemma/definition for the automorphic periods of $\itSigma$.

\begin{lemma}
Let $I\subset S_\infty$ be admissible with respect to $\underline{\lambda}$. There exists a sequence of non-zero complex numbers $\left(p^{{}^\sigma\!I}({}^\sigma\!\itSigma)\right)_{\sigma \in {\rm Aut}(\C)}$ such that
\[
T_\sigma\left(\frac{[\varphi]_I}{p^I(\itSigma)}\right) = \frac{[{}^\sigma\!\varphi]_{{}^\sigma\!I}}{p^{{}^\sigma\!I}({}^\sigma\!\itSigma)}
\]
for all $\sigma \in {\rm Aut}(\C)$ and $\varphi \in V_\itSigma^+$. Here $T_\sigma : H_!^{d+{}^\sharp I}([\mathcal{V}_{(\underline{\lambda}(I);\,-{\sf u})}]) \rightarrow H_!^{d+{}^\sharp I}([\mathcal{V}_{({}^\sigma\!\underline{\lambda}({}^\sigma\!I);\,-{\sf u})}])$ is the $\sigma$-linear isomorphism in (\ref{E:Galois equiv. class}).
\end{lemma}

In the following theorem, we prove a period relation for product of automorphic periods and critical value of adjoint $L$-function. 
The result is an analogue of Lemma \ref{L:period relation 1}-(2).
Let $L(s,\itSigma,{\rm Ad})$ be the adjoint $L$-function of $\itSigma$, where ${\rm Ad}$ is the adjoint representation of $\GSp_4(\C)$ on $\frak{pgsp}_4(\C)$.

\begin{thm}\label{T:Ichino-Chen}
Let $I \subset S_\infty$ be admissible with respect to $\underline{\lambda}$.
For $\sigma \in {\rm Aut}(\C)$, we have
\[
\sigma \left( \frac{L(1,\itSigma,{\rm Ad})}{\pi^{3d}\cdot p^I(\itSigma)\cdot p^{S_\infty\smallsetminus I}(\itSigma^\vee)} \right) = \frac{L(1,{}^\sigma\!\itSigma,{\rm Ad})}{\pi^{3d}\cdot p^{{}^\sigma\!I}({}^\sigma\!\itSigma)\cdot p^{S_\infty\smallsetminus {}^\sigma\!I}({}^\sigma\!\itSigma^\vee)}.
\]
\end{thm}

\begin{proof}
The key ingredients of the proof are the Serre duality for coherent cohomology and our previous result \cite{Chen2021c}.
Let $J=S_\infty \smallsetminus I$. 
Note that $\itSigma^\vee$ is $C$-algebraic of weight $(\underline{\lambda};\,-{\sf u})$ and $J$ is also admissible with respect to $\underline{\lambda}$.
We have $K_2^{S_\infty}$-equivariant pairings
\[
\exterior{d+{}^\sharp I}(\frak{p}_2^+)^{S_\infty} \times \exterior{d+{}^\sharp J}(\frak{p}_2^+)^{S_\infty} \longrightarrow \exterior{3d}(\frak{p}_2^+)^{S_\infty},\quad (X,Y)\longmapsto X\wedge Y
\]
and
\begin{align*}
V_{(\underline{\lambda}(I);\,-{\sf u})} \times V_{(\underline{\lambda}(J);\,{\sf u})} &\longrightarrow \bigotimes_{v \in S_\infty}V_{(-3,-3;\,0)},\quad
({\bf v},{\bf w})\longmapsto \<{\bf v},{\bf w}\>_{\underline{\lambda}(I)}.
\end{align*}
Identify $\exterior{3d}(\frak{p}_2^+)^{S_\infty}$ with $\bigotimes_{v \in S_\infty}V_{(3,3;\,0)}$ by the isomorphism in (\ref{E:p_2}).
We then have a $K_2^{S_\infty}$-equivariant pairing
\[
\left(\exterior{d+{}^\sharp I}(\frak{p}_2^+)^{S_\infty}\otimes V_{(\underline{\lambda}(I);\,-{\sf u})}\right) \times \left( \exterior{d+{}^\sharp J}(\frak{p}_2^+)^{S_\infty}\otimes V_{(\underline{\lambda}(J);\,{\sf u})} \right) \longrightarrow \C.
\]
When restrict to $V_{(\underline{\lambda}^I;\,-{\sf u})} \times V_{(\underline{\lambda}^J;\,{\sf u})}$ via the embedding $\xi_{\underline{\lambda}}^I \times \xi_{\underline{\lambda}}^J$, the above pairing is equal to $C\cdot\<\,,\,\>_{\underline{\lambda}^I}$ for some $C \in \Q^\times$.
We also have the Petersson bilinear pairing
\[
V_\itSigma \times V_{\itSigma^\vee} \longrightarrow \C,\quad(\varphi_1,\varphi_2)\longmapsto  \<\varphi_1,\varphi_2\>=\int_{\A_\F^\times\GSp_4(\F)\backslash\GSp_4(\A_\F)}\varphi_1(g)\varphi_2(g)\,dg^{\rm Tam}.
\]
Here $dg^{\rm Tam}$ is the Tamagawa measure on $\A_\F^\times\backslash\GSp_4(\A_\F)$.
We thus obtain a bilinear pairing
\[
H_{\rm cusp}^{d+{}^\sharp I}([\mathcal{V}_{(\underline{\lambda}(I);\,-{\sf u})}])[\itSigma_f] \times H_{\rm cusp}^{d+{}^\sharp J}([\mathcal{V}_{(\underline{\lambda}(J);\,{\sf u})}])[\itSigma_f^\vee] \longrightarrow \C,\quad (c_1,c_2)\longmapsto \int_{{\rm Sh}(G_2,X_2)}c_1\wedge c_2.
\]
The pairing satisfies the Galois equivariant property:
\[
\sigma \left( \int_{{\rm Sh}(G_2,X_2)}c_1\wedge c_2\right) = \int_{{\rm Sh}(G_2,X_2)}T_\sigma c_1\wedge T_\sigma c_2
\]
for all $\sigma \in {\rm Aut}(\C)$.
Indeed, the pairing is the restriction of the Serre duality pairing
\[
H_!^{d+{}^\sharp I}([\mathcal{V}_{(\underline{\lambda}(I);\,-{\sf u})}])\times H_!^{d+{}^\sharp J}([\mathcal{V}_{(\underline{\lambda}(J);\,{\sf u})}]) \longrightarrow \C
\]
in \cite[Proposition 3.8 and Remark 3.8.4]{Harris1990}.
For $\varphi \in V_\itSigma^+$, let $\varphi^I : \GSp_4(\A_\F) \rightarrow V_{(\underline{\lambda}^I;\,-{\sf u})}$ be the vector-valued cusp form defined by
\[
\varphi^I(g) = c_{\underline{\lambda}}^I\circ\varphi\left(g\cdot\prod_{v \in S_\infty \smallsetminus I}{\rm diag}(-1,-1,1,1)\right).
\]
Here $c_{\underline{\lambda}}^I : V_{\underline{\lambda}}\rightarrow V_{\underline{\lambda}^I}$ is the homomorphism given by
\[
c_{\underline{\lambda}}^I = \left(\bigotimes_{v \in I}{\rm id}\right)\otimes \left(\bigotimes_{v \notin I}c_{\underline{\lambda}_v}\right).
\]
For $\underline{i} = (i_v)_{v \in S_\infty}$ with $0 \leq i_v \leq \lambda_{1,v}-\lambda_{2,v}$, let ${\rm pr}_{\underline{i}}(\varphi^I) \in V_\itSigma$ be defined as in (\ref{E:i-th component}) with $\underline{\lambda}$ replaced by $\underline{\lambda}^I$.
By (\ref{E:U(2) relation}), we have
\[
{\rm pr}_{\underline{i}}(\varphi^I) = \prod_{v \notin I}(-1)^{i_v}\cdot{\rm pr}_{\underline{i}^I}(\varphi)^I.
\]
Here $\underline{i}^I = (i_v)_{v \in I} \times (\lambda_{1,v}-\lambda_{2,v}-i_v)_{v \notin I}$.
In particular, we have
\[
[\varphi]_I = \sum_{\underline{i}}\prod_{v \in S_\infty}(-1)^{i_v}{\lambda_{1,v}-\lambda_{2,v} \choose i_v}\cdot{\rm pr}_{\underline{i}}(\varphi^I)\otimes \xi_{\underline{\lambda}}^I\left(\bigotimes_{v \in S_\infty}x^{\lambda_{1,v}-\lambda_{2,v}-i_v}y^{i_v}\right).
\]
Similarly we define $\varphi^J$ for $\varphi \in V_{\itSigma^\vee}^+$. Then it is clear that
\[
\int_{{\rm Sh}(G_2,X_2)}[\varphi_1]_I\wedge[\varphi_2]_J =C \cdot\int_{\A_\F^\times\GSp_4(\F)\backslash\GSp_4(\A_\F)}\<\varphi_1^I(g),\varphi_2^J(g)\>_{\underline{\lambda}^I}\,dg^{\rm Tam}
\]
for $\varphi_1 \in V_\itSigma^+$ and $\varphi_2 \in V_{\itSigma^\vee}^+$.
By Schur's orthogonal relation, we have
\begin{align*}
&\int_{\A_\F^\times\GSp_4(\F)\backslash\GSp_4(\A_\F)}\<\varphi_1^I(g),\varphi_2^J(g)\>_{\underline{\lambda}^I}\,dg^{\rm Tam}\\
& = {\rm dim}\, V_{\underline{\lambda}}\cdot\<{\bf v},{\bf w}\>_{\underline{\lambda}^I}^{-1}\cdot\int_{\A_\F^\times\GSp_4(\F)\backslash\GSp_4(\A_\F)}\<\varphi_1^I(g),{\bf v}\>_{\underline{\lambda}^I}\<\varphi_2^J(g),{\bf w}\>_{\underline{\lambda}^J}\,dg^{\rm Tam}
\end{align*}
for any ${\bf v} \in V_{(\underline{\lambda}^I;\,-{\sf u})}$ and ${\bf w} \in V_{(\underline{\lambda}^J;\,{\sf u})}$ such that $\<{\bf v},{\bf w}\>_{\underline{\lambda}^I} \neq 0$.
We take
\[
{\bf v} = \left(\bigotimes_{v \in I}y^{\lambda_{1,v}-\lambda_{2,v}}\right)\otimes \left(\bigotimes_{v \in J}x^{\lambda_{1,v}-\lambda_{2,v}}\right),\quad {\bf w} = \left(\bigotimes_{v \in I}x^{\lambda_{1,v}-\lambda_{2,v}}\right)\otimes \left(\bigotimes_{v \in J}y^{\lambda_{1,v}-\lambda_{2,v}}\right).
\]
Then we have
\begin{align*}
\int_{{\rm Sh}(G_2,X_2)}[\varphi_1]_I\wedge[\varphi_2]_J &=C\cdot{\rm dim}\, V_{\underline{\lambda}}\cdot(-1)^{{\sf u}{}^\sharp J}\cdot \<{\rm pr}_{\underline{0}}(\varphi_1)^I,{\rm pr}_{\underline{0}}(\varphi_2)^J\>\\
&=C\cdot{\rm dim}\, V_{\underline{\lambda}}\cdot(-1)^{{\sf u}{}^\sharp J}\cdot \<{\rm pr}_{\underline{0}}(\varphi_1)^\varnothing,{\rm pr}_{\underline{0}}(\varphi_2)^{S_\infty}\>.
\end{align*}
Similarly, we have
\[
\int_{{\rm Sh}(G_2,X_2)}[{}^\sigma\!\varphi_1]_{{}^\sigma\!I}\wedge[{}^\sigma\!\varphi_2]_{{}^\sigma\!J} =C\cdot{\rm dim}\, V_{\underline{\lambda}}\cdot(-1)^{{\sf u}{}^\sharp J}\cdot \<{\rm pr}_{\underline{0}}({}^\sigma\!\varphi_1)^\varnothing,{\rm pr}_{\underline{0}}({}^\sigma\!\varphi_2)^{S_\infty}\>
\]
for all $\sigma \in {\rm Aut}(\C)$.
We thus conclude from the Galois equivariant property of the Serre duality pairing and the definition of automorphic periods that
\[
\sigma \left( \frac{ \<{\rm pr}_{\underline{0}}(\varphi_1)^\varnothing,{\rm pr}_{\underline{0}}(\varphi_2)^{S_\infty}\>}{p^I(\itSigma)\cdot p^J(\itSigma^\vee)}  \right) = \frac{ \<{\rm pr}_{\underline{0}}({}^\sigma\!\varphi_1)^\varnothing,{\rm pr}_{\underline{0}}({}^\sigma\!\varphi_2)^{S_\infty}\>}{p^{{}^\sigma\!I}({}^\sigma\!\itSigma)\cdot p^{{}^\sigma\!J}({}^\sigma\!\itSigma^\vee)} 
\]
for all $\sigma \in {\rm Aut}(\C)$.
On the other hand, by our main result \cite[Theorem 1.1]{Chen2021c}, we have
\[
\sigma\left(\frac{\pi^{3d}\cdot\<{\rm pr}_{\underline{0}}(\varphi_1)^\varnothing,{\rm pr}_{\underline{0}}(\varphi_2)^{S_\infty}\>}{L(1,\itPi,{\rm Ad})}\right) = \frac{\pi^{3d}\cdot\<{\rm pr}_{\underline{0}}({}^\sigma\!\varphi_1)^\varnothing,{\rm pr}_{\underline{0}}({}^\sigma\!\varphi_2)^{S_\infty}\>}{L(1,{}^\sigma\!\itPi,{\rm Ad})}
\]
for all $\sigma \in {\rm Aut}(\C)$. This completes the proof.
\end{proof}

\begin{rmk}
In \cite{Chen2021c}, the result was proved under the assumption that $\itSigma_v$ is a discrete series representation for all $v \in S_\infty$. However, the proof goes without change for limits of discrete series representations, as long as we can extend \cite[Lemma 3.1]{Chen2021c}. This is done in Proposition \ref{P:rational structure 2}.
\end{rmk}

\section{Proof of main result}\label{S:proof}

In this section, we prove our main result Theorem \ref{T:main}.

\subsection{The Kim--Ramakrishnan--Shahidi lifts}\label{SS:KRS lift}

Let $\itPi$ be a regular $C$-algebraic irreducible cuspidal automorphic representation of $\GL_2(\A_\F)$ with central character $\omega_\itPi$ and satisfying condition (\ref{E:discrete condition 1}). Let $(\underline{\kappa};\,{\sf w}) \in \Z^{S_\infty} \times \Z$ be the weight of $\itPi$. Assume $\itPi$ is non-CM. Let ${\rm Sym}^3(\itPi)$ be the functorial lift of $\itPi$ to $\GL_4(\A_\F)$ with respect to the symmetric cube representation of $\GL_2(\C)$.
The existence of the lift was proved by Kim and Shahidi \cite{KS2002}. Since $\itPi$ is non-CM, ${\rm Sym}^3(\itPi)$ is cuspidal automorphic.
We have the factorization of twisted exterior square $L$-function of ${\rm Sym}^3(\itPi)$ by $\omega_\itPi^{-3}$:
\[
L(s,{\rm Sym}^3(\itPi),\exterior{2}\otimes\,\omega_\itPi^{-3}) = L(s,\itPi,{\rm Sym}^4\otimes\omega_\itPi^{-2})\cdot\zeta_\F(s).
\]
Thus $L(s,{\rm Sym}^3(\itPi),\exterior{2}\otimes\,\omega_\itPi^{-3})$ has a pole at $s=1$, as the twisted symmetric fourth $L$-function is holomorphic and non-zero at $s=1$ by \cite[Corollary 7.1.5]{BLGG2011}.
By the result of Gan and Takeda \cite[Theorem 12.1]{GT2011}, ${\rm Sym}^3(\itPi)$ strongly descend to an irreducible globally generic cuspidal automorphic representation $\itSigma$ of $\GSp_4(\A_\F)$. We call it the Kim--Ramakrishnan--Shahidi lift of $\itPi$ (cf.\,\cite[Theorem C]{RS2007b} for $\omega_\itPi=1$).
By the functoriality of $\itSigma$, we have the equality between twisted spin (resp.\,standard) $L$-function of $\itSigma$ and twisted symmetric cube (resp.\,fourth) $L$-function of $\itPi$:
\begin{align}
L(s,\itSigma,{\rm spin}\otimes\chi) & = L(s,\itPi,{\rm Sym}^3\otimes\chi),\label{E:twisted spin}\\
L(s,\itSigma,{\rm std}\otimes\chi) &= L(s,\itPi,{\rm Sym}^4\otimes\chi\omega_\itPi^{-2})\label{E:twisted standard}
\end{align}
for Hecke character $\chi$ of $\A_\F^\times$.
Moreover, $\itSigma$ is regular $C$-algebraic, satisfying condition (\ref{E:discrete condition 2}) with weight $(\underline{\lambda};\,{\sf u})$, and has central character $\omega_\itPi^3$, where
\[
\underline{\lambda}_v = (2\kappa_v-1,1-\kappa_v),\quad {\sf u} = 3{\sf w}
\]
for $v \in S_\infty$. 

We have the following period relations for the automorphic periods of $\itSigma$ and powers of the Petersson norm of the normalized newform of $\itPi$, which is a crucial ingredient in the proof of Theorem \ref{T:main}.
\begin{thm}\label{P:period relation main}
Let $\itSigma$ be the Kim--Ramakrishnan--Shahidi lift of $\itPi$.
\begin{itemize}
\item[(1)] Assume $\kappa_v \geq 6$ for all $v \in S_\infty$. For $\sigma \in {\rm Aut}(\C)$, we have
\begin{align*}
& \sigma \left( \frac{p^\varnothing(\itSigma)}{|D_\F|^{1/2}\cdot(2\pi\sqrt{-1})^{\sum_{v \in S_\infty}(-\kappa_v+4-3{\sf w})/2}\cdot(\sqrt{-1})^d\cdot G(\omega_\itPi)^3\cdot\Vert f_\itPi\Vert^3}\right)\\
& = \frac{p^\varnothing({}^\sigma\!\itSigma)}{|D_\F|^{1/2}\cdot(2\pi\sqrt{-1})^{\sum_{v \in S_\infty}(-\kappa_v+4-3{\sf w})/2}\cdot(\sqrt{-1})^d\cdot G({}^\sigma\!\omega_\itPi)^3\cdot\Vert f_{{}^\sigma\!\itPi}\Vert^3}.
\end{align*}
\item[(2)] Assume $\kappa_v \geq 3$ for all $v \in S_\infty$. For $\sigma \in {\rm Aut}(\C/\F^{\Gal})$, we have
\begin{align*}
& \sigma \left( \frac{p^{S_\infty}(\itSigma)}{|D_\F|^{1/2}\cdot(2\pi\sqrt{-1})^{\sum_{v \in S_\infty}(\kappa_v+4-3{\sf w})/2}\cdot(\sqrt{-1})^d\cdot G(\omega_\itPi)^3\cdot\Vert f_\itPi\Vert^4}\right)\\
& = \frac{p^{S_\infty}({}^\sigma\!\itSigma)}{|D_\F|^{1/2}\cdot(2\pi\sqrt{-1})^{\sum_{v \in S_\infty}(\kappa_v+4-3{\sf w})/2}\cdot(\sqrt{-1})^d\cdot G({}^\sigma\!\omega_\itPi)^3\cdot\Vert f_{{}^\sigma\!\itPi}\Vert^4}.
\end{align*}
Here $\F^{\Gal}$ is the Galois closure of $\F$ in $\C$.
\item[(3)] Assume $\kappa_v = \kappa_w\geq3$ for all $v,w \in S_\infty$. Then the assertion in (2) holds with ${\rm Aut}(\C/\F^{\Gal})$ replaced by ${\rm Aut}(\C)$.
\end{itemize}
\end{thm}

\begin{proof}
The assertions will be proved in \S\,\ref{S:period relations} below.
\end{proof}

\subsection{Proof of Theorem \ref{T:main}}

First we recall the result due to Liu \cite{Liu2019b} on the algebraicity of the critical values of the twisted standard $L$-functions for $\GSp_{2n}(\A_\F)$.

\begin{thm}[Liu]\label{T:Liu}
Let $\itPsi$ be an $C$-algebraic irreducible cuspidal automorphic representation of $\GSp_{2n}(\A_\F)$ such that $\itPsi_v$ is a holomorphic discrete series representation for each $v \in S_\infty$. 
There exists a sequence of non-zero complex numbers $(\Omega({}^\sigma\!\itPsi))_{\sigma \in {\rm Aut}(\C)}$ satisfying the following property: Let $\chi$ be a finite order Hecke character of $\A_\F^\times$ with parallel signature and $m \in \Z_{\geq 1}$ be a critical point of the twisted standard $L$-function $L(s,\itPsi , \,{\rm std}\otimes \chi)$ such that $m\neq1$ if $\F=\Q$ and $\chi^2=1$. 
For $\sigma \in {\rm Aut}(\C)$, we have
\begin{align*}
&\sigma\left(  \frac{L^S(m,\itPsi,{\rm std}\otimes\chi)}{|D_\F|^{(n+1)/2}\cdot(2\pi\sqrt{-1})^{(n+1)dm}\cdot (\sqrt{-1})^{d{\sf v}}\cdot G(\chi)^{n+1}\cdot \Omega(\itPsi)}\right)\\
& = \frac{L^S(m,{}^\sigma\!\itPsi,{\rm std}\otimes{}^\sigma\!\chi)}{|D_\F|^{(n+1)/2}\cdot(2\pi\sqrt{-1})^{(n+1)dm}\cdot (\sqrt{-1})^{d{\sf v}}\cdot  G({}^\sigma\!\chi)^{n+1}\cdot \Omega({}^\sigma\!\itPsi)}.
\end{align*}
Here ${\sf v}$ is the integer so that $|\omega_\itPsi| = |\mbox{ }|_{\A_\F}^{\sf v}$ and $S$ is a sufficiently large set of places containing $S_\infty$.
\end{thm}

\begin{rmk}
The result in \cite[Corollary 0.0.2]{Liu2019b} is stated for $\F=\Q$ with $\GSp_{2n}$ and ${\rm Aut}(\C)$ replaced by $\Sp_{2n}$ and ${\rm Aut}(\C/\Q(\zeta_N))$ for some cyclotomic field $\Q(\zeta_N)$. The generalization to arbitrary $\F$ is straightforward. By considering the Hilbert--Siegel modular varieties associated to $\GSp_{2n,\F}$ instead of the connected Shimura varieties associated to $\Sp_{2n,\F}$, we can prove the full Galois equivariance property over ${\rm Aut}(\C)$ by following the standard arguments as in \cite[\S\,3]{GL2016}. The crucial point is that Liu is able to explicitly compute the archimedean doubling local zeta integrals. 
\end{rmk}

\begin{rmk}
The other results in the literatures on the algebraicity of twisted standard $L$-functions were mainly for scalar-valued Hilbert--Siegel cusp forms (cf.\,\cite{Harris1981}, \cite{Sturm1981}, \cite{Mizumoto1991}, \cite{Shimura2000}, \cite{BS2000}). For our purpose here, we need to consider vector-valued Hilbert--Siegel cusp forms. We refer to \cite{Kozima2000} for certain vector-valued cases and \cite{PSS2020} for $n=2$ and ${\rm sgn}(\chi)=-1$. We also refer to the recent result \cite{HPSS2021} for sufficiently large critical points $m$.
\end{rmk}

We have weak automorphic descent from $\GL_7(\A_\F)$ to $\GSp_{6}(\A_\F)$ with respect to the symmetric sixth representation of $\GL_2(\C)$. More precisely, we have following result.

\begin{prop}\label{P:descent}
Assume $\itPi$ is non-CM and $\F\cap\Q(\zeta_5)=\Q$. There exists an $C$-algebraic irreducible cuspidal automorphic representation $\itPsi$ of $\GSp_{6}(\A_\F)$ satisfying the following conditions:
\begin{itemize}
\item $\itPsi_v$ is a holomorphic discrete series representation for all $v \in S_\infty$.
\item $L(s,\itPsi_v,{\rm std}\otimes\chi_v) = L(s,\itPi_v,{\rm Sym}^6 \otimes \omega_{\itPi,v}^{-3}\chi_v)$ for all $v \in S_\infty$.
\item $L^S(s,\itPsi,{\rm std}\otimes\chi) = L^S(s,\itPi,{\rm Sym}^6 \otimes \omega_{\itPi}^{-3}\chi)$ for some sufficiently large set $S$ of places containing $S_\infty$.
\end{itemize}
\end{prop}

\begin{proof}
Let ${\rm Sym}^6(\itPi)$ be the functorial lift of $\itPi$ to $\GL_7(\A_\F)$ with respect to the symmetric sixth representation of $\GL_2(\C)$. Then ${\rm Sym}^6(\itPi)$ is a regular $C$-algebraic irreducible cuspidal automorphic representation of $\GL_7(\A_\F)$. Note that the existence of the lift was proved by Clozel and Thorne \cite{CT2017} under the assumption $\F\cap\Q(\zeta_5)=\Q$. 
Since $\itPi^\vee = \itPi \otimes \omega_\itPi^{-1}$, it is easy to see that ${\rm Sym}^6(\itPi) \otimes \omega_{\itPi}^{-3}$ is self-dual with trivial central character.
Moreover, we have the factorization of twisted symmetric square $L$-function of ${\rm Sym}^6(\itPi)$:
\begin{align*}
L(s,{\rm Sym}^6(\itPi),{\rm Sym}^2 \otimes \omega_\itPi^{-3}) = L(s,\itPi,{\rm Sym}^{12}\otimes \omega_\itPi^{-6})\cdot L(s,\itPi,{\rm Sym}^{8}\otimes \omega_\itPi^{-4})\cdot L(s,\itPi,{\rm Sym}^{4}\otimes \omega_\itPi^{-2})\cdot \zeta_\F(s).
\end{align*}
Since the twisted symmetric power $L$-functions are holomorphic and non-zero at $s=1$ by \cite[Corollary 7.1.5]{BLGG2011}, we see that $L(s,{\rm Sym}^6(\itPi),{\rm Sym}^2 \otimes \omega_\itPi^{-3})$ has a pole at $s=1$.
By Arthur's multiplicity formula \cite[Theorem 1.5.2]{Arthur2013}, the automorphic representation ${\rm Sym}^6(\itPi) \otimes \omega_{\itPi}^{-3}$ descend, with respect to the standard representation ${\rm SO}_{7}(\C) \rightarrow \GL_7(\C)$, weakly to an irreducible cuspidal automorphic representation $\itPsi'$ of $\Sp_6(\A_\F)$ such that $\itPsi_v'$ is a discrete series representation for all $v \in S_\infty$. Moreover, the descent is strong at $v \in S_\infty$ since ${\rm Sym}^6(\itPi_v)$ is essentially tempered for $v \in S_\infty$. By the results of Patrikis \cite[Corollary 3.1.6 and Proposition 3.1.14]{Patrikis2019}, there exists an $C$-algebraic irreducible cuspidal automorphic representation $\itPsi$ of $\GSp_6(\A_\F)$ satisfying the following conditions:
\begin{itemize}
\item For all $v \in S_\infty$, $\itPsi_v$ is a holomorphic discrete series representation such that $\itPsi_v \vert_{\Sp_6(\F_v)}$ contains $\itPsi_v'$.
\item $\itPsi'$ is a weak functorial lift of $\itPsi$ with respect to the $L$-homomorphism ${\rm GSpin}_7(\C) \rightarrow {\rm SO}_{7}(\C)$. 
\end{itemize}
The automorphic representation $\itPsi$ then clearly satisfies the conditions we want. This completes the proof.
\end{proof}

\subsubsection{Case $\F \neq \Q$}\label{SS:3.2.1}

We begin with the case $\F\neq \Q$. We assume $\itPi$ is non-CM, $\F\cap\Q(\zeta_5)=\Q$, and $\kappa_v \geq 6$ for all $v \in S_\infty$. Let $\chi$ be a finite order Hecke character of $\A_\F^\times$ and $m \in {\rm Crit}(\itPi,{\rm Sym}^6\otimes\chi)$ a right-half critical point.
For any finite place $v$, by the Ramanujan conjecture for $\itPi$ proved by Blasius \cite{Blasius2006}, we have $L(m,\itPi_v,{\rm Sym}^6\otimes\chi_v)\neq0$. Moreover, it is easy to show that (cf.\,\cite[Proposition 3.17]{Raghuram2009})
\[
{}^\sigma\! L(s,\itPi_v,{\rm Sym}^6\otimes\chi_v) = L(s,{}^\sigma\!\itPi_v,{\rm Sym}^6\otimes{}^\sigma\!\chi_v)
\]
as rational functions in $q_v^{-s}$.
Therefore, Conjecture \ref{C:DC 2} holds for $L^{(\infty)}(m,\itPi,{\rm Sym}^6\otimes\chi)$ if and only if it holds for $L^S(m,\itPi,{\rm Sym}^6\otimes\chi)$ for any finite set of places $S$ containing $S_\infty$.
Note that $L^S(m,\itPi,{\rm Sym}^6\otimes\chi)$ is non-zero by \cite[Corollary 7.1.5]{BLGG2011}.
By Theorem \ref{T:Liu} and Proposition \ref{P:descent}, we conclude that Conjecture \ref{C:DC 2} holds for all $\chi$ and right-half critical $m$ if and only if it holds for some $\chi$ and some right-half critical $m$.
Now we verify the conjecture for 
\[
\chi = |\mbox{ }|_{\A_\F}^{3{\sf w}}\omega_\itPi^{-3},\quad m=1-3{\sf w}.
\]
Consider the Kim--Ramakrishnan--Shahidi lift $\itSigma$ of $\itPi$. By the period relations in Theorem \ref{P:period relation main}, we have
\[
\sigma \left( \frac{p^\varnothing(\itSigma)\cdot p^{S_\infty}(\itSigma^\vee)}{(2\pi\sqrt{-1})^{4d}\cdot \Vert f_\itPi \Vert^7} \right) = \frac{p^\varnothing({}^\sigma\!\itSigma)\cdot p^{S_\infty}({}^\sigma\!\itSigma^\vee)}{(2\pi\sqrt{-1})^{4d}\cdot \Vert f_{{}^\sigma\!\itPi} \Vert^7}
\]
for all $\sigma \in {\rm Aut}(\C/\F^{\Gal})$. Comparing with Theorem \ref{T:Ichino-Chen}, we conclude that
\begin{align*}
\sigma \left( \frac{L(1,\itSigma,{\rm Ad})}{\pi^{7d}\cdot \Vert f_\itPi \Vert^7} \right) = \frac{L(1,{}^\sigma\!\itSigma,{\rm Ad})}{{\pi^{7d}\cdot \Vert f_{{}^\sigma\!\itPi} \Vert^7}}
\end{align*}
for all $\sigma \in {\rm Aut}(\C/\F^{\Gal})$.
Moreover, the Galois equivariance holds for ${\rm Aut}(\C)$ if we assume further that $\kappa_v=\kappa_w$ for all $v,w \in S_\infty$.
On the other hand, we have the factorization of adjoint $L$-function:
\[
L(s,\itSigma,{\rm Ad}) = L(s,\itPi,{\rm Sym}^6\otimes\omega_\itPi^{-3})\cdot L(s,\itPi,{\rm Sym}^2\otimes\omega_\itPi^{-1}).
\]
By Deligne's conjecture for $L^{(\infty)}(1,\itPi,{\rm Sym}^2\otimes \omega_\itPi^{-1})$ proved by Sturm \cite{Sturm1980}, \cite{Sturm1989}, we have
\[
\sigma \left( \frac{L^{(\infty)}(1,\itPi,{\rm Sym}^2\otimes \omega_\itPi^{-1})}{(2\pi\sqrt{-1})^{2d+\sum_{v \in S_\infty}\kappa_v}\cdot(\sqrt{-1})^{d{\sf w}}\cdot\Vert f_\itPi\Vert} \right) = \frac{L^{(\infty)}(1,{}^\sigma\!\itPi,{\rm Sym}^2\otimes {{}^\sigma\!\omega_\itPi^{-1}})}{(2\pi\sqrt{-1})^{2d+\sum_{v \in S_\infty}\kappa_v}\cdot(\sqrt{-1})^{d{\sf w}}\cdot\Vert f_{{}^\sigma\!\itPi}\Vert}
\]
for all $\sigma \in {\rm Aut}(\C)$.
Note that
\[
L(1,\itSigma_v,{\rm Ad}) \in \pi^{-3\lambda_{1,v}+\lambda_{2,v}-3}\cdot\Q^\times = \pi^{-7\kappa_v+1}\cdot\Q^\times
\]
for all $v \in S_\infty$.
We conclude that Conjecture \ref{C:DC 2} holds for $L^{(\infty)}(1,\itPi,{\rm Sym}^6\otimes\omega_\itPi^{-3})$ with ${\rm Aut}(\C)$ replaced by ${\rm Aut}(\C/\F^{\Gal})$.
Moreover, it holds for ${\rm Aut}(\C)$ if we assume further that $\kappa_v = \kappa_w$ for all $v,w\in S_\infty$.
We thus prove Theorem \ref{T:main} for the right-half critical points. As for the left-half critical points, it follows from the global functional equation for $\GL_7(\A_\F)$ applied to the regular $C$-algebraic irreducible cuspidal automorphic representation ${\rm Sym}^6(\itPi)$.
We refer to \cite[Proposition A.4]{Chen2021} for the precise statement. 

\subsubsection{Case $\F=\Q$}\label{SS:3.2.2}
Now we prove the case when $\F=\Q$ by a base change trick. Let $\kappa {\geq 6}$ be the weight of $\itPi_\infty$. 
We assume $\itPi$ is non-CM.
Let $\L$ be a totally real cyclic extension over $\Q$ such that $r = [\L:\Q] \geq 2$ is a prime number and $\L\cap\Q(\zeta_5)=\Q$.
Let $\omega_{\mathbb{L}/\Q}$ be a non-trivial Hecke character of $\A_\Q^\times$ that vanishing on ${\rm N}_{\mathbb{L}/\Q}(\A_{\mathbb{L}}^\times)$. Note that by class field theory we have $\omega_{\mathbb{L}/\Q}^r=1$. 
Let ${\rm BC}_{\mathbb{L}}(\itPi)$ be the base change lift of $\itPi$ to $\GL_2(\A_{\mathbb{L}})$. Then ${\rm BC}_{\mathbb{L}}(\itPi)$ is regular $C$-algebraic and cuspidal with 
\[
{\rm BC}_{\mathbb{L}}(\itPi)_v = \itPi_\infty
\]
for all archimedean places $v$ of $\mathbb{L}$. 
In particular, the assumption in Theorem \ref{T:main}-(2) is satisfied by ${\rm BC}_{\mathbb{L}}(\itPi)$.
Note that $\omega_{{\rm BC}_{\mathbb{L}}(\itPi)} = \omega_\itPi\circ{\rm N}_{\mathbb{L}/\Q}$.
Let $f_{{\rm BC}_{\mathbb{L}}(\itPi)}$ be the normalized newform of ${{\rm BC}_{\mathbb{K}}(\itPi)}$.
By \cite[Theorem 4.3]{Shimura1978}, one can deduce that
\begin{align}\label{E:Main proof 5}
\sigma \left( \frac{\Vert f_{{\rm BC}_{\mathbb{L}}(\itPi)}\Vert}{\Vert f_\itPi\Vert^r}\right) = \frac{\Vert f_{{\rm BC}_{\mathbb{L}}({}^\sigma\!\itPi)}\Vert}{\Vert f_{{}^\sigma\!\itPi}\Vert^r}
\end{align}
for all $\sigma \in {\rm Aut}(\C)$.
Let $\chi$ be a finite order Hecke character of $\A_\Q^\times$. It follows easily from the definition of Gauss sum that
\begin{align}\label{E:Main proof 6}
\sigma \left( \frac{G(\chi\circ {\rm N}_{\mathbb{L}/\Q})}{G(\chi)^r} \right) =  \frac{G({}^\sigma\!\chi\circ {\rm N}_{\mathbb{L}/\Q})}{G({}^\sigma\!\chi)^r}
\end{align}
for all $\sigma \in {\rm Aut}(\C)$. We have the following factorization of twisted symmetric sixth $L$-function:
\begin{align}\label{E:Main proof 7}
\begin{split}
L(s,{\rm BC}_{\mathbb{L}}(\itPi),{\rm Sym}^6\otimes\chi\circ{\rm N}_{\mathbb{L}/\Q})  = \prod_{i=1}^r L(s,\itPi,{\rm Sym}^6\otimes\chi\omega_{\mathbb{L}/\Q}^i).
\end{split}
\end{align}
Assume $\chi$ is chosen so that ${\rm sgn}(\chi)=1$ and $\chi^2 \notin \<\omega_{\L/\Q}\>$.
Consider the critical value 
\[
L^{(\infty)}(1-3{\sf w},{\rm BC}_{\mathbb{L}}(\itPi),{\rm Sym}^6\otimes \chi\circ{\rm N}_{\L/\Q}).
\]
As we have proved in \S\,\ref{SS:3.2.1}, Conjecture \ref{C:DC 2} holds for this critical value. 
Therefore, by Theorem \ref{T:Liu} and (\ref{E:Main proof 5})-(\ref{E:Main proof 7}), we obtain the period relation
\[
\sigma\left( \frac{\Omega(\itPsi)}{(2\pi\sqrt{-1})^{6\kappa}\cdot \Vert f_\itPi\Vert^{6}}\right)^r = \left(\frac{\Omega({}^\sigma\!\itPsi)}{(2\pi\sqrt{-1})^{6\kappa}\cdot \Vert f_{{}^\sigma\!\itPi}\Vert^{6}}\right)^r
\]
for all $\sigma \in {\rm Aut}(\C)$. Here $\itPsi$ is any fixed irreducible cuspidal automorphic representation of $\GSp_6(\A_\Q)$ satisfying conditions in Proposition \ref{P:descent}.
Taking $r=2,3$, we then deduce the period relation
\[
\sigma\left( \frac{\Omega(\itPsi)}{(2\pi\sqrt{-1})^{6\kappa}\cdot \Vert f_\itPi\Vert^{6}}\right) = \frac{\Omega({}^\sigma\!\itPsi)}{(2\pi\sqrt{-1})^{6\kappa}\cdot  \Vert f_{{}^\sigma\!\itPi}\Vert^{6}}
\]
for all $\sigma \in {\rm Aut}(\C)$.
By this period relation and Theorem \ref{T:Liu} again, we conclude that Conjecture \ref{C:DC 2} holds for all right-half critical $m-3{\sf w} \in {\rm Crit}(\itPi,{\rm Sym}^6\otimes\chi)$ with $m \neq 1$ if $\chi^2=1$. Now we consider the remaining cases in the right-half critical region. Assume $\chi^2=1$ and ${\rm sgn}(\chi)=1$. Take $r=3$. Since Conjecture \ref{C:DC 2} holds for the critical values
\[
L^{(\infty)}(1-3{\sf w},{\rm BC}_{\L}(\itPi),{\rm Sym}^6\otimes\chi\circ{\rm N}_{\L/\Q}),\quad L^{(\infty)}(1-3{\sf w},\itPi,{\rm Sym}^6\otimes\chi\omega_{\mathbb{L}/\Q}),\quad L^{(\infty)}(1-3{\sf w},\itPi,{\rm Sym}^6\otimes\chi\omega_{\mathbb{L}/\Q}^2),
\]
we deduce from (\ref{E:Main proof 5})-(\ref{E:Main proof 7}) again that Conjecture \ref{C:DC 2} also holds for $L^{(\infty)}(1-3{\sf w},\itPi,{\rm Sym}^6\otimes\chi)$. 
Similarly as in \S\,\ref{SS:3.2.1}, we then conclude from the global functional equation that Conjecture \ref{C:DC 2} also holds for the left-half critical points. This completes the proof of Theorem \ref{T:main}.

\section{Algebraicity for $L$-functions on $\GSp_4 \times \GL_2$}

Let $(\itSigma,V_\itSigma)$ and $(\itPi',V_{\itPi'})$ be irreducible globally generic cuspidal automorphic representations of $\GSp_4(\A_\F)$ and $\GL_2(\A_\F)$ with central characters $\omega_\itSigma$ and $\omega_{\itPi'}$, respectively. 
Let 
\[
L(s,\itSigma \times \itPi')
\]
be the Rankin--Selberg $L$-function of $\itSigma \times \itPi'$. We recall the integral representation of the Rankin--Selberg $L$-function due to Novodvorsky \cite[\S\,3]{Novodvorsky1979} (see also \cite[Theorem 1.1]{PSS1984} and \cite[Theorem A]{LNM1254B}) in \S\,\ref{SS:integral rep.}. 
Under the assumptions that $\itSigma$ and $\itPi'$ are $C$-algebraic and satisfying conditions (\ref{E:discrete condition 1}) and (\ref{E:discrete condition 2}), in \S\,\ref{SS:coho. inter.} we cohomologically interpret the global zeta integral defining the integral representation. 
The main result of this section is Theorem \ref{T:algebraicity GSp_4 x GL_2}, where we prove the algebraicity of critical values of the Rankin--Selberg $L$-function in terms of the automorphic periods defined in \S\,\ref{S:automorphic period}.
In \S\,\ref{SS:Morimoto}, we also recall the result of Morimoto \cite{Morimoto2018} on the algebraicity for $L(s,\itSigma \times \itPi')$ and establish a period relation in Proposition \ref{P:period relation 1}.

\subsection{Integral representation of $L(s,\itSigma \times \itPi')$}\label{SS:integral rep.}

Let ${\bf G}$ be the reductive group over $\Q$ defined by
\[
{\bf G} = \{(g_1,g_2) \in \GL_{2} \times \GL_{2}\,\vert\,\nu(g_1) = \nu(g_2)\}.
\]
We regard ${\bf G}$ as a closed subgroup of $\GSp_{4}$ by the embedding
\begin{align}\label{E:embedding}
\left(\begin{pmatrix}
  a_1 & b_1 \\
  c_1 & d_1
 \end{pmatrix},
 \begin{pmatrix}
  a_2 & b_2 \\
  c_2 & d_2
 \end{pmatrix}
 \right)  \longmapsto
 \begin{pmatrix}
  a_1 & 0    & b_1 & 0 \\
  0   & a_2  & 0   & b_2 \\
  c_1 & 0    & d_1 & 0 \\
  0   & c_2 & 0   & d_2 
 \end{pmatrix}.
\end{align}
Let $N_{\bf G}$ be the maximal unipotent subgroup of $\bf G$ consisting of upper unipotent matrices.
By abuse of notation, we write ${\bf G}$ for the base change ${\bf G}_\F = \G\times_\Q\F$ to algebraic group over $\F$.
Let $\psi = \bigotimes_v\psi_v$ be a non-trivial additive character of $\F\backslash\A_\F$. 
For a finite place $v$ such that $\itPi_v'$ and $\itSigma_v$ are unramified and $\psi_v$ has conductor $\o_{v}$, let 
\[
W_{v,\psi_{v}}^{\circ} \in \mathcal{W}(\itPi_v',\psi_{v}),\quad W_{v,\psi_{U,v}}^{\circ} \in \mathcal{W}(\itSigma_v,\psi_{U,v})
\]
be the $\GL_2(\o_v)$-invariant and $\GSp_4(\o_v)$-invariant Whittaker functions, respectively, normalized so that $W_{v,\psi_v}^\circ(1) = W_{v,\psi_{U,v}}^{\circ}(1)=1$.
For $\varphi' \in V_{\itPi'}$, $\varphi \in V_\itSigma$, and holomorphic section $f^{(s)}$ of $I(\omega_\itSigma^{-1}\omega_{\itPi'}^{-1},s)$, we define the global zeta integral
\[
Z(E(f^{(s)}),\,\varphi',\,\varphi) = \int_{\A_\F^\times\G(\F)\backslash\G(\A_\F)}(E(f^{(s)}) \otimes \varphi')(g)\varphi(g)\,dg^{\rm Tam}.
\]
Here $dg^{\rm Tam}$ is the Tamagawa measure on $\A_\F^\times\backslash\G(\A_\F)$.
By an unfolding argument (cf.\,\cite[\S\,1.1]{LNM1254B}), we have
\[
Z(E(f^{(s)}),\,\varphi',\,\varphi) = \int_{\A_\F^\times N_{\G}(\A_\F)\backslash\G(\A_\F)}(f^{(s)} \otimes W_{\varphi',\psi})(g)W_{\varphi,\psi_U}(g)\,d\overline{g}^{\rm Tam}
\]
for ${\rm Re}(s)$ sufficiently large.
Here $d\overline{g}^{\rm Tam}$ is the quotient of $dg^{\rm Tam}$ by the Tamagawa measure on $N_\G(\A_\F)$.
For each place $v$, $W_v' \in \mathcal{W}(\itPi_v',\psi_v)$, $W_v \in \mathcal{W}(\itSigma_v,\psi_v)$, and $f_v^{(s)}$ a meromorphic section of $I(\omega_{\itSigma,v}^{-1}\omega_{\itPi',v}^{-1},s)$, let $Z(f_v^{(s)},W_v',W_v)$ be the local zeta integral defined by
\begin{align}\label{E:local zeta}
Z(f_v^{(s)},W_v',W_v) = \int_{\F_v^\times N_{\G}(\F_v)\backslash \G(\F_v)}(f_v^{(s)}\otimes W_v')(g_v)W_v(g_v)\,d\overline{g}_v.
\end{align}
Here $d\overline{g}_v$ is the quotient measure on $\F_v^\times N_{\G}(\F_v)\backslash \G(\F_v)$ defined as follows: 
\begin{itemize}
\item If $v$ is finite, then $d\overline{g}_v$ is the quotient of the Haar measures on $\F_v^\times\backslash\G(\F_v)$ and $N_\G(\F_v)$ with \[{\rm vol}(\o_v^\times\backslash\G(\o_v)) = {\rm vol}(N_\G(\o_v))=1.\]
\item If $v \in S_\infty$, then
\[
d\overline{g}_v = |a_{1,v}a_{2,v}|_v^{-3}da_{1,v}\,da_{2,v}\,dk_{1,v}\,dk_{2,v}
\]
for
\[
g_v = \left({\rm diag}(a_{1,v}a_{2,v},a_{2,v}^{-1}) k_{1,v}, {\rm diag}(a_{1,v},1) k_{2,v}\right)
\]
with $a_{1,v},a_{2,v}\in\F_v^\times$ and $k_{1,v},k_{2,v} \in {\rm SO}(2)$, where $da_{i,v}$ is the Lebesgue measure and ${\rm Vol}({\rm SO}(2),dk_{i,v})=1$ for $i=1,2$.
\end{itemize}
Note that we have (cf.\,\cite[\S\,6.1]{IP2018})
\[
d\overline{g}^{\rm Tam} = |D_\F|^{-2}\zeta_{\F}(2)^{-2}\cdot \prod_v d\overline{g}_v.
\]
Now we recall the integral representation.
Let $S$ be a finite set of places of $\F$ containing $S_\infty$ so that for $v \notin S$, $\itPi_v'$ and $\itSigma_v$ are unramified and $\psi_v$ has conductor $\o_{v}$.
We write 
\[
\itPi'_S = \bigotimes_{v \in S}\itPi'_v\,\quad\itSigma_{S} = \bigotimes_{v \in S}\itSigma_{v},\quad  \psi_S = \bigotimes_{v \in S}\psi_v.
\]
By \cite[Theorem 3.1]{PSR1987}, we have
\[
Z(f_{v,\circ}^{(s)},\,W_{v,\psi_v}^\circ,\,W_{v,\psi_{U,v}}^\circ) = L(s,\itSigma_v\times\itPi_v')
\]
for all $v \notin S$. 
Thus we obtain the following

\begin{prop}[Novodvorsky]\label{P:integral rep.}
Let $W_S' \in \mathcal{W}(\itPi_S',\psi_S)$, $W_S \in \mathcal{W}(\itSigma_S,\psi_{U,S})$, and $f_S^{(s)}$ a meromorphic section of $I(\omega_{\itSigma,S}^{-1}\omega_{\itPi',S}^{-1},s)$.
Let $\varphi' \in V_{\itPi'}$, $\varphi \in V_\itSigma$, and $f^{(s)}$ the meromorphic section of $I(\omega_{\itSigma}^{-1}\omega_{\itPi'}^{-1},s)$ defined by
\[
W_{\varphi',\psi} = \prod_{v \notin S}W_{v,\psi_v}^\circ \cdot W_S',\quad W_{\varphi,\psi_U} = \prod_{v \notin S}W_{v,\psi_{U,v}}^\circ \cdot W_S,\quad f^{(s)} = \bigotimes_{v \notin S} f_{v,\circ}^{(s)}\otimes f_S^{(s)}.
\]
Then we have
\[
Z(E(f^{(s)}),\,\varphi',\,\varphi) = |D_\F|^{-2}\zeta_{\F}(2)^{-2}\cdot L^S(s,\itSigma\times\itPi')\cdot Z(f_S^{(s)},W_S',W_S).
\]
\end{prop}

\subsection{Cohomological interpretation of global zeta integral}\label{SS:coho. inter.}

Let $K = \R_+ \times {\rm U}(2) \times {\rm U}(2)$ and regard it as a closed subgroup of ${\bf G}(\R)$ by the homomorphism
\[
(a,a_1+\sqrt{-1}\,b_1,a_1+\sqrt{-1}\,b_1 )\longmapsto \left( a \bp a_1 & b_1 \\ -b_1 & a_1\ep, a\bp a_2 & b_2 \\ -b_2 & a_2\ep\right).
\]
We write $\frak{g}\subset\frak{g}_1 \oplus \frak{g}_1$ for the Lie algebra of ${\bf G}(\R) \subset \GL_2(\R) \times \GL_2(\R)$. Define $\frak{k} \subset \frak{g}$ and $\frak{p}^\pm,\frak{P} \subset \frak{g}_\C$ by
\[
\frak{k} = {\rm Lie}(K)=\R\cdot ({\bf 1}_2 \oplus {\bf 1}_2)\oplus \left( \frak{k}_1 \oplus \frak{k}_1\right),\quad \frak{p}^\pm = \frak{p}_1^\pm \oplus \frak{p}_1^\pm,\quad \frak{P} = \frak{k}_\C \oplus \frak{p}^-.
\]
For $({\sf v};\,\kappa_1,\kappa_2)\in\Z\times\Z\times\Z$ with $\kappa_1+\kappa_2 \equiv {\sf v}\,({\rm mod}\,2)$, let $(\rho_{({\sf v};\,\kappa_1,\kappa_2)},V_{({\sf v};\,\kappa_1,\kappa_2)})$ be the algebraic character of $K$ defined by $V_{({\sf v};\,\kappa_1,\kappa_2)}=\C$ and 
\[
\rho_{({\sf v};\,\kappa_1,\kappa_2)}(a(u_1,u_2))\cdot z = a^{\sf v}u_1^{\kappa_1}u_2^{\kappa_2}\cdot z
\]
for $a \in \R^\times$ and $u_1,u_2 \in {\rm U}(1)$.
Let $(G,X)$ be the Shimura datum defined by
\[
G = {\rm Res}_{\F/\Q}{\bf G}_\F,
\]
and $X$ is the $G(\R)$-conjugacy class containing the morphism $h : \mathbb{S}\rightarrow G_\R$ with
\[
h(x+\sqrt{-1}\,y) = \left( \bp x & y \\ -y & x\ep,\cdots,\bp x & y \\ -y & x\ep\right) \times \left( \bp x & y \\ -y & x\ep,\cdots,\bp x & y \\ -y & x\ep\right)
\]
on $\R$-points.
Under the identification of $\F_v$ with $\R$ for each $v \in S_\infty$, we have
\[
K_{h} = K^{S_\infty},\quad \frak{k}_{h} = \frak{k}^{S_\infty},\quad \frak{p}_{h}^\pm = (\frak{p}^\pm)^{S_\infty},\quad \frak{P}_{h} = \frak{P}^{S_\infty}.
\]
For $({\sf v};\,\underline{\kappa}_1,\underline{\kappa}_2) \in \Z\times\Z^{S_\infty}\times\Z^{S_\infty}$ with $\underline{\kappa}_1 = (\kappa_{1,v})_{v \in S_\infty}$, $\underline{\kappa}_2 = (\kappa_{2,v})_{v \in S_\infty}$, and $\kappa_{1,v}+\kappa_{2,v} \equiv {\sf v}\,({\rm mod}\,2)$, let $(\rho_{({\sf v};\,\underline{\kappa}_1,\underline{\kappa}_2)},V_{({\sf v};\,\underline{\kappa}_1,\underline{\kappa}_2)})$ be the motivic algebraic character of $K^{S_\infty}$ defined by
\[
(\rho_{({\sf v};\,\underline{\kappa}_1,\underline{\kappa}_2)},V_{({\sf v};\,\underline{\kappa}_1,\underline{\kappa}_2)}) = \left( \bigotimes_{v \in S_\infty}\rho_{({\sf v};\,\underline{\kappa}_{1,v},\underline{\kappa}_{2,v})},\bigotimes_{v \in S_\infty}V_{({\sf v};\,\underline{\kappa}_{1,v},\underline{\kappa}_{2,v})}\right).
\]
Let 
\[
[\mathcal{V}_{({\sf v};\,\underline{\kappa}_1,\underline{\kappa}_2)}]
\]
be the automorphic line bundle on ${\rm Sh}(G,X)$ defined by $(\rho_{({\sf v};\,\underline{\kappa}_1,\underline{\kappa}_2)},V_{({\sf v};\,\underline{\kappa}_1,\underline{\kappa}_2)})$.

We assume $\itSigma$ and $\itPi'$ are $C$-algebraic and satisfying conditions (\ref{E:discrete condition 1}) and (\ref{E:discrete condition 2}). 
Let $(\underline{\lambda};\,{\sf u})$ and $(\underline{\ell};\,{\sf w}')$ be the weights of $\itSigma$ and $\itPi'$, respectively.
A critical point for $L(s,\itSigma \times \itPi')$ is an integer $m$ which is not a pole of the archimedean local factors $L(s,\itSigma_v \times \itPi_v')$ and $L(1-s,\itSigma_v^\vee \times (\itPi_v')^\vee)$ for all $v \in S_\infty$.
Note that we have
\begin{align*}
L(s-\tfrac{{\sf u}+{\sf w}'}{2},\itSigma_v \times \itPi_v') &= \Gamma_\C\left(s+\tfrac{|\lambda_{1,v}+\lambda_{2,v}-\ell_v|}{2}\right)\Gamma_\C\left(s+\tfrac{|\lambda_{1,v}-\lambda_{2,v}-\ell_v|}{2}\right)\\
&\times \Gamma_\C\left(s+\tfrac{\lambda_{1,v}+\lambda_{2,v}+\ell_v-2}{2}\right)\Gamma_\C\left(s+\tfrac{\lambda_{1,v}-\lambda_{2,v}+\ell_v-2}{2}\right).
\end{align*}
If we assume further that $1 \leq \ell_v \leq \lambda_{1,v}$
for all $v \in S_\infty$, then the set of critical points is given by
\[
{\rm Crit}(\itSigma \times \itPi') = \left\{m \in \Z \, \left\vert \, 1-\tfrac{|\lambda_{1,v}+\lambda_{2,v}-\ell|}{2} \leq m+\tfrac{{\sf u}+{\sf w}'}{2}\leq \tfrac{|\lambda_{1,v}+\lambda_{2,v}-\ell|}{2} \right\}\right..
\]
Put
\[
\kappa = \min_{v \in S_\infty}\left\{|\lambda_{1,v}+\lambda_{2,v}-\ell_v|\right\}.
\]
Let $I,J$ be subsets of $S_\infty$ defined by
\begin{align}\label{E:period set}
I  = \left\{v \in S_\infty\,\vert\,1 \leq \ell_v \leq \lambda_{1,v}+\lambda_{2,v}-1\right\},\quad
J  = \left\{v \in S_\infty\,\vert\,\lambda_{1,v}+\lambda_{2,v}+1 \leq \ell_v \leq \lambda_{1,v}\right\}.
\end{align}
It is clear that $I$ and $J$ are admissible with respect to $\underline{\lambda}$ and $\underline{\ell}$, respectively.
Note that $I\cup J=S_\infty$ if and only if ${\rm Crit}(\itSigma \times \itPi')$ is non-empty. 
In the following proposition, we prove Galois equivariant property of the global zeta integral in terms of the automorphic periods.

\begin{prop}\label{P:Galois equiv. global}
Assume the following conditions are satisfied:
\begin{itemize}
\item $I\cup J=S_\infty$.
\item $\kappa = |\lambda_{1,v}+\lambda_{2,v}-\ell_v|$ for all $v \in S_\infty$. 
\item If $\F=\Q$, then $\kappa \neq 2$ or $\omega_\itSigma\omega_{\itPi'} \neq |\mbox{ }|_{\A_\F}^{{\sf u}+{\sf w}'}$.
\end{itemize}
Let $\varphi' \in V_{\itPi'}^+$, $\varphi \in V_\itSigma^+$, and $f^{(s)}$ be a meromorphic section of $I(\omega_{\itSigma,f}^{-1}\omega_{\itPi',f}^{-1},s)$ satisfying conditions in Proposition \ref{P:Eisenstein series}. For $\sigma \in {\rm Aut}(\C)$, we have
\begin{align*}
&\sigma \left( \frac{Z\left( E^{[\kappa]}(f^{(s)}),\,(\varphi')^J,\,{\rm pr}_{\underline{\ell}-\underline{\lambda}_2}(\varphi)^I\right)}{|D_\F|^{1/2}\cdot(2\pi\sqrt{-1})^{d(\kappa+{\sf u}+{\sf w}')/2}\cdot G(\omega_\itSigma\omega_{\itPi'})\cdot p^I(\itSigma)\cdot p^J(\itPi')}\right)\\
& =  \frac{Z\left(E^{[\kappa]}({}^\sigma\!f^{(s)}),\,({}^\sigma\!\varphi')^{{}^\sigma\!J},\, {\rm pr}_{{}^\sigma\!\underline{\ell}-{}^\sigma\!\underline{\lambda}_2}({}^\sigma\!\varphi)^{{}^\sigma\!I}\right)}{|D_\F|^{1/2}\cdot(2\pi\sqrt{-1})^{d(\kappa+{\sf u}+{\sf w}')/2}\cdot G({}^\sigma\!\omega_\itSigma{}^\sigma\!\omega_{\itPi'})\cdot p^{{}^\sigma\!I}({}^\sigma\!\itSigma)\cdot p^{{}^\sigma\!J}({}^\sigma\!\itPi')}.
\end{align*}
\end{prop}

\begin{proof}
The natural inclusion ${\bf G} \subset \GL_2 \times \GL_2$ induces a morphism $(G,X)\rightarrow (G_1\times G_1,X_1\times X_1)$ of Shimura data. Note that
\[
(\rho_{(-\underline{\kappa};\,{\sf u}+{\sf w})}\boxtimes\rho_{(\underline{\ell}(J);\,-{\sf w})})\vert_{K^{S_\infty}} = \rho_{({\sf u};\,-\underline{\kappa},\,\underline{\ell}(J))}.
\]
Then we have a canonical homomorphism
\[
H^{{}^\sharp J}([\mathcal{V}_{(-\underline{\kappa};\,{\sf u}+{\sf w})}]^{\rm can}\times[\mathcal{V}_{(\underline{\ell}(J);\,-{\sf w})}]^{\rm can}) \longrightarrow H^{{}^\sharp J}([\mathcal{V}_{({\sf u};\,-\underline{\kappa},\,\underline{\ell}(J))}]^{\rm can}).
\]
By the K\"unneth formula, we have a canonical isomorphism
\[
H^{{}^\sharp J}([\mathcal{V}_{(-\underline{\kappa};\,{\sf u}+{\sf w}})]^{\rm can}\times[\mathcal{V}_{(\underline{\ell}(J);\,-{\sf w})}]^{\rm can}) \simeq \bigoplus_{q_1+q_2 = {}^\sharp J}H^{q_1}([\mathcal{V}_{(-\underline{\kappa};\,{\sf u}+{\sf w}})]^{\rm can})\otimes H^{q_2}([\mathcal{V}_{(\underline{\ell}(J);\,-{\sf w})}]^{\rm can}).
\]
In particular, we obtain a homomorphism
\[
F_1 : H^0([\mathcal{V}_{(-\underline{\kappa};\,{\sf u}+{\sf w}})]^{\rm can})\otimes H^{{}^\sharp J}([\mathcal{V}_{(\underline{\ell}(J);\,-{\sf w})}]^{\rm can})\longrightarrow H^{{}^\sharp J}([\mathcal{V}_{({\sf u};\,-\underline{\kappa},\,\underline{\ell}(J))}]^{\rm can}).
\]
For $\varphi' \in V_{\itPi'}^+$ and meromorphic section $f^{(s)}$ of $I(\omega_{\itSigma,f}^{-1}\omega_{\itPi',f}^{-1},s)$ satisfying conditions in Proposition \ref{P:Eisenstein series}, it is clear that $F_1([E^{[\kappa]}(f^{(s)})] \otimes [\varphi']_J)$ is the class represented by
\begin{align}\label{Galois equiv. global proof 1}
\begin{split}
&\left( E^{[\kappa]}(f^{(s)})\otimes (\varphi')^J \right) \otimes \text{\LARGE$\wedge$}_{v \in J}(0\oplus X_{+,v}) \otimes \bigotimes_{v \in S_\infty}1\\
& \in \mathcal{A}({\bf G}(\A_\F))\otimes \exterior{{}^\sharp J}(\frak{p}^+)^{S_\infty}\otimes V_{({\sf u};\,-\underline{\kappa},\,\underline{\ell}(J))}.
\end{split}
\end{align}
For $\sigma \in {\rm Aut}(\C)$, in a similarly way we obtain a homomorphism
\[
F_{1,\sigma}:H^0([\mathcal{V}_{(-\underline{\kappa};\,{\sf u}+{\sf w}})]^{\rm can})\otimes H^{{}^\sharp J}([\mathcal{V}_{({}^\sigma\!\underline{\ell}({}^\sigma\!J);\,-{\sf w})}]^{\rm can})\longrightarrow H^{{}^\sharp J}([\mathcal{V}_{({\sf u};\,-\underline{\kappa},\,{}^\sigma\!\underline{\ell}({}^\sigma\!J))}]^{\rm can})
\]
satisfying the Galois equivariant property:
\begin{align}\label{Galois equiv. global proof 2}
T_\sigma\circ F_1 = F_{1,\sigma}\circ(T_\sigma\otimes T_\sigma).
\end{align}
The embedding (\ref{E:embedding}) induces a morphism $(G,X) \rightarrow (G_2,X_2)$ of Shimura data. 
Note that
\[
\rho_{(\underline{\lambda}(I);\,-{\sf u})}\vert_{K^{S_\infty}} = \bigoplus_{\underline{i}}\left(\bigotimes_{v \in I} \rho_{(-{\sf u};\,{\lambda}_{1,v}-{3}-{i}_v,\,{\lambda}_{2,v}-1+{i}_v)}\right) \otimes \left(\bigotimes_{v \in J} \rho_{(-{\sf u};\,-{\lambda}_{2,v}-{2}-{i}_v,\,-{\lambda}_{1,v}+{i}_v)}\right),
\]
where $\underline{i} = (i_v)_{v \in S_\infty} \in \Z^{ S_\infty}$ with $0 \leq i_v\leq \lambda_{1,v}-\lambda_{2,v}-2$. 
The assumption on $\kappa$ implies that $\rho_{(-{\sf u};\,\underline{\kappa-2},\,\underline{\ell}(I))}$ appears in the right-hand side of the decomposition with 
\[
i_v = \begin{cases}
-\lambda_{2,v}+\ell_v-1 & \mbox{ if $v \in I$},\\
\lambda_{1,v}-\ell_v & \mbox{ if $v \in J$}.
\end{cases}
\]
Consider the projection
\[
\rho_{(\underline{\lambda}(I);\,-{\sf u})} \longrightarrow \rho_{(-{\sf u};\,\underline{\kappa-2},\,\underline{\ell}(I))}.
\]
Then we have a canonical homomorphism
\[
F_2 : H^{d+{}^\sharp I}([\mathcal{V}_{(\underline{\lambda}(I);\,-{\sf u})}]^{\rm sub}) \longrightarrow H^{d+{}^\sharp I}([\mathcal{V}_{(-{\sf u};\,\underline{\kappa-2},\,\underline{\ell}(I))}]^{\rm sub}).
\]
Let $J_0 \subset J$ be the subset consisting of $v \in J$ such that $\ell_v \leq \lambda_{1,v}+2$. For each $J' \subset J_0$, define $\underline{i}(J') \in \Z^{S_\infty}$ defined by 
\[
i_v(J') = \begin{cases}
\ell_v - \lambda_{2,v}-2 & \mbox{ if $v \in J'$},\\
\ell_v - \lambda_{2,v}   & \mbox{ if $v \notin J'$}.
\end{cases}
\]
For $\varphi \in V_\itSigma^+$, by the definition of $[\varphi]_I$ in Lemma \ref{L:Eichler-Shimura 2}, $F_2([\varphi]_I)$ is represented by
\begin{align}\label{Galois equiv. global proof 3}
\begin{split}
&\sum_{J' \subset J_0}\left[C_{J'}\cdot{\rm pr}_{\underline{i}(J')}(\varphi)^I\vert_{{\bf G}(\A_\F)} \right.\\
&\left.\quad\quad\otimes \left(\text{\LARGE$\wedge$}_{v \in I}(X_{+,v}\oplus 0)\wedge(0\oplus X_{+,v})\right)\wedge\left(\text{\LARGE$\wedge$}_{v \in J\smallsetminus J'}(X_{+,v}\oplus 0)\right) \wedge\left(\text{\LARGE$\wedge$}_{v \in J'}(0\oplus X_{+,v})\right) \otimes \bigotimes_{v \in S_\infty} 1\right]\\
& \in C_{\rm rda}^\infty({\bf G}(\A_\F)) \otimes \exterior{d+{}^\sharp I}(\frak{p}^+)^{S_\infty} \otimes V_{(-{\sf u};\,\underline{\kappa-2},\,\underline{\ell}(I))}
\end{split}
\end{align}
for some $C_{J'}\in\Q$ independent of $\varphi$.
Here $C_{\rm rda}^\infty({\bf G}(\A_\F))$ is the space of smooth functions $f$ on ${\bf G}(\F)\backslash {\bf G}(\A_\F)$ such that $X\cdot f$ is rapidly decreasing for all $X \in \frak{g}_\C^{S_\infty}$.
Note that $\underline{i}(\varnothing) = \underline{\ell}-\underline{\lambda}_2$ and by Lemma \ref{L:auxiliary} we have 
\begin{align}\label{Galois equiv. global proof 5}
C_\varnothing \neq 0.
\end{align}
For $\sigma \in {\rm Aut}(\C)$, in a similarly way we obtain a homomorphism
\[
F_{2,\sigma} : H^{d+{}^\sharp I}([\mathcal{V}_{({}^\sigma\!\underline{\lambda}(I);\,-{\sf u})}]^{\rm sub}) \longrightarrow H^{d+{}^\sharp I}([\mathcal{V}_{(-{\sf u};\,\underline{\kappa-2},\,{}^\sigma\!\underline{\ell}({}^\sigma\!I))}]^{\rm sub})
\]
satisfying the Galois equivariant property:
\begin{align}\label{Galois equiv. global proof 4}
T_\sigma\circ F_2 = F_{2,\sigma}\circ T_\sigma.
\end{align}
Now we consider  the Serre duality pairing
\[
H^{{}^\sharp J}([\mathcal{V}_{({\sf u};\,-\underline{\kappa},\,\underline{\ell}(J))}]^{\rm can})\times H^{d+{}^\sharp I}([\mathcal{V}_{(-{\sf u};\,\underline{\kappa-2},\,\underline{\ell}(I))}]^{\rm sub}) \longrightarrow \C,\quad (c_1,c_2)\longmapsto \int_{{\rm Sh}(G,X)}c_1\wedge c_2
\]
in \cite[Proposition 3.8 and Remark 3.8.4]{Harris1990}.
Similarly as in the proof of Theorem \ref{T:Ichino-Chen}, by (\ref{Galois equiv. global proof 1}), (\ref{Galois equiv. global proof 3}), and (\ref{Galois equiv. global proof 5}) we have
\[
\int_{{\rm Sh}(G,X)}F_1([E^{[\kappa]}(f^{(s)})] \otimes [\varphi']_J)\wedge F_2([\varphi]_I) = C'\cdot Z\left( E^{[\kappa]}(f^{(s)}),\,(\varphi')^J,\,{\rm pr}_{\underline{\ell}-\underline{\lambda}_2}(\varphi)^I\right)
\] 
for some $C'\in\Q^\times$ independent of $\varphi',\varphi,f^{(s)}$.
For $\sigma \in {\rm Aut}(\C)$, we also have
\[
\int_{{\rm Sh}(G,X)}F_{1,\sigma}([E^{[\kappa]}({}^\sigma\!f^{(s)})] \otimes [{}^\sigma\!\varphi']_{{}^\sigma\!J})\wedge  F_{2,\sigma}([{}^\sigma\!\varphi]_{{}^\sigma\!I}) = C'\cdot Z\left(E^{[\kappa]}({}^\sigma\!f^{(s)}),\,({}^\sigma\!\varphi')^{{}^\sigma\!J},\, {\rm pr}_{{}^\sigma\!\underline{\ell}-{}^\sigma\!\underline{\lambda}_2}({}^\sigma\!\varphi)^{{}^\sigma\!I}\right).
\] 
By Proposition \ref{P:Eisenstein series}, (\ref{Galois equiv. global proof 2}), \ref{Galois equiv. global proof 4}), for $\sigma \in {\rm Aut}(\C)$ we have
\begin{align*}
&T_\sigma\left(\frac{F_1([E^{[\kappa]}(f^{(s)})] \otimes [\varphi']_J)}{|D_\F|^{1/2}\cdot(2\pi\sqrt{-1})^{d(\kappa+{\sf u}+{\sf w}')/2}\cdot G(\omega_\itSigma\omega_{\itPi'})\cdot p^J(\itPi')}\right)\\
&= \frac{F_{1,\sigma}([E^{[\kappa]}({}^\sigma\!f^{(s)})] \otimes [{}^\sigma\!\varphi']_{{}^\sigma\!J})}{|D_\F|^{1/2}\cdot(2\pi\sqrt{-1})^{d(\kappa+{\sf u}+{\sf w}')/2}\cdot G({}^\sigma\!\omega_\itSigma{}^\sigma\!\omega_{\itPi'})\cdot p^{{}^\sigma\!J}({}^\sigma\!\itPi')}
\end{align*}
and
\begin{align*}
T_\sigma \left(  \frac{F_2([\varphi]_I)}{p^I(\itSigma)}   \right) & = \frac{F_{2,\sigma}([{}^\sigma\!\varphi]_{{}^\sigma\!I})}{p^{{}^\sigma\!I}({}^\sigma\!\itSigma)}.
\end{align*}
The assertion then follows from the Galois equivarient property of the Serre duality pairing.
This completes the proof.
\end{proof}

\subsection{Local zeta integrals}

Let $v$ be a place of $\F$.
In this section, we prove Galois equivariant property of the non-archimedean local zeta integrals and compute archimedean local zeta integral.

\begin{lemma}\label{L:non-archi. local zeta}
Assume $v$ is finite and lying above a rational prime $p$.
Let $W_v' \in \mathcal{W}(\itPi_v',\psi_v)$, $W_v \in \mathcal{W}(\itSigma_v,\psi_{U,v})$, and $f_v^{(s)}$ be a rational section of $I(\omega_{\itSigma,v}^{-1}\omega_{\itPi',v}^{-1},s)$. Then $Z(f_v^{(s)},W_v',W_v)$ 
is a rational function in $q_v^{-s}$, and we have
\[
{}^\sigma\!Z(f_v^{(s)},W_v',W_v) = {}^\sigma\!\omega_{\itPi',v}(u_{\sigma,p})^{-1}\cdot Z({}^\sigma\!f_v^{(s)},t_{\sigma,v}W_v',t_{\sigma,v}W_v)
\]
for all $\sigma \in {\rm Aut}(\C)$. Here $u_{\sigma,p} \in \Z_p^\times\subset \o_v^\times$ is the unique element such that $\sigma(\psi_v(x_v)) = \psi_v(u_{\sigma,p}x_v)$ for all $x_v \in \o_v$.
Moreover, there exists a triplet $(f_v^{(s)},W_v',W_v)$ such that $f_v^{(s)}$ is holomorphic for ${\rm Re}(s)>-\tfrac{{\sf u}+{\sf w}'}{2}$ and $Z(f_v^{(s)},W_v',W_v) = 1$.
\end{lemma}

\begin{proof}
We drop the subscript $v$ for brevity.
We have
\begin{align*}
Z(f^{(s)},W',W) &= \int_{{\bf G}(\o )}dk \int_{\F ^\times}\frac{d^\times a_{1 }}{|a_{1 }| ^2}\int_{\F ^\times}\frac{d^\times a_{2 }}{|a_{2 }| ^2}\, (f ^{(s)}\otimes W ')\left({\rm diag}(a_{1 }a_{2 },a_{1 },a_{2 }^{-1},1)k \right)W \left({\rm diag}(a_{1 }a_{2 },a_{1 },a_{2 }^{-1},1)k \right)\\
& = \int_{{\bf G}(\o )} f ^{(s)}(k_{1 })I_s(W ,W ';k )\,dk ,
\end{align*}
where ${\rm vol}({\bf G}(\o ),dk )=1$ and
\begin{align*}
I_s(W' ,W ;k ) &= \int_{\F ^\times}d^\times a_{1 }\int_{\F ^\times}d^\times a_{2 }\,W'({\rm diag}(a_{1 },1)k_{2 })W \left({\rm diag}(a_{1 }a_{2 },a_{1 },a_{2 }^{-1},1)k \right)|a_{1 }| ^{s-2}|a_{2 }| ^{2s-2}\omega_{\itSigma }\omega_{\itPi '}(a_{2 })
\end{align*}
for $k =(k_{1 },k_{2 }) \in {\bf G}(\o )$.
By the asymptotic behavior of the Whittaker functions (cf.\,\cite[Theorem 3.1]{LM2009}), we see that $I_s(W' ,W ;k )$ converges absolutely for ${\rm Re}(s)$ sufficiently large.
Moreover, by \cite[Lemma 6.2]{Chen2021c}, the integral is actually a generalized Tate integral as in \cite[Proposition A]{Grobner2018}. In particular, this implies that $I_s(W' ,W ;k )$ defines a rational function in $q ^{-s}$ and 
\[
{}^\sigma\!I_s(W' ,W ;k ) = I_s(\sigma W' ,\sigma W ;k )
\]
for all $\sigma \in {\rm Aut}(\C)$. 
Making a change of variables from $(a_{1 },a_{2 })$ to $(u_{\sigma,p}^{-1}a_{1 },u_{\sigma,p}^{-1}a_{2 })$, we see that
\[
I_s(\sigma W' ,\sigma W ;k ) = {}^\sigma\!\omega_{\itPi' }(u_{\sigma,p})^{-1}\cdot I_s(t_{\sigma } W' ,t_{\sigma }W ;k ).
\]
Hence we see that $Z(W' ,W ,f ^{(s)})$ is a rational function in $q ^{-s}$ and 
\begin{align*}
{}^\sigma\! Z(f ^{(s)}, W' ,W)&= Z({}^\sigma\!f ^{(s)}, \sigma W' ,\sigma W)\\
& = \omega_{\itPi' }(u_{\sigma,p})^{-1}\cdot Z({}^\sigma\!f ^{(s)}, t_{\sigma } W' ,t_{\sigma } W)
\end{align*}
for all $\sigma \in {\rm Aut}(\C)$. 
Now we verify the second assertion. 
Without lose of generality, we assume $\psi$ has conductor $\o$.
Let $W'$ be the Whittaker function of $\itPi'$ such that 
\[
W'({\rm diag}(a,1)) = \mathbb{I}_{\o^\times}(a).
\] 
Let $N_2$ be the upper unipotent matrices of $\GL_2$ and $K(n)\subset K_1(n)\subset K_0(n)$ be open compact subgroups of $\GL_2(\o)$ defined by
\begin{align*}
K(n) & = \left\{k \in \GL_2(\o)\,\left\vert\,k \equiv 1 \,({\rm mod}\,\varpi^n) \right\}\right.,\quad
K_1(n) =  \left\{k \in \GL_2(\o)\,\left\vert\,k \equiv \bp * & * \\ 0 & 1\ep \,({\rm mod}\,\varpi^n) \right\}\right.,\\\quad K_0(n) &= \left\{k \in \GL_2(\o)\,\left\vert\,k \equiv \bp * & * \\ 0 & *\ep \,({\rm mod}\,\varpi^n) \right\}\right..
\end{align*}
Let $P_3$ and $Q$ be the mirabolic subgroup and Klingen subgroup of $\GL_3$ and $\GSp_4$, respectively, defined by
\[
P_3 = \left \{ \bp * & * &* \\ *&*&* \\ 0&0&1\ep \in \GL_3 \right\},\quad Q = \left\{ \bp * & * & *&* \\ 0&*&*&*\\0&0&*&0\\0&*& *&*  \ep\in \GSp_4  \right\}.
\]
We have a homomorphism $i : Q \rightarrow P_3$ defined by
\[
i\left(\bp t & 0 & 0&0 \\ 0&a&0&b\\0&0&\nu t^{-1}&0\\0&c&0&d  \ep\bp 1 & x & y&z \\0&1&z&0\\0&0&1&0\\0&0&-x&1 \ep\right) =  \bp g\cdot t\nu^{-1}& 0 \\ 0 & 1\ep\bp 1 & 0 & z \\ 0 & 1 & -x \\ 0 & 0 & 1\ep,
\]
where $g = \bp a&b\\c&d \ep$. Let $A_3$ be the diagonal matrices of $P_3$, $N_3$ the unipotent radical of $P_3$, and $\psi_{N_3}$ the character of $N_3$ defined by $\psi_{N_3}(u) = \psi(-u_{12}+u_{23})$ for $u=(u_{ij})$.
Consider the compact induction ${\rm c}\mbox{-}{\rm ind}_{N_3(\F)}^{P_3(\F)}(\psi_{N_3})$. 
Since $\itSigma$ is generic, by \cite[Proposition 2.35]{BZ1976} and \cite[Theorem 2.5.3-(ii)]{RS2007}, every function in $\rm c$-${\rm ind}_{N_3(\F)}^{P_3(\F)}(\psi_{N_3})$ can be lifted to a Whittaker function of $\itSigma$ with respect to $\psi_U$ through the homomorphism $i$.
We refer to the proof of \cite[Proposition 2.6.4]{RS2007} for similar arguments.
Now we take a sufficiently large $n$ such that the following conditions are satisfied:
\begin{itemize}
\item $W'$ is right $K(n)$-invariant.
\item $\omega_\itSigma$ and $\omega_{\itPi'}$ are trivial on $1+\varpi^n\o$.
\end{itemize}
Define $f_n \in {\rm c}\mbox{-}{\rm ind}_{N_3(\F)}^{P_3(\F)}(\psi_{N_3})$ by
\[
f_{n}(g) = \begin{cases}
0 & \mbox{ if $g \notin N_3(\F)K_1(n)$},\\
\psi_{N_3}(u) & \mbox{ if $g = u\cdot k\in N_3(\F)K_1(n)$}.
\end{cases}
\]
Here we regard $\GL_2(\F)$ as subgroup of $P_3(\F)$ by the map $g \mapsto \bp g & 0 \\ 0 & 1\ep$.
Note that $f_n$ is well-defined as we assume $\psi$ has conductor $\o$.
Fix a lift $W \in \mathcal{W}(\itSigma,\psi_U)$ of $f_{n}$. Let $n' \geq n$ be another sufficiently large integer such that
$W$ is right invariant by $(K(n') \times K(n'))\cap {\bf G}(\o)$.
Let $f^{(s)} = f_\Phi^{(s)}$ be the rational section of $I(\omega_\itSigma^{-1}\omega_{\itPi'}^{-1},s)$ defined by
\[
f^{(s)}(g) = |\det g|^s\int_{\F^\times} \Phi((0,t)g)\omega_\itSigma\omega_{\itPi'}(t)|t|^{2s}\,d^\times t,
\]
where $\Phi$ is the Schwartz function on $\F^2$ given by
\[
\Phi(x,y) = \omega_\itSigma^{-1}\omega_{\itPi'}^{-1}(y)\mathbb{I}_{\varpi^{n'}\o}(x)\mathbb{I}_{\o^\times}(y).
\]
It is clear that $f^{(s)}$ is holomorphic for ${\rm Re}(s)>-\tfrac{{\sf u}+{\sf w}'}{2}$.
Let $k=(k_1,k_2)\in{\bf G}(\o)$.
Note that $f^{(s)}(k_1)=0$ unless $k_1 \in K_0(n')$. Assume $k_1 \in K_0(n')$ and write
\[
k_1 \in {\rm diag}(t_1,t_2)N_2(\o)K(n')
\]
for some $t_1,t_2 \in \o^\times$ with $t_1t_2=\det k_1$.
Then $f^{(s)}(k_1) = \omega_\itSigma^{-1}\omega_{\itPi'}^{-1}(t_2)$ and we have
\begin{align*}
W({\rm diag}(a_1a_2,a_1,a_2^{-1},1)k)
 & = \omega_\itSigma(a_2^{-1}t_2) f_n({\rm diag}(a_1a_2t_2^{-1},a_2t_2^{-1},1)k_2).
\end{align*}
In particular, $W({\rm diag}(a_1a_2,a_1,a_2^{-1},1)k)=0$ unless $a_2 \in \o^\times$ and $k_2 \in K_0(n)$. Assume $k_2 \in K_0(n)$ and write
\[
k_2 \in {\rm diag}(u_1,u_2)N_2(\o)K(n)
\]
for some $u_1,u_2 \in \o^\times$ with $u_1u_2=\det k_2$.
Then $W'({\rm diag}(a_1,1)k_2)  =  \omega_{\itPi'}(u_2)\mathbb{I}_{\o^\times}(a_1)$ and we have
\begin{align*}
I_s(W',W;k) &= \omega_{\itSigma}(t_2)\omega_{\itPi'}(u_2)\int_{\F^\times}d^\times a_1\int_{\F^\times}d^\times a_2\, \mathbb{I}_{\o^\times}(a_1)\mathbb{I}_{1+\varpi^{n}\o}(a_2t_2^{-1}u_2)\omega_{\itPi'}(a_2)\\
& =  [\o^\times:1+\varpi^n\o]^{-1} \omega_{\itSigma}\omega_{\itPi'}(t_2).
\end{align*}
We conclude that 
\[
Z(f^{(s)},W',W) = [\o^\times:1+\varpi^n\o]^{-1}[{\bf G}(\o):(K_0(n') \times K_0(n))\cap {\bf G}(\o)]^{-1}.
\]
This completes the proof.
\end{proof}

\begin{lemma}\label{L:archi. local zeta}
Assume $v \in S_\infty$, $\kappa = |\lambda_{1,v}+\lambda_{2,v}-\ell_v| \geq 1$, and $\psi_v$ is the standard additive character of $\F_v$. Let $\epsilon_v = 1$ if $v \in I$ and $\epsilon_v = -1$ if $v \in J$. We have
\begin{align*}
&Z\left(\rho({\rm diag}(\epsilon_v,1))f_{v,\kappa}^{(s)},\, W_{(\ell_v;\,{\sf w}')}^+,\,W_{((\underline{\lambda}_v;\,{\sf u}),\,\ell_v-\lambda_{2,v})}^+\right)\\
& = 
2^{-s+2+(3\lambda_{1,v}-\lambda_{2,v}-\ell-{\sf u}-{\sf w}')/2} \cdot (\sqrt{-1})^{\lambda_{2,v}-\ell_v}\cdot \Gamma_\C\left(s+\tfrac{\lambda_{1,v}-\lambda_{2,v}-\ell_v}{2}+\tfrac{{\sf u}+{\sf w}'}{2}\right)\\
&\times \Gamma_\C\left(s+\tfrac{\lambda_{1,v}-\lambda_{2,v}+\ell_v-2}{2}+\tfrac{{\sf u}+{\sf w}'}{2}\right)\Gamma_\C\left(s+\tfrac{\lambda_{1,v}+\lambda_{2,v}+\ell_v-2}{2}+\tfrac{{\sf u}+{\sf w}'}{2}\right).
\end{align*}
\end{lemma}

\begin{proof}
We identify $\F_v=\R$ and drop the subscript $v$ for brevity. 
By explicit formulas (\ref{E:GL_2 Whittaker}) and Theorem \ref{T:Moriyama}, we have
\begin{align*}
&Z\left(\rho({\rm diag}(\epsilon,1))f_{\kappa}^{(s)},\, W_{(\ell;\,{\sf w}')}^+,\,W_{((\underline{\lambda};\,{\sf u}),\,\ell-\lambda_{2})}^+\right)\\
& = \int_{{\rm SO}(2)}dk_{1}\int_{{\rm SO}(2)}dk_{2}\int_{\R^\times}\frac{d^\times a_1}{|a_1|^2}\int_{\R^\times}\frac{d^\times a_2}{|a_2|^2}\,f_\kappa^{(s)}({\rm diag}((-1)^\delta a_1a_2,a_2^{-1})k_1) W_{(\ell;\,{\sf w}')}^+({\rm diag}(a_1,1)k_2)\\
&\quad\quad\quad\quad\quad\quad\quad\quad\quad\quad\quad\quad\quad\quad\quad\quad\quad\quad\quad\quad\quad\times W_{((\underline{\lambda};\,{\sf u}),\,\ell-\lambda_{2})}^+({\rm diag}(a_1a_2,a_1,a_2^{-1},1)(k_1,k_2))\\
& = 2(2\pi\sqrt{-1})^{\lambda_2-\ell}\int_0^\infty d^\times a_1\int_0^\infty d^\times a_2\, a_1^{s+({\sf u}+{\sf w}'+\ell-4)/2} a_2^{2s+{\sf u}+{\sf w}'+\lambda_2-\ell-2}e^{-4\pi a_1}\\
&\quad\quad\quad\times\int_{c_1-\sqrt{-1}\infty}^{c_1+\sqrt{-1}\infty}\frac{ds_1}{2\pi \sqrt{-1}}\,\int_{c_2-\sqrt{-1}\infty}^{c_2+\sqrt{-1}\infty}\frac{ds_2}{2\pi \sqrt{-1}}\,2^{-s_1-s_2}\Gamma_\R(s_1+\lambda_{1}+1)\Gamma_\R(-s_2-\lambda_{2})\\
&\quad\quad\quad\quad\times\Gamma_\R(s_1+s_2+\lambda_{1}-\lambda_{2}+2)\Gamma_\R(s_1+s_2+\lambda_{1}+\lambda_{2}+2)\frac{\Gamma(s_1+\lambda_{1}+1+\ell-\lambda_2)}{\Gamma(s_1+\lambda_{1}+1)} a_1^{(-s_1-s_2)/{2}}a_2^{-s_1}.
\end{align*}
Note that
\begin{align*}
\int_0^\infty a_1^{s+({\sf u}+{\sf w}'+\ell-s_1-s_2-4)/2}e^{-4\pi a_1}\,d^\times a_1 = 2^{-2s-(u+{\sf w}'+\ell-s_1-s_2-4)}\Gamma_\R(2s+{\sf u}+{\sf w}'+\ell-s_1-s_2-4).
\end{align*}
Put $t = s+\tfrac{{\sf u}+{\sf w}'}{2}$. Then we have
\begin{align*}
&Z\left(\rho({\rm diag}(\epsilon,1))f_{\kappa}^{(s)},\, W_{(\ell;\,{\sf w}')}^+,\,W_{((\underline{\lambda};\,{\sf u}),\,\ell-\lambda_{2})}^+\right)\\
& = 2^{-2t-\ell+5}(2\pi\sqrt{-1})^{\lambda_2-\ell}\\
&\times\int_{0}^\infty d^\times a_2 \, a_2^{2t+\lambda_2-\ell-2}\int_{c_1-\sqrt{-1}\infty}^{c_1+\sqrt{-1}\infty} \frac{ds_1}{2\pi\sqrt{-1}}\,a_2^{-s_1}\Gamma_\R(s_1+\lambda_1+1)\frac{\Gamma(s_1+\lambda_1+1+\ell-\lambda_2)}{\Gamma(s_1+\lambda_1+1)}\\
&\quad\quad\quad\quad\quad\quad\quad\quad\quad\quad\times\int_{c_2-\sqrt{-1}\infty}^{c_2+\sqrt{-1}\infty} \frac{ds_2}{2\pi\sqrt{-1}}\,\Gamma_\R(s_1+s_2+\lambda_{1}-\lambda_{2}+2)\Gamma_\R(s_1+s_2+\lambda_{1}+\lambda_{2}+2)\\
&\quad\quad\quad\quad\quad\quad\quad\quad\quad\quad\quad\quad\quad\quad\quad\quad\quad\quad\quad\quad\quad\quad\quad\quad\times\Gamma_\R(-s_2-\lambda_2)\Gamma_\R(2t+\ell-s_1-s_2-4).
\end{align*}
By Barne's first lemma, the above integration in $s_2$ is equal to 
\begin{align*}
2\cdot\frac{\Gamma_\R(s_1+\lambda_1-2\lambda_2+2)\Gamma_\R(2t+\lambda_1-\lambda_2+\ell-2)\Gamma_\R(s_1+\lambda_1+2)\Gamma_\R(2t+\lambda_1+\lambda_2+\ell-2)}{\Gamma_\R(2t+s_1+2\lambda_1-\lambda_2+\ell)}.
\end{align*}
Therefore, we have
\begin{align*}
&Z\left(\rho({\rm diag}(\epsilon,1))f_{\kappa}^{(s)},\, W_{(\ell;\,{\sf w}')}^+,\,W_{((\underline{\lambda};\,{\sf u}),\,\ell-\lambda_{2})}^+\right)\\
& = 2^{-2t-\ell+6}(2\pi\sqrt{-1})^{\lambda_2-\ell}\Gamma_\R(2t+\lambda_1-\lambda_2+\ell-2)\Gamma_\R(2t+\lambda_1+\lambda_2+\ell-2)\\
&\times\int_{0}^\infty d^\times a_2 \, a_2^{2t+\lambda_2-\ell-2}\int_{c_1-\sqrt{-1}\infty}^{c_1+\sqrt{-1}\infty} \frac{ds_1}{2\pi\sqrt{-1}}\,a_2^{-s_1}\Gamma_\C(s_1+\lambda_1+1)\frac{\Gamma(s_1+\lambda_1+1+\ell-\lambda_2)}{\Gamma(s_1+\lambda_1+1)}\\
&\quad\quad\quad\quad\quad\quad\quad\quad\quad\quad\quad\quad\quad\quad\quad\quad\quad\quad\quad\quad\quad\quad\quad\quad\quad\times \frac{\Gamma_\R(s_1+\lambda_1-2\lambda_2+2)}{\Gamma_\R(2t+s_1+2\lambda_1-\lambda_2+\ell)}\\
& = 2^{-4t-\ell+7}\cdot(\sqrt{-1})^{\lambda_2-\ell}\cdot\Gamma_\R(2t+\lambda_1-\lambda_2-\ell)\Gamma_\R(2t+\lambda_1-\lambda_2+\ell-2)\Gamma_\R(2t+\lambda_1+\lambda_2+\ell-2).
\end{align*}
Here the last equality follows from the Mellin inversion formula. Note that $\Gamma_\R(2s) = 2^{s-1}\Gamma_\C(s)$. This completes the proof.
\end{proof}

\subsection{Algebraicity of critical values}

We keep the notation and assumptions as in \S\,\ref{SS:coho. inter.}. Let $S$ be a finite set of places containing $S_\infty$ such that $\itPi'_v$ and $\itSigma_v$ are unramified for $v \notin S_\infty$.
Let $L^S(s,\itSigma \times \itPi')$ be the partial $L$-function with respect to $S$.
In the following theorem, we prove algebraicity of the critical values for $L^S(s,\itSigma \times \itPi')$ in terms of the automorphic periods $p^I(\itSigma)$ and $p^J(\itPi')$.

\begin{thm}\label{T:algebraicity GSp_4 x GL_2}
Let $I,J$ be subsets of $S_\infty$ defined in (\ref{E:period set}).
Assume the following conditions are satisfied:
\begin{itemize}
\item $I\cup J=S_\infty$.
\item $\kappa = |\lambda_{1,v}+\lambda_{2,v}-\ell_v|$ for all $v \in S_\infty$. 
\item If $\F=\Q$, then $\kappa \neq 2$ or $\omega_\itSigma\omega_{\itPi'} \neq |\mbox{ }|_{\A_\F}^{{\sf u}+{\sf w}'}$.
\end{itemize}
Let
\begin{align*}
p(\itSigma \times \itPi')
& = |D_\F|^{1/2}\cdot(2\pi\sqrt{-1})^{\sum_{v \in S_\infty}(3\lambda_{1,v}-\lambda_{2,v}+\ell_v+5{\sf u}+5{\sf w}')/2} \cdot (\sqrt{-1})^{\sum_{v \in I}\lambda_{1,v}+\sum_{v \in J}(\lambda_{2,v}+\ell_v)}\\
& \times G(\omega_\itSigma\omega_{\itPi'}^2) \cdot p^I(\itSigma)\cdot p^J(\itPi').
\end{align*}
\begin{itemize}
\item[(1)] Let $m = \tfrac{\kappa}{2}-\tfrac{{\sf u}+{\sf w}'}{2}$ be the rightmost critical point. We have
\[
\sigma \left( \frac{L^S(m,\itSigma \times \itPi')}{(2\pi\sqrt{-1})^{4dm}\cdot p(\itSigma \times \itPi')} \right) = \frac{L^S(m,{}^\sigma\!\itSigma \times {}^\sigma\!\itPi')}{(2\pi\sqrt{-1})^{4dm}\cdot p({}^\sigma\!\itSigma \times {}^\sigma\!\itPi')}
\]
for all $\sigma \in {\rm Aut}(\C)$. In particular, we have
\[
 \frac{L^S(m,\itSigma \times \itPi')}{(2\pi\sqrt{-1})^{4dm}\cdot p(\itSigma \times \itPi')} \in \Q(\itSigma)\Q(\itPi')\Q(I).
\]
\item[(2)] Assume $\itSigma_v$ and $\itPi_v'$ are discrete series representations for all $v \in S_\infty$. Then the assertion in (1) holds for all critical points.
\end{itemize}
\end{thm}

\begin{proof}
We begin with the first assertion. 
Let $\psi=\psi_\F$ be the standard additive character of $\F$ 
and write $T=S\smallsetminus S_\infty$.
By Lemma \ref{L:non-archi. local zeta}, there exist $W_T' \in \mathcal{W}(\itPi_T',\psi_T)$, $W_T \in \mathcal{W}(\itSigma_T,\psi_{U,T})$, and a rational section $f_T^{(s)}$ of $I(\omega_{\itSigma,T}^{-1}\omega_{\itPi',T}^{-1},s)$ holomorphic for ${\rm Re}(s)>-\tfrac{{\sf u}+{\sf w}'}{2}$ such that
\[
Z(f_T^{(s)},W_T',W_T) = 1.
\]
Moreover, combine with (\ref{E:Galois Gauss sum}), for all $\sigma \in {\rm Aut}(\C)$ we have
\begin{align}\label{E:algebraicity proof 1}
Z({}^\sigma\!f_T^{(s)},t_{\sigma,T}W_T',t_{\sigma,T}W_T) = \frac{\sigma (G(\omega_{\itPi'}))}{G({}^\sigma\!\omega_{\itPi'})}.
\end{align}
Let $\varphi' \in V_{\itPi'}^+$, $\varphi \in V_\itSigma^+$, and $f^{(s)}$ the rational section of $I(\omega_{\itSigma,f}^{-1}\omega_{\itPi',f}^{-1},s)$ defined so that
\[
W_{\varphi'}^{(\infty)} = \prod_{v \notin S}W_{v,\psi_v}^\circ \cdot W_T',\quad W_{\varphi}^{(\infty)} = \prod_{v \notin S}W_{v,\psi_{U,v}}^\circ \cdot W_T,\quad f^{(s)} = \bigotimes_{v \notin S} f_{v,\circ}^{(s)}\otimes f_T^{(s)}.
\]
For $v \in S_\infty$, let $\epsilon_v = 1$ if $v \in I$ and $\epsilon_v = -1$ otherwise.
By Proposition \ref{P:integral rep.}, we have
\begin{align*}
Z\left( E^{[\kappa]}(f^{(s)}),\, (\varphi')^J,\, {\rm pr}_{\underline{\ell}-\underline{\lambda}_2}(\varphi)^I\right) &= |D_\F|^{-2}\zeta_\F(2)^{-2}\cdot L^S(m,\itSigma \times \itPi')\\
&\times\prod_{v \in S_\infty}Z\left.\left(\rho({\rm diag}(\epsilon_v,1))f_{v,\kappa}^{(s)},\, W_{(\ell_v;\,{\sf w}')}^+,\,W_{((\underline{\lambda}_v;\,{\sf u}),\,\ell_v-\lambda_{2,v})}^+\right)\right\vert_{s=m}.
\end{align*}
By Lemma \ref{L:archi. local zeta}, up to a non-zero rational number, the above product of archimedean local zeta integrals is equal to 
\[
\pi^{-3dm-\sum_{v \in S_\infty}(3\lambda_{1,v}-\lambda_{2,v}+\ell_v+3{\sf u}+3{\sf w}'-4)/2}\cdot(\sqrt{-1})^{\sum_{v \in S_\infty}(\lambda_{2,v}-\ell_v)}.
\]
Also note that $\zeta_\F(2) \in \pi^d\cdot\Q^\times$. Therefore we can rewrite the above equality as follows:
\begin{align*}
Z\left( E^{[\kappa]}(f^{(s)}),\, (\varphi')^J,\, {\rm pr}_{\underline{\ell}-\underline{\lambda}_2}(\varphi)^I\right)  &= C\cdot \pi^{-3dm-\sum_{v \in S_\infty}(3\lambda_{1,v}-\lambda_{2,v}+\ell_v+3{\sf u}+3{\sf w}')/2}\cdot(\sqrt{-1})^{\sum_{v \in S_\infty}(\lambda_{2,v}-\ell_v)}\\
&\times L^S(m,\itSigma\times\itPi')
\end{align*}
for some $C \in \Q^\times$ depending only on $\F$, $\itSigma_\infty$ and $\itPi'_\infty$.
Similarly, together with (\ref{E:algebraicity proof 1}), we have
\begin{align*}
Z\left(E^{[\kappa]}({}^\sigma\!f^{(s)}),\, ({}^\sigma\!\varphi')^{{}^\sigma\!J},\, {\rm pr}_{{}^\sigma\!\underline{\ell}-{}^\sigma\!\underline{\lambda}_2}({}^\sigma\!\varphi)^{{}^\sigma\!I}\right)  &= C\cdot \pi^{-3dm-\sum_{v \in S_\infty}(3\lambda_{1,v}-\lambda_{2,v}+\ell_v+3{\sf u}+3{\sf w}')/2}\cdot(\sqrt{-1})^{\sum_{v \in S_\infty}(\lambda_{2,v}-\ell_v)}\\
&\times \frac{\sigma (G(\omega_{\itPi'}))}{G({}^\sigma\!\omega_{\itPi'})}\cdot L^S(m,{}^\sigma\!\itSigma\times{}^\sigma\!\itPi')
\end{align*}
for all $\sigma \in {\rm Aut}(\C)$.
The first assertion then follows from Proposition \ref{P:Galois equiv. global}.

Now we consider the second assertion. We may assume $L(s,\itSigma \times \itPi')$ has more than one critical point, that is, $\kappa \geq 2$.
For any finite place $v$ of $\F$, since $\itSigma_v$ and $\itPi_v'$ are essentially unitary, the local factor $L(s,\itSigma_v \times \itPi'_v)$ is holomorphic for ${\rm Re}(s) \geq 1-\tfrac{{\sf u}+{\sf w}'}{2}$. In particular $L(\tfrac{\kappa}{2}-\tfrac{{\sf u}+{\sf w}'}{2},\itSigma_v\times\itPi'_v) \neq0$. Also note that (cf.\,\cite[Proposition 5.4]{Morimoto2014})
\[
{}^\sigma\! L(s,\itSigma_v \times \itPi_v') = L(s,{}^\sigma\!\itSigma_v \times {}^\sigma\!\itPi_v')
\]
as rational functions in $q_v^{-s}$ for all $\sigma \in {\rm Aut}(\C)$. 
Therefore, assertion (1) actually holds for $L^{(\infty)}(\tfrac{\kappa}{2}-\tfrac{{\sf u}+{\sf w}'}{2},\itSigma\times\itPi')$. 
Since we have assume $\itSigma_v$ is a discrete series representation for all $v \in S_\infty$, the functorial lift of $\itSigma$ to $\GL_4(\A_\F)$ is a regular $C$-algebraic isobaric automorphic representation of $\GL_4(\A_\F)$. Together with the assumption on $\itPi'_v$ for $v \in S_\infty$, by the result of Harder and Raghuram \cite[Theorem 7.21]{HR2020}, we have
\begin{align}\label{E:algebraicity proof 2}
\sigma \left( \frac{L^{(\infty)}(m-1,\itSigma \times \itPi')}{(2\pi\sqrt{-1})^{-4d}\cdot L^{(\infty)}(m,\itSigma \times \itPi')}\right) = \frac{L^{(\infty)}(m-1,{}^\sigma\!\itSigma \times {}^\sigma\!\itPi')}{(2\pi\sqrt{-1})^{-4d}\cdot L^{(\infty)}(m,{}^\sigma\!\itSigma \times {}^\sigma\!\itPi')}
\end{align}
for all $\sigma \in {\rm Aut}(\C)$ and $m \in {\rm Crit}(\itSigma \times \itPi')$ such that $L^{(\infty)}(m,\itSigma \times \itPi') \neq 0$.
By \cite[Theorem 5.2]{Shahidi1981}, the non-vanishing condition is satisfied for all critical points $m \geq 1-\tfrac{{\sf u}+{\sf w}'}{2}$.
We thus conclude from (\ref{E:algebraicity proof 2}) that assertion (1) holds for all right-half critical points $m$, that is, $m \geq \tfrac{1}{2}-\tfrac{{\sf u}+{\sf w}'}{2}$.
For the left-half critical points, the assertion follows from the global functional equation.
This completes the proof.
\end{proof}

\subsection{A period relation}\label{SS:Morimoto}

We keep the notation and assumptions as in \S\,\ref{SS:coho. inter.}.
Assume further the following conditions are satisfied:
\begin{itemize}
\item The functorial lift of $\itSigma$ to $\GL_4(\A_\F)$ is cuspidal.
\item $\itSigma_v$ is a discrete series representation for all $v \in S_\infty$.
\end{itemize}
For $v \in S_\infty$, let $\itSigma_v^{\rm hol}$ be the holomorphic discrete series representation of $\GSp_4(\F_v)$ which belongs to the $L$-packet containing $\itSigma_v$. 
Let $\itSigma^{\rm hol}$ be the irreducible admissible representation of $\GSp_4(\A_\F)$ defined by
\begin{align}\label{E:holomorphic}
\itSigma^{\rm hol} = \bigotimes_{v\in S_\infty} \itSigma_v^{\rm hol} \otimes \itSigma_f.
\end{align}
Since the transfer of $\itSigma$ to $\GL_4(\A_\F)$ is cuspidal, it follows from Arthur's multiplicity formula proved by Gee and Ta\"{i}bi \cite[Theorem 7.4.1]{GT2019} that $\itSigma^{\rm hol}$ appears in the automorphic discrete spectrum of $\GSp_4(\A_\F)$. It then follows from the temperedness of $\itSigma_v^{{\rm hol}}$ for all $v \in S_\infty$ and the result of Wallach \cite[Theorem 4.3]{Wallach1984} that $\itSigma^{\rm hol}$ is cuspidal. 
Recall the period $\Omega(\itSigma^{\rm hol}) \in \C^\times$ in Theorem \ref{T:Liu}.
In this section, we establish a period relation between the periods $p^\varnothing(\itSigma)$ and $\Omega(\itSigma^{\rm hol})$ in Proposition \ref{P:period relation 1} below.
The period relation is a consequence of Theorem \ref{T:algebraicity GSp_4 x GL_2} for $I=\varnothing$ and a result of Morimoto \cite{Morimoto2018} recalled in the following theorem.

\begin{thm}[Morimoto]\label{T:Morimoto}
Let $\itPi'$ be an irreducible cuspidal automorphic representation of $\GL_2(\A_\F)$.
Assume $\itPi'$ is regular $C$-algebraic of weight $(\underline{\ell};\,{\sf w}')$ and 
\[
\lambda_{1,v}+\lambda_{2,v}+5 \leq \ell_v \leq \lambda_{1,v}-\lambda_{2,v}+5
\]
for all $v \in S_\infty$.
Let $m$ be a critical point of $L(s,\itSigma \times \itPi')$ such that $m +\tfrac{{\sf u}+{\sf w}'}{2} > 2 $. 
For $\sigma \in {\rm Aut}(\C)$, we have
\begin{align*}
&\sigma \left( \frac{L^{(\infty)}(m,\itSigma \times \itPi')}{(2\pi\sqrt{-1})^{4dm+\sum_{v \in S_\infty}(2{\sf u}+(\ell_v+5{\sf w}')/2)}\cdot (\sqrt{-1})^{d{\sf u}+d{\sf w}'}\cdot G(\omega_\itSigma\omega_{\itPi'})^2\cdot \Omega(\itSigma^{\rm hol})\cdot p^{S_\infty}(\itPi')}\right)\\
& = \frac{L^{(\infty)}(m,{}^\sigma\!\itSigma \times {}^\sigma\!\itPi')}{(2\pi\sqrt{-1})^{4dm+\sum_{v \in S_\infty}(2{\sf u}+(\ell_v+5{\sf w}')/2)}\cdot (\sqrt{-1})^{d{\sf u}+d{\sf w}'}\cdot G({}^\sigma\!\omega_\itSigma{}^\sigma\!\omega_{\itPi'})^2\cdot \Omega({}^\sigma\!\itSigma^{\rm hol})\cdot p^{S_\infty}({}^\sigma\!\itPi')}.
\end{align*}
\end{thm}

\begin{rmk}
In the notation of \cite[Theorem 1]{Morimoto2018}, we have
\[
\Omega(\itSigma^{\rm hol}) = \<\Phi,\Phi\>,\quad p^{S_\infty}(\itPi') = (2\pi\sqrt{-1})^{-\sum_{v \in S_\infty}(\ell_v+{\sf w}')/2}\cdot \<\Psi,\Psi\>.
\]
\end{rmk}

The result of Morimoto was proved when ${\sf u}$ is even. Proceeding similarly as in \cite[\S\,6]{Morimoto2018}, we see that in order to extend the result to arbitrary ${\sf u}$, it suffices to generalize \cite[Proposition 6.2]{Morimoto2018}.
More precisely, we have the following period relation for Yoshida lifts on $\GSp_4(\A_\Q)$.

\begin{prop}\label{P:Yoshida lifts}
Let $\itPsi$ be an $C$-algebraic irreducible cuspidal automorphic representation of $\GSp_4(\A_\Q)$ such that $\itPsi_\infty$ is a holomorphic discrete series representation. 
Assume there exist non-isomorphic irreducible cuspidal automorphic representations $\itPi_1$ and $\itPi_2$ of $\GL_2(\A_\Q)$ such that $\omega_{\itPi_1} = \omega_{\itPi_2} = \omega_\itPsi$ and the functorial lift of $\itPsi$ to $\GL_4(\A_\Q)$ is the isobaric sum $\itPi_1 \boxplus \itPi_2$.
Let $(\kappa_1;\,{\sf u})$ and $(\kappa_2;\,{\sf u})$ be the weights of $\itPi_1$ and $\itPi_2$, respectively, with $\kappa_1\geq\kappa_2$. 
For $\sigma \in {\rm Aut}(\C)$, we have
\[
\sigma \left( \frac{\Omega(\itPsi)}{(2\pi\sqrt{-1})^{\kappa_1}\cdot \Vert f_{\itPi_1}\Vert} \right) = \frac{\Omega({}^\sigma\!\itPsi)}{(2\pi\sqrt{-1})^{\kappa_1}\cdot \Vert f_{{}^\sigma\!\itPi_1}\Vert}.
\]
Here $\Omega(\itPsi) \in \C^\times$ is the period in Theorem \ref{T:Liu}.
\end{prop}

\begin{proof}

We have the factorization of $L$-functions:
\[
L(s,\itPsi,{\rm std}\otimes \chi) =  L(s,\itPi_1 \times \itPi_2^\vee\otimes \chi)\cdot L(s,\chi)
\]
for any Hecke character $\chi$ of $\A_\Q^\times$. Here $L(s,\itPi_1 \times \itPi_2^\vee\otimes \chi)$ is the Rankin--Selberg $L$-function of $\itPi_1 \times \itPi_2^\vee\otimes \chi$.
Assume $\chi$ has finite order with signature $-1$. Then $m=1$ is a critical point for $L(s,\itPsi,{\rm std}\otimes \chi)$.
For the Dirichlet $L$-function $L(s,\chi)$, we have the well-known algebraicity result that
\[
\sigma \left ( \frac{L^{(\infty)}(1,\chi)}{ (2\pi\sqrt{-1})\cdot G(\chi)}\right ) = \frac{L^{(\infty)}(1,{}^\sigma\!\chi)}{(2\pi\sqrt{-1})\cdot G({}^\sigma\!\chi)}
\]
for all $\sigma \in {\rm Aut}(\C)$.
Note that the condition $\itPsi_\infty$ is a discrete series representation implies that $\kappa_1 > \kappa_2$.
We have the algebraicity result proved by Shimura \cite[Theorem 3]{Shimura1976} that
\[
\sigma \left(\frac{L^{(\infty)}(1,\itPi_1 \times \itPi_2^\vee \otimes \chi)}{(2\pi\sqrt{-1})^{\kappa_1+2}\cdot(\sqrt{-1})^{{\sf u}}\cdot G(\chi)^2 \cdot \Vert f_{\itPi_1}\Vert} \right) = \frac{L^{(\infty)}(1,{}^\sigma\!\itPi_1 \times {}^\sigma\!\itPi_2^\vee \otimes {}^\sigma\!\chi)}{(2\pi\sqrt{-1})^{\kappa_1+2}\cdot(\sqrt{-1})^{{\sf u}}\cdot G({}^\sigma\!\chi)^2 \cdot \Vert f_{{}^\sigma\!\itPi_1}\Vert}
\]
for all $\sigma \in {\rm Aut}(\C)$.
Since $\itPi_1$ and $\itPi_2\otimes \chi^{-1}$ are non-isomorphic and $\chi$ is non-trivial, it follows that $L^{(\infty)}(1,\chi)$ and $L^{(\infty)}(1,\itPi_1 \times \itPi_2^\vee \otimes \chi)$ are non-zero. 
The period relation thus follows from Theorem \ref{T:Liu} by taking $\chi$ so that $\chi^2 \neq 1$.
This completes the proof.
\end{proof}

\begin{rmk}
When $\omega_{\itPi_1} = \omega_{\itPi_2} =1$, and $\itPi_1$ and $\itPi_2$ have square-free levels, the period relation follows from the explicit inner product formula for Yoshida lifts established by B\"ocherer--Dummigan--Schulze-Pillot \cite[Corollary 8.8]{BDS2012}.
Based on \cite{Morimoto2014}, Saha proved the period relation in \cite[Theorem 5.1]{Saha2015} assuming $\kappa_{1}\geq 12$ and $\kappa_2=2$.
\end{rmk}

\begin{prop}\label{P:period relation 1}
Assume the following conditions are satisfied:
\begin{itemize}
\item The functorial lift of $\itSigma$ to $\GL_4(\A_\F)$ is cuspidal.
\item $\itSigma_v$ is a discrete series representation for all $v \in S_\infty$.
\item $\lambda_{2,v} \leq -5$ for all $v \in S_\infty$.
\end{itemize}
Let $\itSigma^{\rm hol}$ be the $C$-algebraic irreducible cuspidal automorphic representation of $\GSp_4(\A_\F)$ defined in (\ref{E:holomorphic}).
For $\sigma \in {\rm Aut}(\C)$, we have
\begin{align*}
&\sigma \left(\frac{p^{\varnothing}(\itSigma)}{|D_\F|^{1/2}\cdot(2\pi\sqrt{-1})^{\sum_{v \in S_\infty}(-3\lambda_{1,v}+\lambda_{2,v}-{\sf u})/2}\cdot(\sqrt{-1})^{\sum_{v \in S_\infty}\lambda_{1,v}}\cdot G(\omega_{\itSigma}) \cdot \Omega(\itSigma^{\rm hol})}\right)\\
& = \frac{p^{\varnothing}({}^\sigma\!\itSigma)}{|D_\F|^{1/2}\cdot(2\pi\sqrt{-1})^{\sum_{v \in S_\infty}(-3\lambda_{1,v}+\lambda_{2,v}-{\sf u})/2}\cdot(\sqrt{-1})^{\sum_{v \in S_\infty}\lambda_{1,v}}\cdot G({}^\sigma\!\omega_{\itSigma}) \cdot \Omega({}^\sigma\!\itSigma^{\rm hol})}.
\end{align*}
\end{prop}

\begin{proof}
Let $\itPi'$ be an $C$-algebraic irreducible cuspidal automorphic representation of $\GL_2(\A_\F)$ with weight $(\underline{\ell};\,{\sf w}')$ such that
\begin{align}\label{E:4.4}
\lambda_{1,v}+\lambda_{2,v}+5 \leq \ell_v \leq \lambda_{1,v}
\end{align}
for all $v \in S_\infty$. We refer to \cite[Theorem 1.1]{Weinstein2009} for the existence of $\itPi'$. Note that the rightmost critical value of $L(s,\itSigma \times \itPi')$ is non-zero by condition (\ref{E:4.4}).
Therefore, the period relation follows immediately from Theorem \ref{T:algebraicity GSp_4 x GL_2}-(1) and Theorem \ref{T:Morimoto}.
\end{proof}

\section{Period relations}\label{S:period relations}

We keep the notation of \S\,\ref{SS:KRS lift}.
The purpose of this section is to prove the period relations in Theorem \ref{P:period relation main}. The assertions are proved based on Theorem \ref{T:algebraicity GSp_4 x GL_2}, Proposition \ref{P:period relation 1}, and algebraicity results in the literature.
Let $\itSigma$ be the Kim--Ramakrishnan--Shahidi lift of $\itPi$. Recall $\itSigma$ is regular $C$-algebraic, satisfying condition (\ref{E:discrete condition 2}) with weight $(\underline{\lambda};\,{\sf u})$, and has central character $\omega_\itPi^3$, where
\[
\underline{\lambda}_v = (2\kappa_v-1,1-\kappa_v),\quad {\sf u} = 3{\sf w}
\]
for $v \in S_\infty$.

\subsection{Proof of Theorem \ref{P:period relation main}-(1)}

In this section, we prove assertion (1) of Theorem \ref{P:period relation main}. We assume $\kappa_v \geq 6$ for all $v \in S_\infty$.
By Deligne's conjecture for symmetric fourth $L$-function proved in \cite[Theorem 1.4]{Chen2021}, for all critical points $m \geq 1$ of $L(s,\itPi,{\rm Sym}^4\otimes\omega_\itPi^{-2})$ (see also \cite{Morimoto2021}), we have
\begin{align*}
&\sigma \left(\frac{L^{(\infty)}(m,\itPi,{\rm Sym}^4 \otimes \omega_\itPi^{-2})}{|D_\F|^{1/2}\cdot(2\pi\sqrt{-1})^{3dm+3\sum_{v \mid \infty}\kappa_v}\cdot (\sqrt{-1})^{d{\sf w}}\cdot \Vert f_\itPi \Vert^3} \right) = \frac{L^{(\infty)}(m,{}^\sigma\!\itPi,{\rm Sym}^4 \otimes {}^\sigma\!\omega_\itPi^{-2})}{|D_\F|^{1/2}\cdot(2\pi\sqrt{-1})^{3dm+3\sum_{v \mid \infty}\kappa_v}\cdot (\sqrt{-1})^{d{\sf w}}\cdot \Vert f_{{}^\sigma\!\itPi} \Vert^3}
\end{align*}
for all $\sigma \in {\rm Aut}(\C)$.
Note that the critical values $L^{(\infty)}(m,\itPi,{\rm Sym}^4\otimes\omega_\itPi^{-2})$ are non-zero for all critical points $m$ as we explained in \S\,\ref{S:proof} for symmetric sixth $L$-functions.
Therefore, it follows from (\ref{E:twisted standard}) and Theorem \ref{T:Liu} for $L(s,\itSigma,{\rm std}) = L(s,\itSigma^{\rm hol},{\rm std})$ that
\[
\sigma \left( \frac{\Omega(\itSigma^{\rm hol})}{(2\pi\sqrt{-1})^{3\sum_{v \in S_\infty}\kappa_v}\cdot \Vert f_\itPi\Vert^3} \right) = \frac{\Omega({}^\sigma\!\itSigma^{\rm hol})}{(2\pi\sqrt{-1})^{3\sum_{v \in S_\infty}\kappa_v}\cdot \Vert f_{{}^\sigma\!\itPi}\Vert^3} 
\]
for all $\sigma \in {\rm Aut}(\C)$.
Note that the assumption $\kappa_{v} \geq 6$ is equivalent to $\lambda_{2,v} \leq -5$ for $v \in S_\infty$.
We thus conclude (1) of Theorem \ref{P:period relation main} from the above period relation and Proposition \ref{P:period relation 1}.

\subsection{Proof of Theorem \ref{P:period relation main}-(2) and (3)}
In this section, we prove assertions (2) and (3) of Theorem \ref{P:period relation main}. 
We assume $\kappa_v \geq 3$ for all $v \in S_\infty$.
Firstly we explain the idea of the proof.
We consider the Rankin--Selberg $L$-function $L(s,\itSigma \times \itPi')$ for some auxiliary $C$-algebraic irreducible cuspidal automorphic representation $\itPi'$ of $\GL_2(\A_\F)$ with weight $(\underline{\ell};\,{\sf w}')$ which satisfies the following conditions:
\begin{itemize}
\item[(i)] $\itPi'$ is the automorphic induction $I_\K^\F(\chi')$ of some Hecke character $\chi'$ of $\A_\K^\times$ for some CM-extension $\K/\F$.
\item[(ii)] The set $I$ in (\ref{E:period set}) is equal to $S_\infty$, that is, $1 \leq \ell_v \leq \kappa_v-1$ for all $v \in S_\infty$.
\item[(iii)] We have $\kappa_v-\ell_v = \kappa_w-\ell_w$ for all $v,w \in S_\infty$.
\item[(iv)] If $\F=\Q$, then $\omega_\itSigma\omega_{\itPi'} \neq |\mbox{ }|_{\A_\F}^{{\sf u}+{\sf w}'}$.
\end{itemize}
By conditions (ii)-(iv), the algebraicity of the rightmost critical value of $L(s,\itSigma \times \itPi')$ can be expressed in terms of $p^{S_\infty}(\itSigma)$ and some fudge factors by Lemma \ref{L:period relation 1}-(1) and Theorem \ref{T:algebraicity GSp_4 x GL_2}-(1).
On the other hand, consider the functorial lift ${\rm Sym}^3(\itPi)$ of $\itPi$ to $\GL_4(\A_\F)$ with respect to the symmetric cube representation of $\GL_2(\C)$. Note that ${\rm Sym}^3(\itPi)$ is cuspidal since $\itPi$ is non-CM.
Let ${\rm BC}_\K({\rm Sym}^3(\itPi))$ be the base change lift of ${\rm Sym}^3(\itPi)$ to $\GL_4(\A_\K)$. By condition (i) and the adjointness property between automorphic induction and base change, we have
\begin{align}\label{E:period relation proof 3}
L(s,\itSigma \times \itPi') = L(s,{\rm BC}_\K({\rm Sym}^3(\itPi)) \otimes \chi').
\end{align}
We then use results in the literature \cite{GL2016}, \cite{Harris2021}, and \cite{JST2021}, together with Deligne's conjecture for symmetric cube $L$-functions of $\itPi$ \cite{GH1993} and \cite{Chen2021d} to study the algebraicity of the $L$-function on the right-hand side of (\ref{E:period relation proof 3}).
As a consequence, we obtain period relation between $p^{S_\infty}(\itSigma)$ and $\Vert f_\itPi \Vert^4$.

Let $\K$ be a totally imaginary quadratic extension over $\F$, $c \in \Gal(\K/\F)$ the non-trivial automorphism, and $\omega_{\K/\F}$ the quadratic Hecke character of $\A_\F^\times$ associated to $\K/\F$ by class field theory. 
For each $v \in S_\infty$, fix a complex embedding $j_v$ of $\K$ lying over $v$ and we identify $\K_v$ with $\C$ via $j_v$.
Then $\Phi = \left\{ j_v,\,\vert\,v \in S_\infty\right\}$ is a CM-type of $\K$.
Let $\chi$ be an algebraic Hecke character of $\A_\K^\times$. There exist $\underline{\ell} = (\ell_v)_{v \in S_\infty} \in \Z^{S_\infty}$ and ${\sf w}' \in \Z$ such that
\begin{align}\label{E:infinity type}
\chi_v(z) = z^{(\ell_v+{\sf w}'-2)/2}\overline{z}^{(-\ell_v+{\sf w}')/2},\quad \ell_v \equiv {\sf w}'\,({\rm mod}\,2)
\end{align}
for all $v \in S_\infty$.
Let $I_\K^\F(\chi|\mbox{ }|_{\A_\K}^{1/2})$ be the automorphic induction of $\chi|\mbox{ }|_{\A_\K}^{1/2}$ to $\GL_2(\A_\F)$. It is an isobaric automorphic representation of $\GL_2(\A_\F)$ with central character $\chi\vert_{\A_\F^\times}|\mbox{ }|_{\A_\F}\cdot\omega_{\K/\F}$. If we assume further that $\chi \neq \chi^c$ and $\underline{\ell} \in \Z_{\geq 1}^{S_\infty}$, then $I_\K^\F(\chi|\mbox{ }|_{\A_\K}^{1/2})$ is cuspidal, $C$-algebraic, and satisfying condition (\ref{E:discrete condition 1}) of weight $(\underline{\ell};\,{\sf w}')$.
We have the following result on the algebraicity of critical values of Rankin--Selberg $L$-functions on $\GL_4(\A_\K) \times \GL_1(\A_\K)$, which is a special case of the results of Guerberoff and Lin \cite[Theorem 2]{GL2016} and Harris \cite[Theorem 7.1]{Harris2021} for $n=4$.

\begin{thm}[Guerberoff--Lin, Harris]\label{T:GLH}
Let $\itPsi$ be a regular $C$-algebraic conjugate self-dual irreducible cuspidal automorphic representation of $\GL_4(\A_\K)$. For each $v \in S_\infty$, let $\{z^{a_{i,v}}\overline{z}^{-a_{i,v}}\}_{1 \leq i \leq 4}$ be the infinity type of $\itPsi$ at $v$ arranged so that $a_{1,v} > a_{2,v} > a_{3,v} > a_{4,v}$.
There exists a sequence of non-zero complex numbers $(P({}^\sigma\!\itPsi))_{\sigma \in {\rm Aut}(\C)}$ satisfying the following property:
Let $\chi$ be an algebraic Hecke character of $\A_\K^\times$ such that
\begin{align}\label{E:period relation proof 1}
2a_{2,v} > 1-\ell_v > 2a_{3,v}
\end{align}
for all $v \in S_\infty$, where $\underline{\ell}$ is determined as in (\ref{E:infinity type}). For any critical point $m+\tfrac{1}{2} \in \Z+\tfrac{1}{2}$ of $L(s,\itPsi\otimes\chi)$, we have
\begin{align}\label{E:GLH}
\sigma \left( \frac{L^S(m+\tfrac{1}{2},\itPsi \otimes \chi)}{(2\pi\sqrt{-1})^{4dm+2d{\sf w}'}\cdot G(\chi\vert_{\A_\F^\times})^2 \cdot P(\itPsi)} \right) = \frac{L^S(m+\tfrac{1}{2},{}^\sigma\!\itPsi \otimes {}^\sigma\!\chi)}{(2\pi\sqrt{-1})^{4dm+2d{\sf w}'}\cdot G({}^\sigma\!\chi\vert_{\A_\F^\times})^2 \cdot P({}^\sigma\!\itPsi)}
\end{align}
for all $\sigma \in {\rm Aut}(\C/\K^{\Gal})$. Here $\K^{\Gal}$ is the Galois closure of $\K$ in $\C$ and $S$ is a sufficiently large set of places containing $S_\infty$.
\end{thm}

\begin{proof}
In order to apply \cite[Theorem 2]{GL2016} in our case under assumption (\ref{E:period relation proof 1}), we need to verify that $\itPsi$ can be descend to a cohomological irreducible cuspidal automorphic representation of the quasi-split unitary group ${\rm U}(2,2)(\A_\F)$ and so that the descent is strong at $v \in S_\infty$ (cf.\,the assumption in the beginning of \cite[\S\,4.3]{GL2016}). The existence of descent is guaranteed by Arthur's multiplicity formula \cite{Mok2015}. By \cite[Theorem 2]{GL2016} and \cite[Theorem 7.1]{Harris2021}, we have
\begin{align}\label{E:period relation proof 2}
\sigma \left( \frac{L^S(m+\tfrac{1}{2},\itPsi \otimes \chi)}{(2\pi\sqrt{-1})^{4dm+2d}\cdot P(\itPsi) \cdot p(\widecheck{\chi},\Phi)^2\cdot p(\widecheck{\chi},\Phi^c)^2} \right) = \frac{L^S(m+\tfrac{1}{2},{}^\sigma\!\itPsi \otimes {}^\sigma\!\chi)}{(2\pi\sqrt{-1})^{4dm+2d}\cdot P({}^\sigma\!\itPsi) \cdot p({}^\sigma\!\widecheck{\chi},{}^\sigma\!\Phi)^2\cdot p({}^\sigma\!\widecheck{\chi},{}^\sigma\!\Phi^c)^2}
\end{align}
for all $\sigma \in {\rm Aut}(\C/\K^{\Gal})$ and critical points $m+\tfrac{1}{2}$.
Here $\widecheck{\chi} = (\chi^c)^{-1}$, and $p(\widecheck{\chi},\Phi)$ and $p(\widecheck{\chi},\Phi^c)$ are the CM-periods of $\widecheck{\chi}$ with respect to the CM-types $\Phi$ and $\Phi^c$, respectively, defined in \cite[\S\,1]{Harris1993}.
By the properties of CM-periods \cite[Proposition 1.4 and Lemma 1.6]{Harris1993} and \cite[(1.10.9) and (1.10.10)]{Harris1997}, we have
\[
\sigma \left( \frac{p(\widecheck{\chi},\Phi)\cdot p(\widecheck{\chi},\Phi^c)}{(2\pi\sqrt{-1})^{d({\sf w}'-1)}\cdot G(\chi\vert_{\A_\F^\times})}  \right) = \frac{p({}^\sigma\!\widecheck{\chi},{}^\sigma\!\Phi)\cdot p({}^\sigma\!\widecheck{\chi},{}^\sigma\!\Phi^c)}{(2\pi\sqrt{-1})^{d({\sf w}'-1)}\cdot G({}^\sigma\!\chi\vert_{\A_\F^\times})} 
\]
for all $\sigma \in {\rm Aut}(\C)$.
Thus we can replace $p(\widecheck{\chi},\Phi)^2\cdot p(\widecheck{\chi},\Phi^c)^2$ by $(2\pi\sqrt{-1})^{2d({\sf w}'-1)}\cdot G(\chi\vert_{\A_\F^\times})^2$ in (\ref{E:period relation proof 2}).
\end{proof}

\begin{rmk}\label{R:5.2}
As explained in \cite[Remarks 3.7.1 and 3.7.2]{GL2016}, in order to prove the Galois equivariance property (\ref{E:GLH}) over ${\rm Aut}(\C)$, it suffices to refine the Galois equivariance property in \cite[Lemma 3.3.1]{GL2016} from ${\rm Aut}(\C/\K^{\Gal})$ to ${\rm Aut}(\C)$.
\end{rmk}

We prove assertion (2) first. Assume $\kappa_v \geq 3$ for all $v \in S_\infty$. 
Assume $\K$ is chosen so that ${\rm BC}_\K({\rm Sym}^3(\itPi))$ is cuspidal. 
Fix an algebraic Hecke character $\eta$ of $\A_\K^\times$ such that the following conditions are satisfied:
\begin{itemize}
\item $\eta_v(z) = z^{-3{\sf w}}$ for all $v \in S_\infty$.
\item $\eta\vert_{\A_\F^\times} = \omega_\itPi^{-3}$.
\end{itemize}
We refer to \cite[Lemma 4.1.4]{CHT2008} for the existence of $\eta$. Let
\[
\itPsi = {\rm BC}_\K({\rm Sym}^3(\itPi)) \otimes \eta.
\]
Then $\itPsi$ is $C$-algebraic and conjugate-self dual. Indeed, we have ${\rm Sym}^3(\itPi)^\vee = {\rm Sym}^3(\itPi) \otimes \omega_\itPi^{-3}$.
Note that ${\rm BC}_\K({\rm Sym}^3(\itPi))^c = {\rm BC}_\K({\rm Sym}^3(\itPi))$. Thus we have 
\[
{\rm BC}_\K({\rm Sym}^3(\itPi))^\vee = {\rm BC}_\K({\rm Sym}^3(\itPi))^c \otimes \omega_\itPi^{-3}\circ{\rm N}_{\K/\F}.
\]
Also the infinity type of $\itPsi$ at $v \in S_\infty$ is given by $\{z^{a_{i,v}}\overline{z}^{-a_{i,v}}\}_{1 \leq i \leq 4}$ with 
\[
a_{1,v} = \tfrac{3\kappa_v-3-3{\sf w}}{2},\quad a_{2,v} = \tfrac{\kappa_v-1-3{\sf w}}{2},\quad a_{3,v} = \tfrac{1-\kappa_v-3{\sf w}}{2},\quad a_{4,v} = \tfrac{3-3\kappa_v-3{\sf w}}{2}.
\]
Fix $\underline{\ell} \in \Z^{S_\infty}$ and ${\sf w}' \in \Z$ such that the following conditions are satisfied:
\begin{itemize}
\item $\ell_v \equiv {\sf w}'\,({\rm mod}\,2)$ for all $v \in S_\infty$.
\item $1 \leq \ell_v \leq \kappa_v -2$ for all $v \in S_\infty$. 
\item $\kappa_v -\ell_v = \kappa_w-\ell_w$ for all $v,w \in S_\infty$.
\end{itemize}
Let $\chi$ be an algebraic Hecke character of $\A_\K^\times$ such that $\chi_v$ is given by (\ref{E:infinity type}) for $v \in S_\infty$ and $\chi \neq \chi^c$. Let $\itPi' = I_\K^\F(\chi|\mbox{ }|_{\A_\K}^{1/2})$.
We assume further that $\omega_\itPi^3\omega_{\itPi'} \neq |\mbox{ }|_{\A_\F}^{3{\sf w}+{\sf w}'}$ if $\F=\Q$.
By (\ref{E:period relation proof 3}), we have
\[
L(s,\itSigma \times \itPi') = L(s+\tfrac{1}{2},\itPsi \otimes \eta^{-1}\chi).
\]
Let $m \in \Z$ be the rightmost critical point of $L(s,\itSigma \times \itPi')$. The second assumption on $\underline{\ell}$ implies that the Rankin--Selberg $L$-function has more than one critical point. In particular, we have $L^{(\infty)}(m,\itSigma \times \itPi') \neq 0$.
It is clear that the pair $(\itPsi,\eta^{-1}\chi)$ satisfies condition (\ref{E:period relation proof 1}).
By Lemma \ref{L:period relation 1}-(1), Theorem \ref{T:algebraicity GSp_4 x GL_2}-(1), and Theorem \ref{T:GLH} applied to the rightmost critical value $L^{(\infty)}(m+\tfrac{1}{2},\itPsi \otimes \eta^{-1}\chi)$, we obtain the following period relation:
\begin{align}\label{E:period relation proof 4}
\begin{split}
&\sigma \left(\frac{p^{S_\infty}(\itSigma)}{|D_\F|^{1/2}\cdot(2\pi\sqrt{-1})^{\sum_{v \in S_\infty}(-7\kappa_v+4-3{\sf w})/2}\cdot(\sqrt{-1})^d\cdot G(\omega_\itPi)^3\cdot P(\itPsi)}\right)\\
& = \frac{p^{S_\infty}({}^\sigma\!\itSigma)}{|D_\F|^{1/2}\cdot(2\pi\sqrt{-1})^{\sum_{v \in S_\infty}(-7\kappa_v+4-3{\sf w})/2}\cdot(\sqrt{-1})^d\cdot G({}^\sigma\!\omega_\itPi)^3\cdot P({}^\sigma\!\itPsi)}
\end{split}
\end{align}
for all $\sigma \in {\rm Aut}(\C/\K^{\Gal})$.
On the other hand, we have
\[
L(s,\itPsi \otimes \eta^{-1}) = L(s,\itPi,{\rm Sym}^3)\cdot L(s,\itPi,{\rm Sym}^3 \otimes \omega_{\K/\F}).
\]
By Deligne's conjecture for symmetric cube $L$-function proved in \cite[Theorem 6.2]{GH1993} and \cite[Theorem 1.6]{Chen2021d}, together with Theorem \ref{T:GLH} applied to the rightmost critical point of $L(s,\itPsi\otimes\,\eta^{-1})$, we obtain the following period relation:
\begin{align}\label{E:period relation proof 5}
 \sigma \left( \frac{P(\itPsi)}{(2\pi\sqrt{-1})^{\sum_{v \in S_\infty}4\kappa_v}\cdot\Vert f_\itPi \Vert^4} \right)
 = \frac{P({}^\sigma\!\itPsi)}{(2\pi\sqrt{-1})^{\sum_{v \in S_\infty}4\kappa_v}\cdot\Vert f_{{}^\sigma\!\itPi} \Vert^4}
\end{align}
for all $\sigma \in {\rm Aut}(\C/\K^{\Gal})$.
It follows from the period relations (\ref{E:period relation proof 4}) and (\ref{E:period relation proof 5}) that \begin{align}\label{E:period relation proof 6}
\begin{split}
& \sigma \left( \frac{p^{S_\infty}(\itSigma)}{|D_\F|^{1/2}\cdot(2\pi\sqrt{-1})^{\sum_{v \in S_\infty}(\kappa_v+4-3{\sf w})/2}\cdot(\sqrt{-1})^d\cdot G(\omega_\itPi)^3\cdot\Vert f_\itPi\Vert^4}\right)\\
& = \frac{p^{S_\infty}({}^\sigma\!\itSigma)}{|D_\F|^{1/2}\cdot(2\pi\sqrt{-1})^{\sum_{v \in S_\infty}(\kappa_v+4-3{\sf w})/2}\cdot(\sqrt{-1})^d\cdot G({}^\sigma\!\omega_\itPi)^3\cdot\Vert f_{{}^\sigma\!\itPi}\Vert^4}
\end{split}
\end{align}
for all ${\rm Aut}(\C/\K^{\Gal})$. 
Let $\K'/\F$ be another CM-extension such that ${\rm BC}_{\K'}({\rm Sym}^3(\itPi))$ is cuspidal and 
\[
\K^{\Gal}\cdot\Q(\itPi)\cap (\K')^{\Gal}\cdot\Q(\itPi) = \F^{\Gal}\cdot\Q(\itPi).
\]
Then (\ref{E:period relation proof 6}) also holds for ${\rm Aut}(\C/(\K')^{\Gal})$.
In particular, we have
\begin{align}\label{E:period relation proof 7}
\frac{p^{S_\infty}(\itSigma)}{|D_\F|^{1/2}\cdot(2\pi\sqrt{-1})^{\sum_{v \in S_\infty}(\kappa_v+4-3{\sf w})/2}\cdot(\sqrt{-1})^d\cdot G(\omega_\itPi)^3\cdot\Vert f_\itPi\Vert^4} \in \F^{\Gal}\cdot\Q(\itPi).
\end{align}
We assume further that $\F^{\Gal}\cdot\Q(\itPi) \cap \K^{\Gal} = \F^{\Gal}$. Then it is easy to see that (2) of Theorem \ref{P:period relation main} follows from (\ref{E:period relation proof 6}) and (\ref{E:period relation proof 7}).

Now we prove assertion (3).  Assume $\kappa_v = \kappa_w\geq3$ for all $v,w \in S_\infty$. We also assume $\K$ is chosen so that ${\rm BC}_\K({\rm Sym}^3(\itPi))$ is cuspidal. We have the factorization of twisted exterior square $L$-function of ${\rm BC}_\K({\rm Sym}^3(\itPi))$ by $\omega_\itPi^{-3}\circ{\rm N}_{\K/\F}$:
\[
L(s,{\rm BC}_\K({\rm Sym}^3(\itPi)),\exterior{2} \otimes \,\omega_\itPi^{-3}\circ{\rm N}_{\K/\F}) = L(s,{\rm BC}_\K(\itPi),{\rm Sym}^4\otimes \omega_\itPi^{-2}\circ{\rm N}_{\K/\F})\cdot \zeta_\K(s).
\]
Thus $L(s,{\rm BC}_\K({\rm Sym}^3(\itPi)),\exterior{2} \otimes\, \omega_\itPi^{-3}\circ{\rm N}_{\K/\F})$ has a pole at $s=1$.
Let $\chi$ be a finite order Hecke character of $\A_\K^\times$ such that $\chi \neq \chi^c$. Let $m+\tfrac{1}{2}$ be the rightmost critical point of $L(s,{\rm BC}_\K({\rm Sym}^3(\itPi)))$. Note that $L(m+\tfrac{1}{2},{\rm BC}_\K({\rm Sym}^3(\itPi))) \neq 0$ by the assumption $\min_{v \in S_\infty}\{\kappa_v\} \geq 3$.
By the result of Jiang, Sun, and Tian \cite[Theorem 1.1]{JST2021}, we have
\begin{align}\label{E:period relation proof 8}
\begin{split}
\sigma \left( \frac{L^{(\infty)}(m+\tfrac{1}{2},{\rm BC}_\K({\rm Sym}^3(\itPi)) \otimes \chi)}{G(\chi\vert_{\A_\F^\times})^2\cdot L^{(\infty)}(m+\tfrac{1}{2},{\rm BC}_\K({\rm Sym}^3(\itPi)))} \right) = \frac{L^{(\infty)}(m+\tfrac{1}{2},{\rm BC}_\K({\rm Sym}^3({}^\sigma\!\itPi)) \otimes {}^\sigma\!\chi)}{G({}^\sigma\!\chi\vert_{\A_\F^\times})^2\cdot L^{(\infty)}(m+\tfrac{1}{2},{\rm BC}_\K({\rm Sym}^3({}^\sigma\!\itPi)))}
\end{split}
\end{align}
for all $\sigma \in {\rm Aut}(\C)$.
On the other hand, by Deligne's conjecture for symmetric cube $L$-function, we have
\begin{align}\label{E:period relation proof 9}
\begin{split}
 \sigma \left( \frac{L^{(\infty)}(m+\tfrac{1}{2},{\rm BC}_\K({\rm Sym}^3(\itPi)))}{(2\pi\sqrt{-1})^{4dm+2d+\sum_{v\in S_\infty}(4\kappa_v+6{\sf w})}\cdot G(\omega_\itPi)^6\cdot \Vert f_\itPi \Vert^4} \right)
 = \frac{L^{(\infty)}(m+\tfrac{1}{2},{\rm BC}_\K({\rm Sym}^3({}^\sigma\!\itPi)))}{(2\pi\sqrt{-1})^{4dm+2d+\sum_{v\in S_\infty}(4\kappa_v+6{\sf w})}\cdot G({}^\sigma\!\omega_\itPi)^6\cdot \Vert f_{{}^\sigma\!\itPi} \Vert^4}
\end{split}
\end{align}
for all $\sigma \in {\rm Aut}(\C)$.
Let $\itPi' = I_\K^\F(\chi|\mbox{ }|_{\A_\K}^{1/2})$. Then $\itPi'$ is $C$-algebraic of weight $(\underline{1};\,1)$.
We assume further that $\chi$ is chosen so that $\omega_\itPi^3\omega_{\itPi'} \neq |\mbox{ }|_{\A_\F}^{3{\sf w}+1}$ if $\F=\Q$.
In particular, the assumptions in Theorem \ref{T:algebraicity GSp_4 x GL_2} are satisfied for the pair $(\itSigma,\itPi')$ with $I=S_\infty$.
We then easily verify that assertion (3) of Theorem \ref{P:period relation main} follows from (\ref{E:period relation proof 8}), (\ref{E:period relation proof 9}), Lemma \ref{L:period relation 1}-(1), and Theorem \ref{T:algebraicity GSp_4 x GL_2}.
This completes the proof of Theorem \ref{P:period relation main}.

\end{document}